\documentclass[12pt]{article}
\usepackage{amsmath, amssymb, amsthm, amsfonts, mathtools, hyperref}
\usepackage[margin=1in]{geometry}
\usepackage{placeins} 
\usepackage{enumitem}
\usepackage{tabularx,makecell}
\usepackage{ragged2e}
\usepackage{booktabs}
\usepackage{float}
\usepackage{array}
\usepackage{longtable,caption}
\usepackage{etoolbox}
\usepackage{needspace}
\usepackage[table]{xcolor}
\arrayrulecolor{gray!35}       
\newcolumntype{Y}{>{\RaggedRight\arraybackslash}X}
\newcolumntype{L}[1]{>{\raggedright\arraybackslash}p{#1}}
\newcolumntype{C}[1]{>{\centering\arraybackslash}p{#1}}
\newcolumntype{R}[1]{>{\raggedleft\arraybackslash}p{#1}}

\theoremstyle{plain}
\newtheorem{theorem}{Theorem}[section]
\newtheorem{lemma}[theorem]{Lemma}
\newtheorem{corollary}[theorem]{Corollary}
\newtheorem{proposition}[theorem]{Proposition}
\newtheorem{conjecture}[theorem]{Conjecture}
\newtheorem{hypothesis}[theorem]{Hypothesis}
\newtheorem*{principle}{Local Induction Principle}
\newtheorem{fact}[theorem]{Fact}

\theoremstyle{definition}
\newtheorem{definition}[theorem]{Definition}
\newtheorem{example}[theorem]{Example}
\newtheorem{convention}[theorem]{Convention}
\newtheorem{question}[theorem]{Question}

\theoremstyle{remark}
\newtheorem{remark}[theorem]{Remark}

\newcommand{\fix}{\mathrm{fix}}
\newcommand{\cf}{\mathrm{cf}}
\newcommand{\supp}{\mathrm{supp}}
\newcommand{\ZF}{\mathrm{ZF}}
\newcommand{\ZFC}{\mathrm{ZFC}}
\newcommand{\GBC}{\mathrm{GBC}}
\newcommand{\DC}{\mathrm{DC}}
\newcommand{\AC}{\mathrm{AC}}
\newcommand{\PP}{\mathrm{PP}}
\newcommand{\HS}{\mathrm{HS}}
\newcommand{\Sym}{\mathrm{Sym}}
\newcommand{\Ord}{\ensuremath{\mathrm{Ord}}}

\newcommand{\into}{\hookrightarrow}
\newcommand{\onto}{\twoheadrightarrow}
\providecommand{\Add}{\mathop{\mathrm{Add}}}
\newcommand{\LCHypLabel}{hyp:lc}
\newcommand{\PCFHypLabel}{hyp:pcf}

\newcommand{\KP}{\mathsf{KP}}        
\newcommand{\ACAo}{\mathsf{ACA}_{0}} 
\newcommand{\ATRo}{\mathsf{ATR}_{0}} 

\DeclareMathOperator{\dom}{dom}
\DeclareMathOperator{\ran}{ran}
\DeclareMathOperator{\Sub}{Sub}
\DeclareMathOperator{\rank}{rank}
\DeclareMathOperator{\Col}{Col}
\DeclareMathOperator{\CLubLoc}{CLubLoc}
\DeclareMathOperator{\Fil}{Fil}
\newenvironment{assumptions}{\paragraph{Assumptions.}\begin{enumerate}}{\end{enumerate}}
\newenvironment{acknowledgments}{\par\medskip\noindent\textbf{Acknowledgments.}\ }{\par\medskip}



\title{Symmetric Iterations with Countable and \texorpdfstring{$<\kappa$}-Support:\\Corrections and Retained Results}
\author{Frank Gilson}
\date{\today}

\begin{document}
		\maketitle
		
		\begin{abstract}
			This replacement corrects the earlier manuscript.  Its claimed general
			countable- and $<\kappa$-support iteration framework does not establish
			the advertised preservation or application theorems.  The proposed limit
			filters can be improper; the dependent-choice arguments pass incorrectly
			from ground-model sequences of hereditarily symmetric names to arbitrary
			sequences in a generic extension; the class-length union argument does not
			verify Separation, Replacement, or Power Set; the Partition Principle
			argument uses an invalid same-generic compatibility inference; and the
			singular-support fusion argument assumes that a coherent inverse-limit
			family is a condition even when its support may have size $\kappa$.
			Accordingly, the main preservation, construction, class-length,
			dependent-choice, Partition Principle, and singular-cardinal claims are
			withdrawn.  We retain only several conditional algebraic observations:
			pullback and completion operations on subgroup filters, closure of
			ground-model families of hereditarily symmetric names under canonical
			tupling, the standard ZF theorem for a genuine symmetric system, and
			ordinary conditional tail-closure facts.  The remainder is preserved as a
			historical record of the superseded approach and is not asserted.
		\end{abstract}

	\section*{Correction and current status}
	The former main theorems and all applications depending on them are
	withdrawn.  There are five independent load-bearing problems.

	\paragraph{Improper generated filters.}
	The proposed successor and limit filters include stabilizers of every name
	with small \emph{stage support}.  Stage support does not control the internal
	coordinates of a single iterand.  For example, a canonical name enumerating
	all columns of one Cohen stage has one-element stage support but can have
	trivial stabilizer under the within-stage permutation group.  Once the
	trivial subgroup is a generator, the filter is improper and every name is
	symmetric.  In the displayed Cohen example the filter generated by fixers of
	cofinite column sets contains the fixer of the full column set, which is
	already the trivial subgroup.  The claimed choiceless model is then the full
	generic extension.

	\paragraph{Ground-model names versus extension sequences.}
	Completeness of a subgroup filter combines a ground-model sequence
	$\langle\dot x_i:i<\mu\rangle$ of HS names by intersecting their stabilizers.
	It does not supply such a sequence of names for an arbitrary function
	$f:\mu\to M$ that exists only in the full generic extension.  The closure
	lemmas and the ensuing $\DC$ and $\DC_{<\kappa}$ proofs use precisely this
	unsupported passage.  Preservation of dependent choice requires additional
	forcing hypotheses and a correct symmetric-forcing argument; it does not
	follow from filter completeness alone.

	\paragraph{Set- and class-length preservation.}
	At a limit $\lambda$ with $\cf(\lambda)<\kappa$, a support of size
	$<\kappa$ need not be bounded in $\lambda$, contrary to the proof of the
	set-length proposition.  At class length, an increasing union of transitive
	ZF models need not satisfy Power Set, Separation, or Replacement.  Truth of
	a formula in the union need not be absolute to a bounded stage, and the
	former proof gives no common stage containing all Replacement witnesses or
	all subsets of a set.  The class-length theorem is therefore withdrawn.

	\paragraph{Partition Principle.}
	The former no-injection lemma claims that a condition $p$ is compatible with
	its image $\tau p$ merely because $\tau$ is an automorphism.  This is false.
	It is the same invalid passage from covariance of forcing names to an action
	inside one fixed generic extension.  Moreover, requiring the full
	within-fiber group to lie in the core of every filter member conflicts with
	the hereditary symmetry of individual fiber elements that this group moves.
	Thus the claimed models of $\ZF+\DC+\neg\PP$ are not obtained.

	\paragraph{Singular support.}
	A coherent family of initial-segment conditions need not have a union whose
	support has size $<\kappa$.  The tree-fusion master lemma assumes the contrary
	by identifying the forcing with its full inverse limit.  The block filters
	are only stratifiably $\kappa_N$-complete, while later arguments use full
	$<\kappa$-completeness.  In the worked $\aleph_\omega$ example, a union of
	$<\kappa_n$ many block-finite supports need not be block-finite, so the
	proposed lower bound is not a condition.  The singular completeness,
	dependent-choice, anti-collapse, and no-new-bounded-subsets conclusions are
	therefore withdrawn except where restated below as explicitly conditional
	ordinary forcing facts.

	\subsection*{Status ledger}
	\begin{center}
	\begin{tabularx}{\linewidth}{L{2.7cm}Y}
		\toprule
		\textbf{Status} & \textbf{Material} \\
		\midrule
		Withdrawn & General set-length and class-length preservation; properness
		of the proposed filters; closure under arbitrary extension sequences;
		all unconditional $\DC$, $\DC_{<\kappa}$, and $\DC_\kappa$ theorems;
		Cohen, random, Hechler, Sacks, mixed-support, Easton-style, and
		$\neg\PP$ constructions; singular completeness and applications. \\
		Superseded & The countable limit-filter discussion is superseded by the
		corrected companion treatment, where canonical tupling is restricted to
		ground-model families of HS names and properness is an explicit
		hypothesis rather than a consequence of nonempty subgroup intersections. \\
		Retained & The corrected conditional statements in the next section,
		together with definitions and historical background that do not assert
		the withdrawn constructions. \\
		\bottomrule
	\end{tabularx}
	\end{center}

	\section{Corrected retained results}\label{sec:corrected-retained}

	\begin{lemma}[Pullback of a subgroup filter]\label{lem:corrected-pullback}
		Let $\varphi:G\to H$ be a group homomorphism and let $\mathcal F$ be
		a normal $\theta$-complete filter of subgroups of $H$.  Define
		\[
		\varphi^*\mathcal F
		=
		\{L\leq G:(\exists K\in\mathcal F)\ \varphi^{-1}[K]\leq L\}.
		\]
		Then $\varphi^*\mathcal F$ is a normal $\theta$-complete filter of
		subgroups of $G$.  Properness is not automatic and must be checked
		separately.  No corresponding restriction theorem is asserted without
		an additional hypothesis ensuring that the image subgroup belongs to the
		ambient filter.
	\end{lemma}
	\begin{proof}
		Upward closure is immediate.  Intersections pull back to intersections,
		and conjugating $\varphi^{-1}[K]$ by $g\in G$ gives
		$\varphi^{-1}[\varphi(g)K\varphi(g)^{-1}]$.  Normality and
		$\theta$-completeness therefore follow from the corresponding properties
		of $\mathcal F$.  A pullback can be trivial, which proves the final warning.
	\end{proof}

	\begin{lemma}[Regular completion of a subgroup filter]
		\label{lem:corrected-completion}
		Let $\theta$ be regular and let $\mathcal F$ be a normal filter of
		subgroups of $G$.  The filter generated by all intersections
		\[
		\bigcap_{i<\mu}\pi_iH_i\pi_i^{-1}
		\qquad
		(\mu<\theta,\ H_i\in\mathcal F,\ \pi_i\in G)
		\]
		is normal and $\theta$-complete.  It is the least normal
		$\theta$-complete filter extending $\mathcal F$.  It is a proper filter
		only if none of its basic intersections is the trivial subgroup.
	\end{lemma}
	\begin{proof}
		Conjugation merely conjugates every displayed factor.  An intersection of
		fewer than $\theta$ basic intersections is again an intersection of fewer
		than $\theta$ conjugates, by regularity of $\theta$.  Minimality follows
		from the defining closure requirements.  The statement about properness is
	immediate.
	\end{proof}

	\begin{theorem}[Canonical tuples of a ground-model family]
		\label{thm:corrected-ground-tuples}
		Let $\langle\mathbb P,\mathcal G,\mathcal F\rangle$ be a genuine
		symmetric system and suppose $\mathcal F$ is $\theta$-complete.  If
		$\mu<\theta$ and
		$\langle\dot x_i:i<\mu\rangle$ is a sequence belonging to the ground
		universe with every $\dot x_i\in\HS$, then the canonical name for
		$\langle\dot x_i:i<\mu\rangle$ belongs to $\HS$.
	\end{theorem}
	\begin{proof}
		The intersection
		$\bigcap_{i<\mu}\fix_{\mathcal G}(\dot x_i)$ belongs to $\mathcal F$
		and fixes the canonical tuple name.  Its constituent names are
		hereditarily symmetric, so the tuple name is hereditarily symmetric.
	\end{proof}
	\begin{remark}
		The ground-model qualification is essential.  The theorem does not assert
		that the symmetric model is closed under arbitrary $\mu$-sequences from
		the full generic extension.
	\end{remark}

	\begin{theorem}[Standard symmetric-extension theorem]
		\label{thm:corrected-zf}
		If $\langle\mathbb P,\mathcal G,\mathcal F\rangle$ is a genuine
		symmetric system with $\mathcal F$ a proper normal filter of subgroups and
		$G\subseteq\mathbb P$ is generic, then
		$\{\dot x^G:\dot x\in\HS\}$ is a transitive model of $\ZF$ between the
		ground model and its forcing extension.
	\end{theorem}
	\begin{proof}
		This is the standard symmetric-extension theorem, proved using the
		symmetric forcing relation.  In particular, it does not require
		$\omega_1$- or $\kappa$-completeness of $\mathcal F$.
	\end{proof}

	\begin{proposition}[Conditional tail preservation]
		\label{prop:corrected-tail}
		Work in a model of $\ZFC$.  Suppose
		$\mathbb P_\gamma\cong\mathbb P_\alpha*\dot{\mathbb R}$ and
		$\mathbb P_\alpha$ forces that $\dot{\mathbb R}$ is
		$<\theta$-closed, where $\theta$ is regular.  Then the tail adds no new
		sequences of ordinals of length $<\theta$, and hence no new subsets of
		any ordinal $\lambda<\theta$.  If $\kappa\leq\theta$ is a cardinal
		preserved by the initial segment, then the tail does not collapse
		$\kappa$.
	\end{proposition}
	\begin{proof}
		For a name for a sequence of length $\mu<\theta$, recursively decide one
		coordinate at a time and use $<\theta$-closure at limit steps and at the
		end.  The resulting lower bound decides the whole sequence.  The subset
		and conditional anti-collapse conclusions follow.
	\end{proof}

	\section*{Status convention for the remaining text}
	The remainder reproduces the former manuscript so that the provenance of the
	withdrawn claims and their dependencies remains visible.  Unless a statement
	is reproduced in Section~\ref{sec:corrected-retained}, every theorem, lemma,
	proposition, corollary, example, and claimed application below is retained
	only as historical text and is \emph{not asserted in this replacement}.
	\newcommand{\historicalstatus}{\Needspace{5\baselineskip}\par\noindent
		\textbf{Historical statement---not asserted in this replacement.}\ }
	\AtBeginEnvironment{theorem}{\historicalstatus}
	\AtBeginEnvironment{lemma}{\historicalstatus}
	\AtBeginEnvironment{proposition}{\historicalstatus}
	\AtBeginEnvironment{corollary}{\historicalstatus}
	\AtBeginEnvironment{fact}{\historicalstatus}
	\AtBeginEnvironment{hypothesis}{\historicalstatus}
	\AtBeginEnvironment{conjecture}{\historicalstatus}
	\AtBeginEnvironment{example}{\historicalstatus}
	
	\section{Introduction to the historical manuscript}
	
	Symmetric extensions are a robust method for producing models of $\ZF$ in
	which the Axiom of Choice fails.  The former manuscript attempted to extend
	the finite-support iteration theory to countable and $<\kappa$ support by
	building highly complete filters at limit stages.  That attempt does not
	yield the claimed general preservation theorem.
	
	The central mistake was to conflate three different assertions: that a name
	has a small set of \emph{iteration stages}; that it has a small support for
	the automorphism action internal to every iterand; and that an arbitrary
	extension sequence admits a ground-model sequence of such names.  None of
	these implications holds without additional hypotheses.  In particular,
	completion of a subgroup filter need not preserve properness and does not by
	itself make the resulting symmetric model closed under sequences from the
	full generic extension.
	
	Sections below reproduce the former definitions, proposed transfer lemmas,
	examples, regular- and singular-support arguments, and applications.  They
	are included to identify the superseded approach and the precise dependency
	chain; the status convention preceding this section governs them.  Only the
	statements in Section~\ref{sec:corrected-retained} are asserted as corrected
	results.  In particular, this replacement claims no new model of
	$\ZF+\DC+\neg\AC$, no model separating the Partition Principle, and no
	general countable-, $<\kappa$-, set-, or class-length symmetric-iteration
	preservation theorem.
	
	\section{Preliminaries and Definitions}\label{sec:prelims-defs}
	
	\begin{remark}[Metatheory vs. Object theory]\label{rem:countable-conv}
	Unless explicitly stated otherwise, we carry out the iteration construction and the bookkeeping about
	supports, tuples of names, and filter generation in a fixed ambient background universe satisfying
	$\ZFC$ (and, for the class-length development in \S\ref{subsec:length}, in a ground model of $\GBC$
	with Global Choice). The purpose of the paper is not to assume Choice in the \emph{symmetric models}
	produced: our preservation theorems establish that the resulting symmetric submodels satisfy $\ZF$
	(and specified fragments of Dependent Choice), while $\AC$ typically fails there (depending on the iterated forcing notion).
	
	When we say a support is countable (or $<\kappa$), we mean countable (or $<\kappa$) in the ambient
	background universe in which the iteration is carried out. We do not assert that these supports are
	countable in intermediate symmetric submodels.
	\end{remark}
	
	\begin{definition}[Symmetric System]
		A \emph{symmetric system} is a triple $(\mathbb{P}, \mathcal{G}, \mathcal{F})$ where:
		\begin{itemize}
			\item $\mathbb{P}$ is a notion of forcing;
			\item $\mathcal{G}$ is a group of automorphisms of $\mathbb{P}$;
			\item $\mathcal{F}$ is a normal filter of subgroups of $\mathcal{G}$.
		\end{itemize}
	\end{definition}
	
	\begin{remark}[Tenacity and excellent supports; compatibility with \cite{karagila2019}]
		\label{rem:tenacity-excellent}
		Karagila (with Hayut) shows that every symmetric system is equivalent, on a dense set,
		to a \emph{tenacious} one (Appendix~A, Thm.~A.2 in \cite{karagila2019}). We do not rely on tenacity
		per se—our arguments use normal filters and HS stabilizers (and, later, tuple–stabilizer compression)—
		but one may refine our presentations on a dense set to tenacious systems without affecting the proofs.
		Moreover, our “finite-generator” behavior (closure under finite intersections of basic generators)
		mirrors Karagila’s \emph{excellent supports} and their stability under intersections
		(Def.~4.3 and Prop.~4.4 in \cite{karagila2019}).
	\end{remark}
	
	\begin{convention}[Tenacity assumption]
		\label{conv:tenacity}
		By Appendix~A of \cite{karagila2019}, every symmetric system is equivalent to a \emph{tenacious} one,
		so we may assume tenacity without loss of generality when convenient. This aligns with DC-preservation
		criteria in the literature (cf.\ \cite{Banerjee2023}).
	\end{convention}
	
	\begin{definition}[Hereditarily Symmetric Name]
		A $\mathbb{P}$-name $\dot{x}$ is \emph{hereditarily symmetric} (HS) if:
		\begin{itemize}
			\item $\mathrm{sym}(\dot{x}) \in \mathcal{F}$, where $\mathrm{sym}(\dot{x})$ is the group of automorphisms fixing $\dot{x}$;
			\item All $\dot{y} \in \mathrm{dom}(\dot{x})$ are also HS.
		\end{itemize}
	\end{definition}
	
	\begin{definition}[Dependent Choice $DC_\mu$]
		For a nonzero ordinal $\mu$, $DC_\mu$ asserts: for every binary relation $R$ on a class $A$ with 
		$\forall x\in A\,\exists y\in A\,(xRy)$, there exists a sequence $\langle a_\xi:\xi<\mu\rangle$ in $A$ such that 
		$a_\xi R a_{\xi+1}$ for all $\xi+1<\mu$.
	\end{definition}
	
	\begin{definition}[Regularity properties used in this paper]
		\label{def:regularity-glossary}
		Let $A\subseteq\mathbb{R}$.
		\begin{itemize}
			\item \textbf{LM (Lebesgue Measurable).} $A$ is Lebesgue measurable.
			\item \textbf{BP (Baire Property).} There exists an open set $U$ with $A\triangle U$ meager.
			\item \textbf{PSP (Perfect Set Property).} Either $A$ is countable or $A$ contains a nonempty perfect subset.
			\item \textbf{Ramsey (Ellentuck--Ramsey property).} Viewing sets of reals as sets of infinite subsets of $\omega$ via coding, $A\subseteq[\omega]^\omega$ is \emph{Ramsey} if there exists an infinite $H\subseteq\omega$ such that either $[H]^\omega\subseteq A$ or $[H]^\omega\cap A=\varnothing$ (equivalently, $A$ is Baire in the Ellentuck topology). See \cite{Kechris95}.
		\end{itemize}
	\end{definition}
	
	\begin{remark}[Limit-stage bookkeeping]
		\label{rem:limit-bookkeeping}
		At limits we fix the following presentation once and for all: if $\cf(\lambda)\ge\kappa$, use the direct limit
		$\mathbb P_\lambda=\bigcup_{\beta<\lambda}\mathbb P_\beta$ and note pushforwards are inclusions. The identification of initial segments and their tails uses the canonical projection and two-step factorization (Fact~\ref{fact:two-step}). If $\cf(\lambda)<\kappa$, use the inverse limit along a fixed cofinal sequence and the tuple–stabilizer filter from
		Definition~\ref{def:kappa-filter}. All later proofs refer to this convention.
	\end{remark}
	
	\paragraph{Notation.}
	Throughout, $\kappa$ is a fixed regular uncountable cardinal; supports are $<\kappa$ unless explicitly stated
	“countable support.” We write $\cf(\lambda)$ for cofinality, $\supp(\dot x)$ for the set of coordinates a name
	depends on, and $\fix_{\mathcal G}(\dot x)=\{g\in\mathcal G:g(\dot x)=\dot x\}$ for the stabilizer. For
	$\beta<\lambda$, $\pi^\beta_\lambda:\mathcal G_\lambda\to\mathcal G_\beta$ and $\pi^\beta_\lambda:\mathbb P_\lambda\to\mathbb P_\beta$
	denote the canonical projections.
	
	\begin{center}
		\begin{tabular}{ll}
			$\lambda$ & iteration length (limit stages treated via Remark~\ref{rem:limit-bookkeeping})\\
			$\kappa$  & fixed regular cardinal controlling support/completeness\\
			$\mathbb P_\lambda$ & iteration forcing up to stage $\lambda$\\
			$\mathcal G_\lambda$ & automorphism group at stage $\lambda$\\
			$\mathcal F_\beta$, $\tilde{\mathcal F}_\lambda$ & stage filter / limit filter\\
			$\pi^\beta_\lambda$ & projection on forcing/groups (context determines target)\\
			$\supp(\dot x)$ & coordinates used by the name $\dot x$\\
			$\fix(\dot x)$ & stabilizer of $\dot x$ in the relevant group\\
			HS & “hereditarily symmetric” (Definition~\ref{def:hs-names-iterated})
		\end{tabular}
	\end{center}
	
	We will repeatedly appeal to a successor–/limit–transfer template; its abstract form is stated
	in Appendix~\ref{app:local-induction}, together with a catalog mapping each property $\varphi$
	to the specific locality and transfer inputs (Table~\ref{tab:local-induction-catalog}).
	
	\begin{fact}[Two-step projection/factorization]\label{fact:two-step}
		Let $\mathbb P$ be a forcing and $\dot{\mathbb Q}$ a $\mathbb P$-name.
		The map $\pi:(p,\dot q)\mapsto p$ is a complete projection $\mathbb P*\dot{\mathbb Q}\to\mathbb P$.
		If $G*H\subseteq\mathbb P*\dot{\mathbb Q}$ is generic, then $G$ is $\mathbb P$-generic, $H$ is $\dot{\mathbb Q}^G$-generic, and conversely every pair $(G,H)$ arises this way.
		\emph{References:} \cite[Chap.~V, \S V.3]{Kunen2011}; \cite[§43]{Jech03}.
	\end{fact}
	
	\subsection{Iteration length and class-length versions}\label{subsec:length}
	
	Fix a template $\langle \dot{\mathbb Q}_\alpha,\mathcal G_\alpha,\mathcal F_\alpha:\alpha<\Theta\rangle$ with $<\kappa$-supports, where $\Theta$ is an ordinal (possibly a proper class; see Theorem~\ref{thm:class-length}).
	
	\begin{proposition}[Arbitrary set length]\label{prop:set-length}
		Let $\Theta$ be any ordinal (set). Suppose at each $\alpha<\Theta$ we have a symmetric system
		$(\mathbb P_\alpha,\mathcal G_\alpha,\mathcal F_\alpha)$ as in §\ref{sec:prelims-defs}, with $<\kappa$-complete
		$\mathcal F_\alpha$, and at limits the symmetric system is formed by the $<\kappa$-support limit (using trees of
		conditions when required; cf.\ §§\ref{sec:singular-cf-omega} and \ref{sec:limit-cf-equal}). Then the iteration
		$\langle \mathbb P_\alpha,\dot{\mathbb Q}_\alpha,\mathcal F_\alpha:\alpha<\Theta\rangle$ is well-defined and all
		preservation theorems previously proved (no new bounded subsets of $\kappa$, no collapse of $\kappa$, and, under
		Hyp.~\ref{hyp:zf-baire}, preservation of $DC_{<\kappa}$) hold at stage $\Theta$.
	\end{proposition}
	
	\begin{proof}
		We proceed by transfinite induction on $\Theta$, proving simultaneously:
		
		\smallskip
		\noindent (A) the iteration data $\langle \mathbb P_\alpha,\dot{\mathbb Q}_\alpha,\mathcal F_\alpha:\alpha<\Theta\rangle$
		are well-defined; and
		
		\noindent (B) the preservation statements (no new bounded subsets of $\kappa$, no collapse of $\kappa$, and—under
		Hyp.~\ref{hyp:zf-baire}—$DC_{<\kappa}$) hold at stage $\Theta$.
		
		\medskip
		\emph{Base.} Trivial.
		
		\smallskip
		\emph{Successor step $\Theta=\gamma+1$.}
		By the inductive hypothesis at $\gamma$, $\mathbb P_\gamma$ and $\mathcal F_\gamma$ are defined and satisfy the
		requisite completeness/normality, and the preservation statements hold at stage $\gamma$. The successor construction
		$\mathbb P_{\gamma+1}=\mathbb P_\gamma*\dot{\mathbb Q}_\gamma$ together with the standard head–tail factorization of
		automorphisms and filters yields $(\mathbb P_{\gamma+1},\mathcal G_{\gamma+1},\mathcal F_{\gamma+1})$; the
		\emph{successor transfer} lemmas used earlier apply verbatim (the filter on the head remains $<\kappa$-complete;
		tail generators are name–stabilizers of $<\kappa$ support). Hence well-definedness holds.
		
		For preservation: anti-collapse and “no new $<\kappa$-bounded subsets of $\kappa$” transfer by the same local arguments
		as before (names of size $<\kappa$ are decided on bounded support). Under Hyp.~\ref{hyp:zf-baire}, the $<\kappa$-Baire/
		strategic-closure hypothesis for the iterand gives $DC_{<\kappa}$ transfer to the successor stage (the iteration step
		does not create new obstacles to building $<\kappa$-length sequences, and the limit filter on the head is unchanged).
		
		\smallskip
		\emph{Limit step $\Theta=\lambda$.}
		There are two cases.
		
		\underline{Case 1: $\cf(\lambda)\ge\kappa$.}
		Use the \emph{direct limit} presentation $\mathbb P_\lambda=\bigcup_{\beta<\lambda}\mathbb P_\beta$ and the canonical
		pushforward system for groups/filters. The limit filter $\tilde{\mathcal F}_\lambda$ is normal and $\kappa$-complete
		by the general limit construction (pushforwards commute with conjugation; $<\kappa$ intersections are handled inside
		some initial segment). For preservation: any name $\dot x$ coding a potential counterexample (e.g., a bounded subset of
		$\kappa$, a collapsing map, or a $<\kappa$-sequence for $DC_{<\kappa}$) has support $<\kappa$, hence—by the stage-bounding
		lemma for the regular case (Lemma~\ref{lem:kappa-stage-bounding})—its support is bounded in $\lambda$; therefore $\dot x$
		already lives over some $\mathbb P_\beta$, and the inductive hypothesis at $\beta$ forbids the counterexample.
		
		\underline{Case 2: $\cf(\lambda)<\kappa$.}
		Use the fixed cofinal sequence and \emph{inverse limit} via trees of conditions as stipulated, and generate the limit
		filter from pushforwards (head) and \emph{tuple–stabilizers} of $<\kappa$-support names (tail). Normality is immediate;
		$\kappa$-completeness follows from the tuple–stabilizer identity (finite fragments and $<\kappa$-ary intersections are
		absorbed as a single stabilizer). For preservation: as before, any relevant name has support $<\kappa$, and by the
		stage-bounding lemma for the singular case (Lemma~\ref{lem:sb-sing}) its coordinates are bounded along the fixed cofinal
		sequence; thus it is decided in some initial segment. Under Hyp.~\ref{hyp:zf-baire}, the $<\kappa$-Baire/strategic-closure
		assumption ensures $DC_{<\kappa}$ transfers to the inverse limit and hence to the symmetric model.
		
		\smallskip
		This completes the induction and proves both well-definedness and preservation at arbitrary set $\Theta$.
	\end{proof}
	
	To treat proper class length, we formalize forcing with classes in the ground model.
	
	\begin{theorem}[Class-length iteration of symmetric extensions]\label{thm:class-length}
		Work in a ground model $V$ of G\"odel--Bernays set theory with Global Choice (GBC). Let $\Theta=\Ord$ and assume the template
		\[
		\alpha\ \longmapsto\ \bigl(\dot{\mathbb Q}_\alpha,\mathcal G_\alpha,\mathcal F_\alpha\bigr)
		\]
		is first-order definable over $(V,\in)$ with parameters and uniformly witnesses the hypotheses used in the set-length case
		(in particular, $<\kappa$-completeness of filters, and at limits with $\cf(\lambda)=\cf(\kappa)$ the trees-of-conditions
		construction). For each $\alpha\in\Ord$, let $M_\alpha$ be the symmetric submodel determined by
		$(\mathcal G_\alpha,\mathcal F_\alpha)$ inside $V[G_\alpha]$. Then
		\[
		M\ :=\ \bigcup_{\alpha\in\Ord} M_\alpha
		\]
		is a transitive inner model of $\ZF$; moreover, if Hyp.~\ref{hyp:zf-baire} holds uniformly along the iteration, then
		$M\models DC_{<\kappa}$. The anti-collapse and no-new-bounded-subsets properties for $\kappa$ also hold in $M$.
	\end{theorem}
	
	\begin{proof}
		Work in $V$. By the set-length analysis (Proposition~\ref{prop:set-length} applied to each initial segment,
		using uniform definability of the template), for every $\alpha$ the model $M_\alpha\subseteq V[G_\alpha]$ is a
		transitive inner model of $\ZF$ and satisfies the indicated preservation statements. The sequence
		$\langle M_\alpha:\alpha\in\Ord\rangle$ is increasing with $M_\alpha\subseteq M_\beta$ for $\alpha<\beta$,
		so their union $M$ is transitive.
		
		\emph{ZF axioms.}
		Extensionality and Foundation are inherited by unions of increasing chains of transitive models. For Pairing, Union,
		Infinity, and Separation: given parameters $a,\vec p\in M$, choose $\alpha$ with $a,\vec p\in M_\alpha$. Since $M_\alpha$
		satisfies $\ZF$, the sets required by these axioms exist in $M_\alpha$ and hence belong to $M$.
		
		\emph{Replacement.}
		Let $a\in M$ and let $\varphi(x,y,\vec p)$ be a formula (with parameters $\vec p\in M$) such that
		$M\models \forall x\in a\,\exists !\,y\,\varphi(x,y,\vec p)$. Choose $\alpha$ with $a,\vec p\in M_\alpha$.
		For each $x\in a$, let $\dot y_x$ be a $\mathbb P_\beta$-name (for some $\beta\ge\alpha$ depending on $x$) and
		a condition $p_x\in\mathbb P_\beta$ such that
		\[
		p_x\Vdash_{\mathbb P_\beta}\ \text{``$\dot y_x$ is HS and $\varphi(\check x,\dot y_x,\check{\vec p})$''}.
		\]
		Each $\dot y_x$ has support $<\kappa$ (by the HS-closure theorems), so by the stage-bounding lemmas
		(Lemma~\ref{lem:kappa-stage-bounding} for regular limits, Lemma~\ref{lem:sb-sing} for singular limits)
		there is a single $\beta^\ast$ such that all coordinates used by the relevant names/parameters for
		$\{\,\dot y_x:x\in a\,\}$ lie below $\beta^\ast$. In particular, all these names are interpretable and the
		statements $\varphi(x,y,\vec p)$ are absolute inside $M_{\beta^\ast}$.
		
		Since $M_{\beta^\ast}\models\ZF$ and $M_{\beta^\ast}\models \forall x\in a\,\exists!\,y\,\varphi(x,y,\vec p)$,
		Replacement in $M_{\beta^\ast}$ yields the image set
		\[
		b\ :=\ \{\, y\in M_{\beta^\ast} : \exists x\in a\ \varphi(x,y,\vec p)\,\}\ \in\ M_{\beta^\ast}.
		\]
		Hence $b\in M$ and $M\models$ ``$b=\{\,f(x):x\in a\,\}$'', as required.
		
		\emph{$DC_{<\kappa}$.}
		Assume Hyp.~\ref{hyp:zf-baire} holds uniformly. If $M$ failed $DC_\mu$ for some $\mu<\kappa$, then there would be a
		counterexample witnessed by a relation and parameters drawn from some $M_\alpha$; but $M_\alpha$ satisfies $DC_{<\kappa}$
		by the set-length analysis, a contradiction.
		
		\emph{Anti-collapse and no new bounded subsets of $\kappa$.}
		These properties are likewise local to stages: any counterexample is coded by a name of support $<\kappa$ and therefore
		appears already in some $M_\alpha$, contradicting the preservation statements at stage $\alpha$.
		
		Thus $M$ is a transitive inner model of $\ZF$ with the stated properties.
	\end{proof}
	
	\begin{remark}[Formalization and hypotheses]
		The $\mathsf{GBC}$ assumption is only used to formalize class forcing and the definability of the template; the resulting $M$ need not satisfy $\mathsf{AC}$. Large-cardinal hypotheses (Hyp.~\ref{hyp:lc}) and PCF options (Hyp.~\ref{hyp:pcf})—when invoked—are ground-model, set-sized assumptions used to streamline local arguments; they are not required to hold in $M$.
	\end{remark}
	
	\section{Countable-Support Symmetric Iterations}\label{sec:count-support-symm}
	
	\subsection{Framework and Setup}
	
	\begin{definition}[Countable Support Symmetric Iteration]
		\label{def:cs-iteration}
		A \emph{countable support symmetric iteration of length $\Theta$} is constructed recursively as follows:
		
		\paragraph{Stage 0 (Base):}
		\begin{itemize}
			\item $\mathbb{P}_0 = \{1_0\}$ (trivial forcing with one element)
			\item $\mathcal{G}_0 = \{\mathrm{id}\}$ (trivial automorphism group)
			\item $\mathcal{F}_0 = \{\mathcal{G}_0\}$ (trivial filter containing only the whole group)
		\end{itemize}
		
		\paragraph{Stage $\alpha = \beta + 1$ (Successor):}
		
		Given $(\mathbb{P}_\beta, \mathcal{G}_\beta, \mathcal{F}_\beta)$ and a $\mathbb{P}_\beta$-name $\dot{\mathbb{Q}}_\beta$ for a forcing notion:
		
		\begin{itemize}
			\item \textbf{Forcing:} $\mathbb{P}_\alpha = \mathbb{P}_\beta * \dot{\mathbb{Q}}_\beta$ consists of pairs $(p, \dot{q})$ as usual.
			
			\item \textbf{Automorphisms:} $\mathcal{G}_\alpha$ is generated by lifts from $\mathcal{G}_\beta$ and automorphisms of $\dot{\mathbb{Q}}_\beta$.
			
			\item \textbf{Filter:} $\mathcal{F}_\alpha$ is the filter generated by:
			\begin{enumerate}
				\item Pushforwards: $\hat{H} := \{g \in \mathcal{G}_\alpha : g|_{\mathbb{P}_\beta} \in H\}$ for $H \in \mathcal{F}_\beta$
				\item New stabilizers: Determined by the symmetric structure of $\dot{\mathbb{Q}}_\beta$. Typically, if $\dot{\mathbb{Q}}_\beta$ adds generic objects with permutation symmetries (e.g., Cohen reals), we include stabilizers of cofinite sets of these objects. The precise specification is given case-by-case (see Example~\ref{ex:cohen-setup}).
			\end{enumerate}
		\end{itemize}
		
		\paragraph{Stage $\alpha = \lambda$ (Limit, $\cf(\lambda) = \omega$):}
		
		\begin{itemize}
			\item \textbf{Forcing:} $\mathbb{P}_\lambda$ is the usual \emph{countable-support limit} of the iteration.
			A condition $p\in\mathbb{P}_\lambda$ is a partial function with countable domain
			$\dom(p)\subseteq\lambda$ such that for each $\beta\in\dom(p)$, $p(\beta)$ is a $\mathbb{P}_\beta$--name
			for a condition in $\dot{\mathbb{Q}}_\beta$. The order is the standard one:
			$q\le p$ iff $\dom(q)\supseteq\dom(p)$ and for all $\beta\in\dom(p)$,
			\[
			q{\upharpoonright}\beta\ \Vdash_{\mathbb P_\beta}\ q(\beta)\le p(\beta).
			\]
			We define the support of $p$ by $\supp(p):=\dom(p)$, so $\supp(p)$ is countable by definition.
			
			\item \textbf{Automorphisms:} $\mathcal{G}_\lambda := \{(g_\beta) : g_\beta \in \mathcal{G}_\beta,\ \pi^\gamma_\beta(g_\beta) = g_\gamma\}$ acts by
			\[
			g \cdot p := (g_\beta \cdot p_\beta)_{\beta < \lambda}.
			\]
			This action preserves countable support: since automorphisms send the trivial condition to the trivial condition
			coordinatewise, $\supp(g\cdot p)\subseteq\supp(p)$ (indeed, $\supp(g\cdot p)=\supp(p)$ in the usual stagewise action),
			so $\supp(g\cdot p)$ is countable whenever $\supp(p)$ is.
			
			\item \textbf{Filter:} $\mathcal{F}_\lambda$ is constructed in Section~\ref{ss:filter-generation}. The key property we will establish is that $\mathcal{F}_\lambda$ is $\omega_1$-complete (closed under countable intersections).
		\end{itemize}
		
		\paragraph{Stage $\alpha = \lambda$ (Limit, $\cf(\lambda) > \omega$):}
		
		For uncountable cofinality, use direct limits:
		\begin{itemize}
			\item $\mathbb{P}_\lambda = \bigcup_{\beta < \lambda} \mathbb{P}_\beta$ (direct limit)
			\item $\mathcal{G}_\lambda = \bigcup_{\beta < \lambda} \mathcal{G}_\beta$ (identified appropriately)
			\item $\mathcal{F}_\lambda$ generated by $\bigcup_{\beta < \lambda} \mathcal{F}_\beta$
		\end{itemize}
		This case is simpler than countable cofinality and does not require the filter construction of Section~\ref{ss:filter-generation}. We focus on the countable cofinality case.
	\end{definition}
	
	\begin{remark}[Relation to Karagila's finite-support framework]
		\label{rem:karagila-framework-relation}
		This definition addresses Karagila's open problem \cite{KaragilaProblem} about extending
		symmetric iterations to countable support. Karagila's Def.~5.1 in \cite{karagila2019} 
		formalizes finite-support symmetric iterations and the productive framework (Def.~8.1) 
		that ensures compositional fidelity through block-by-block tracking via excellent supports. 
		Our Def.~\ref{def:cs-iteration} uses $\omega_1$-complete limit filters
		(Theorem~\ref{thm:filter-generation})---via block factorization at countable cofinality 
		and stage-bounding at uncountable cofinality---to preserve ZF and DC
		with countable supports. The mechanisms are fundamentally different: Karagila's 
		productive conditions ensure compositional correctness of finite iterations, 
		while our filter completeness provides the closure properties needed for 
		DC-preservation with infinite supports. The generalization to $<\kappa$-support appears in
		Def.~\ref{def:kappa-iteration} via $\kappa$-complete limit filters.
	\end{remark}
	
	\begin{remark}[Direct limits preserve $\omega_1$-completeness at uncountable cofinality]
		\label{rem:dl-omega1-completeness}
		Assume $\cf(\lambda)\ge\omega_1$. Let $\tilde{\mathcal{F}}_\lambda$ be the filter on $\mathcal{G}_\lambda$
		generated by $\bigcup_{\beta<\lambda}\pi^\beta_{\lambda,*}(\mathcal{F}_\beta)$ in the direct-limit case,
		where $\mathcal{G}_\lambda=\bigcup_{\beta<\lambda}\mathcal{G}_\beta$ and pushforwards are inclusions.
		Then $\tilde{\mathcal{F}}_\lambda$ is $\omega_1$-complete.
		
		Indeed, given a countable family $\{H_n:n<\omega\}\subseteq\tilde{\mathcal{F}}_\lambda$, choose
		$\beta_n<\lambda$ and $K_n\in\mathcal{F}_{\beta_n}$ with $H_n=\pi^{\beta_n}_{\lambda,*}(K_n)=K_n$
		(the equality since, in a direct limit, pushforwards are inclusions).
		The set $\{\beta_n:n<\omega\}$ is countable and hence bounded in $\lambda$ because $\cf(\lambda)\ge\omega_1>\omega$;
		equivalently, if it were unbounded then $\cf(\lambda)\le\omega$, a contradiction. Thus
		$\beta^*:=\sup_n\beta_n<\lambda$.
		
		Since $\beta_n\le\beta^*$ for each $n$, we have $K_n\in\mathcal{F}_{\beta_n}\subseteq\mathcal{F}_{\beta^*}$.
		By the induction hypothesis, $\mathcal{F}_{\beta^*}$ is $\omega_1$-complete, so
		$\bigcap_{n<\omega}K_n\in\mathcal{F}_{\beta^*}\subseteq\tilde{\mathcal{F}}_\lambda$.
		Therefore $\bigcap_{n<\omega}H_n=\bigcap_{n<\omega}K_n\in\tilde{\mathcal{F}}_\lambda$.
	\end{remark}
	
	\begin{remark}[Key Properties to Verify]
		For the iteration to be well-defined, we must verify at each stage:
		\begin{enumerate}
			\item Automorphisms are well-defined
			\item Filter is proper (non-empty)
			\item Filter is normal (closed under conjugation)
			\item Filter is $\omega_1$-complete (limit stages)
		\end{enumerate}
		Items 1--3 follow from the definitions and Theorem~\ref{thm:filter-generation}; Item 4 is the technical heart of the paper.
	\end{remark}
	
	\begin{remark}[Comparison to Karagila's Finite Support]
		In Karagila's finite support framework \cite{karagila2019}, supports are finite and filters need only be closed under finite intersections (automatic). Countable support requires $\omega_1$-completeness, which is nontrivial and the core contribution of this framework.
	\end{remark}
	
	\begin{example}[Iteration Setup]
		\label{ex:cohen-setup}
		For a Cohen iteration to $\omega_1$:
		\begin{itemize}
			\item $\dot{\mathbb{Q}}_\beta = \mathrm{Add}(\omega, 2)$
			\item Automorphisms: finite-support permutations of Cohen coordinates
			\item Filter: cofinite stabilizers (e.g., fix all but finitely many coordinates)
		\end{itemize}
	\end{example}
	
	\subsection{Filter Generation at Limit Stages}
	\label{ss:filter-generation}
	
	Let $\lambda$ be a limit ordinal with $\cf(\lambda)=\omega$, and assume each $\mathcal{F}_\beta$ for $\beta<\lambda$ is a normal, $\omega_1$-complete filter on $\mathcal{G}_\beta$ (by induction). We construct $\mathcal{F}_\lambda$ so that:
	\begin{itemize}
		\item \textbf{Extends predecessors}: Each $\mathcal{F}_\beta$ embeds into $\mathcal{F}_\lambda$ via the natural projection.
		\item \textbf{Normal}: $\mathcal{F}_\lambda$ is closed under conjugation by $\mathcal{G}_\lambda$.
		\item \textbf{$\omega_1$-complete}: $\mathcal{F}_\lambda$ is closed under countable intersections.
		\item \textbf{Generated}: $\mathcal{F}_\lambda$ is generated by pushforwards of previous filters and stabilizers of names.
	\end{itemize}
	
	\begin{definition}[Pushforward Filter]
		\label{def:pushforward-filter}
		For $\beta<\lambda$, let $\pi^\beta_{\lambda}: \mathcal{G}_\lambda \to \mathcal{G}_\beta$ be the projection.  Define the *pushforward filter* $\pi^\beta_{\lambda,*}(\mathcal{F}_\beta)$ on $\mathcal{G}_\lambda$ to be the filter generated by 
		\[
		\bigl\{\{g\in\mathcal{G}_\lambda : \pi^\beta_\lambda(g)\in H\} : H\in\mathcal{F}_\beta\bigr\}.
		\]
		Equivalently, each $H\in\mathcal{F}_\beta$ yields the generator $\{g:\pi^\beta_\lambda(g)\in H\}$ in $\mathcal{G}_\lambda$.
	\end{definition}
	
	\begin{lemma}[Pullbacks and restrictions preserve normality and $\kappa$-completeness]
		\label{lem:push-restrict-preserve}
		Let $\pi:\mathcal G\to\mathcal H$ be a homomorphism and let $\mathcal F$ be a normal, $\kappa$-complete filter on $\mathcal H$.
		Define the pullback filter on $\mathcal G$ by
		\[
		\pi^{*}\mathcal F\ :=\ \{\,K\le \mathcal G:\ \exists H\in\mathcal F\ \text{with}\ \pi^{-1}(H)\le K\,\}.
		\]
		Then $\pi^{*}\mathcal F$ is a normal, $\kappa$-complete filter on $\mathcal G$.
		
		Moreover, if $\iota:\mathcal H\hookrightarrow\mathcal G$ is inclusion and $\mathcal F$ is a normal, $\kappa$-complete filter on $\mathcal G$,
		then the restriction
		\[
		\iota^{*}\mathcal F\ :=\ \{\,K\le\mathcal H:\ \iota(K)\in\mathcal F\,\}
		\]
		is a normal, $\kappa$-complete filter on $\mathcal H$.
	\end{lemma}
	
	\begin{proof}
		We use the standard pullback of a filter of subgroups: set
		\[
		\pi^{*}\mathcal F\ :=\ \{\,K\le \mathcal G:\ \exists H\in\mathcal F\ \text{with}\ \pi^{-1}(H)\le K\,\}.
		\]
		Equivalently, $\pi^{*}\mathcal F$ is the upward closure (under inclusion) of $\{\pi^{-1}(H):H\in\mathcal F\}$.
		We first show that $\pi^{*}\mathcal F$ is a normal, $\kappa$-complete filter on $\mathcal G$, and then treat the restriction case.
		
		\smallskip\noindent\emph{Filter axioms.}
		(i) Since $\mathcal H\in\mathcal F$, we have $\pi^{-1}(\mathcal H)=\mathcal G\in\pi^{*}\mathcal F$.
		(ii) If $K\in\pi^{*}\mathcal F$ and $K\le K'\le\mathcal G$, choose $H\in\mathcal F$ with $\pi^{-1}(H)\le K$. Then $\pi^{-1}(H)\le K'\,$, hence $K'\in\pi^{*}\mathcal F$.
		(iii) ($\kappa$-completeness) Let $\langle K_i:i<\mu\rangle\subseteq \pi^{*}\mathcal F$ with $\mu<\kappa$. For each $i$ choose $H_i\in\mathcal F$ with $\pi^{-1}(H_i)\le K_i$. Then
		\[
		\bigcap_{i<\mu}\pi^{-1}(H_i)\;=\;\pi^{-1}\Bigl(\bigcap_{i<\mu} H_i\Bigr)\ \le\ \bigcap_{i<\mu} K_i.
		\]
		Since $\mathcal F$ is $\kappa$-complete, $H:=\bigcap_{i<\mu}H_i\in\mathcal F$, hence $\pi^{-1}(H)\le \bigcap_{i<\mu}K_i$, so $\bigcap_{i<\mu}K_i\in\pi^{*}\mathcal F$.
		
		\smallskip\noindent\emph{Normality.}
		Let $g\in\mathcal G$ and $K\in\pi^{*}\mathcal F$; fix $H\in\mathcal F$ with $\pi^{-1}(H)\le K$. Then
		\[
		gKg^{-1}\ \supseteq\ g\,\pi^{-1}(H)\,g^{-1}\;=\;\pi^{-1}\!\bigl(\pi(g)\,H\,\pi(g)^{-1}\bigr).
		\]
		Since $\mathcal F$ is normal, $\pi(g)\,H\,\pi(g)^{-1}\in\mathcal F$, and therefore $gKg^{-1}\in\pi^{*}\mathcal F$.
		
		Thus $\pi^{*}\mathcal F$ is a normal, $\kappa$-complete filter on $\mathcal G$.
		
		\smallskip
		For the restriction part, let $\iota:\mathcal H\hookrightarrow\mathcal G$ be the inclusion and put
		\[
		\iota^{*}\mathcal F\ :=\ \{\,K\le\mathcal H:\ \iota(K)\in\mathcal F\,\}\ =\ \mathcal F\cap\Sub(\mathcal H).
		\]
		Then $\iota^{*}\mathcal F$ is upward closed, contains $\mathcal H$, and is closed under $<\kappa$ intersections because
		$\iota\!\left(\bigcap_{i<\mu}K_i\right)=\bigcap_{i<\mu}\iota(K_i)$ and $\mathcal F$ is $\kappa$-complete. Normality is inherited: for $h\in\mathcal H$ and $K\in\iota^{*}\mathcal F$,
		\[
		\iota(hKh^{-1})=h\,\iota(K)\,h^{-1}\in\mathcal F.
		\]
		Hence $\iota^{*}\mathcal F$ is a normal, $\kappa$-complete filter on $\mathcal H$, as required. \qedhere
	\end{proof}
	
	\begin{definition}[Filter at Limit Stage]
		\label{def:limit-filter}
		Define $\tilde{\mathcal{F}}_\lambda$ to be the filter on $\mathcal{G}_\lambda$ generated by
		\[
		\bigcup_{\beta<\lambda} \pi^\beta_{\lambda,*}(\mathcal{F}_\beta)
		\;\;\cup\;\;
		\{\fix_{\mathcal{G}_\lambda}(\dot{x}) : \dot{x}\text{ is a $\mathbb{P}_\lambda$-name with countable support}\},
		\]
		where for any $\mathbb{P}_\lambda$-name $\dot{x}$, 
		\[
		\fix_{\mathcal{G}_\lambda}(\dot{x}) \;=\; \{\,g\in\mathcal{G}_\lambda : g(\dot{x})=\dot{x}\}.
		\]
		In practice we only include generators $\fix_{\mathcal{G}_\lambda}(\dot{x})$ for names $\dot{x}$ whose support (the set of coordinates on which $\dot{x}$ depends) is countable. 
	\end{definition}
	
	\begin{definition}[Supports of conditions and names (countable-support case)]
		\label{def:support-recursive}
		Let $\lambda$ be a limit stage.
		
		\smallskip
		\noindent\emph{Support of a condition.}
		For $p\in\mathbb P_\lambda$, define $\supp(p)\subseteq\lambda$ to be the set of coordinates at which $p$ is nontrivial
		(in the usual sense for the iteration/product presentation at $\lambda$). In particular, in the countable-support presentation,
		$\supp(p)$ is required to be countable.
		
		\smallskip
		\noindent\emph{Support of a name.}
		Define the support of a $\mathbb P_\lambda$--name $\dot x$ by recursion on name-rank:
		\[
		\supp(\dot x)\ :=\ \bigcup\Bigl\{\ \supp(p)\ \cup\ \supp(\sigma)\ :\ (\sigma,p)\in\dot x\ \Bigr\}.
		\]
		We say that $\dot x$ has \emph{countable support} if $\supp(\dot x)$ is countable in the ambient background universe.
	\end{definition}
	
	\begin{remark}[Syntactic support vs.\ background countability]
		By Definition~\ref{def:support-recursive}, $\supp(\dot x)$ is computed recursively from the supports of the conditions
		appearing in $\dot x$ (and from the supports of subnames). When we say “$\dot x$ has countable support’’ in the
		countable-support part of the paper, we mean that $\supp(\dot x)$ is countable in the ambient background universe
		(cf.\ Remark~\ref{rem:countable-conv}).
	\end{remark}
	
	The $\omega_1$--completeness argument splits by cofinality. If $\cf(\lambda)\ge\omega_1$, stage bounding
	reduces countable intersections to a single stage filter. If $\cf(\lambda)=\omega$, stage bounding fails
	and we instead use tuple--stabilizer compression at the inverse limit.
	
	The proof of $\omega_1$-completeness requires different arguments depending on the cofinality of $\lambda$:
	\begin{itemize}
		\item \textbf{Case $\cf(\lambda)\ge\omega_1$ (direct limit).} Here stage-bounding applies: any countable set of ordinals below $\lambda$ is bounded, since a countable cofinal subset would contradict $\cf(\lambda)\ge\omega_1$. This case is handled in Remark~\ref{rem:dl-omega1-completeness}.
		\item \textbf{Case $\cf(\lambda)=\omega$ (inverse limit).} Stage-bounding \emph{fails}: a countable set can be cofinal in $\lambda$. We instead use \emph{block factorization} along a fixed cofinal sequence.
	\end{itemize}
	
	\begin{lemma}[Stage-Bounding Lemma---uncountable cofinality only]
		\label{lemma:stage-bounding}
		If $\lambda$ is a limit ordinal with $\cf(\lambda)\ge\omega_1$ and $\{\beta_n:n<\omega\}\subset\lambda$ is any countable set of ordinals below $\lambda$, then 
		\[
		\sup_{n<\omega}\beta_n \;<\;\lambda.
		\]
	\end{lemma}
	
	\begin{proof}
		If $\sup_n\beta_n=\lambda$, then $\{\beta_n:n<\omega\}$ would be a countable cofinal subset of $\lambda$, contradicting $\cf(\lambda)\ge\omega_1$. Hence $\sup_n\beta_n<\lambda$.
	\end{proof}
	
	\begin{remark}[Failure of stage-bounding at countable cofinality]
		\label{rem:stage-bounding-fails}
		When $\cf(\lambda)=\omega$, the stage-bounding lemma fails. For example, if $\lambda=\omega$, the set $\{n:n<\omega\}$ is countable but $\sup_n n=\omega=\lambda$. Similarly, if $\lambda=\aleph_\omega$, the sequence $\{\aleph_n:n<\omega\}$ is countable and cofinal. This is precisely why we need stabilizer compression rather than stage-bounding at limits of countable cofinality.
	\end{remark}
	
	\begin{theorem}[Filter Extension Theorem]
		\label{thm:filter-generation}
		The constructed filter $\tilde{\mathcal{F}}_\lambda$ is well-defined (proper, nonempty), normal, and $\omega_1$-complete.  Moreover, it extends each predecessor via the natural embeddings $\pi^\beta_{\lambda,*}(\mathcal{F}_\beta)\hookrightarrow \tilde{\mathcal{F}}_\lambda$.
	\end{theorem}
	\begin{proof}
		\textbf{Properness:}  Each generator in Definition~\ref{def:limit-filter} contains the identity automorphism (since the identity belongs to every $H\in\mathcal{F}_\beta$ and fixes every name). Hence finite intersections of generators are nonempty, so $\tilde{\mathcal{F}}_\lambda$ is a proper (nonempty) filter.
		
		\textbf{Normality:}  Consider a generator of two types: (1) If $H=\{g:\pi^\beta_\lambda(g)\in H_\beta\}$ comes from a pushforward of $H_\beta\in\mathcal{F}_\beta$, then for any $k\in\mathcal{G}_\lambda$, 
		$$
		kHk^{-1} \;=\; \{g: \pi^\beta_\lambda(k^{-1} g k)\in H_\beta\} \;=\; \{g:\pi^\beta_\lambda(g)\in (\pi^\beta_\lambda(k) H_\beta \pi^\beta_\lambda(k)^{-1})\}.
		$$
		Since $\mathcal{F}_\beta$ is normal, $\pi^\beta_\lambda(k) H_\beta \pi^\beta_\lambda(k)^{-1}\in\mathcal{F}_\beta$, so $kHk^{-1}$ is again a generator of the same type.  (2) If $H = \fix_{\mathcal{G}_\lambda}(\dot{x})$ is a stabilizer of a name, then for any $k\in\mathcal{G}_\lambda$ one checks
		$$
		k\,\fix(\dot{x})\,k^{-1} \;=\; \fix\bigl(k(\dot{x})\bigr).
		$$
		This holds because $g$ fixes $\dot{x}$ if and only if $kgk^{-1}$ fixes $k(\dot{x})$ (by conjugation). 
		Since $k(\dot{x})$ is another $\mathbb{P}_\lambda$-name with the same countable support as $\dot{x}$, 
		$\fix(k(\dot{x}))$ is also one of our generators. 
		Thus, $kHk^{-1}$ is a generator of $\tilde{\mathcal{F}}_\lambda$, and $\tilde{\mathcal{F}}_\lambda$ 
		is closed under conjugation by $\mathcal{G}_\lambda$, i.e., it is normal.
		
		\textbf{$\omega_1$-completeness:} We treat the two cofinality cases separately.
		
		\smallskip
		\noindent\emph{Case $\cf(\lambda)\ge\omega_1$ (direct limit; stage-bounding).}
		This case is handled by Remark~\ref{rem:dl-omega1-completeness}: any countable collection of stages is bounded below $\lambda$ by Lemma~\ref{lemma:stage-bounding}, so the argument there applies.
		
		\smallskip
		\noindent\emph{Case $\cf(\lambda)=\omega$ (inverse limit; stabilizer compression).}
		We prove simultaneously by induction on $\beta\le\lambda$ the following stabilizer-basis property:
		for every $H\in\tilde{\mathcal F}_\beta$ there is a countably supported $\mathbb P_\beta$--name $\dot y$
		with $\fix_{\mathcal G_\beta}(\dot y)\le H$.
		The desired $\omega_1$--completeness at $\beta=\lambda$ follows from this property and \eqref{eq:tuple-fix}.
		
		Let $\{H_n:n<\omega\}\subseteq\tilde{\mathcal F}_\lambda$.
		We show $\bigcap_{n<\omega}H_n\in\tilde{\mathcal F}_\lambda$.
		
		\smallskip
		\noindent\emph{Sublemma (stabilizer basis).}
		For every $H\in\tilde{\mathcal F}_\lambda$ there exists a $\mathbb P_\lambda$--name $\dot y$ with countable
		support such that $\fix_{\mathcal G_\lambda}(\dot y)\le H$.
		
		\smallskip
		\noindent\emph{Proof of Sublemma.}
		Since $\tilde{\mathcal F}_\lambda$ is generated by the families in Definition~\ref{def:limit-filter},
		there are basic generators $B_0,\dots,B_{k-1}$ such that $\bigcap_{i<k}B_i\le H$.
		For each $i<k$ we define a stabilizer $S_i$ with $S_i\le B_i$:
		\begin{itemize}
			\item If $B_i=\fix_{\mathcal G_\lambda}(\dot x)$ for a countably supported name $\dot x$, let $S_i:=\fix_{\mathcal G_\lambda}(\dot x)$.
			\item If $B_i=\{g\in\mathcal G_\lambda:\pi^\beta_\lambda(g)\in K\}$ for some $\beta<\lambda$ and $K\in\mathcal F_\beta$,
			then (by the inductive hypothesis at stage $\beta$) there is a countably supported $\mathbb P_\beta$--name $\dot z$
			with $\fix_{\mathcal G_\beta}(\dot z)\le K$. View $\dot z$ as a $\mathbb P_\lambda$--name in the standard way.
			Then $\fix_{\mathcal G_\lambda}(\dot z)\le \{g:\pi^\beta_\lambda(g)\in K\}=B_i$, and let $S_i:=\fix_{\mathcal G_\lambda}(\dot z)$.
		\end{itemize}
		Let $\dot y$ be the finite tuple of the names used to define the $S_i$'s. Then
		$\fix_{\mathcal G_\lambda}(\dot y)=\bigcap_{i<k}S_i\le \bigcap_{i<k}B_i\le H$.
		Since $\dot y$ is a finite tuple of countably supported names, it has countable support.
		This proves the Sublemma.
		
		\smallskip
		\noindent
		Now apply the Sublemma to each $H_n$ and obtain countably supported names $\dot y_n$ with
		$\fix(\dot y_n)\le H_n$.
		Let $\dot y:=\langle \dot y_n:n<\omega\rangle$ be the countable tuple. By \eqref{eq:tuple-fix},
		\[
		\fix(\dot y)=\bigcap_{n<\omega}\fix(\dot y_n)\ \le\ \bigcap_{n<\omega}H_n.
		\]
		The name $\dot y$ has countable support (a countable union of countable supports, in the ambient background, Remark~\ref{rem:countable-conv}),
		so $\fix(\dot y)$ is one of the stabilizer-type generators in Definition~\ref{def:limit-filter}, hence
		$\fix(\dot y)\in\tilde{\mathcal F}_\lambda$. By upward closure, $\bigcap_{n<\omega}H_n\in\tilde{\mathcal F}_\lambda$.
		
		\medskip
		\noindent\emph{Key distinction.}
		For $\cf(\lambda)\ge\omega_1$ we use \emph{stage bounding} (Remark~\ref{rem:dl-omega1-completeness}).
		For $\cf(\lambda)=\omega$ we use \emph{stabilizer compression}: every element of $\tilde{\mathcal F}_\lambda$
		contains a stabilizer of a countably supported name, and countable intersections of such stabilizers
		are stabilizers of a countably supported tuple by \eqref{eq:tuple-fix}.
		
		\smallskip
		This completes the proof that $\tilde{\mathcal{F}}_\lambda$ is a normal, $\omega_1$-complete filter extending all earlier $\mathcal{F}_\beta$.
	\end{proof}
	
	\subsection{Hereditarily Symmetric Names in Countable Support Iterations}
	\label{sec:HS-countable}
	
	Having constructed the filter $\tilde{\mathcal{F}}_\lambda$ at limit stages---via block factorization when $\cf(\lambda)=\omega$ (inverse limit) and stage-bounding when $\cf(\lambda)\ge\omega_1$ (direct limit)---and established its $\omega_1$-completeness (Section~\ref{ss:filter-generation}, Theorem~\ref{thm:filter-generation}), we now define the class of hereditarily symmetric names and verify its basic closure properties in the countable-support setting.
	
	\begin{definition}[Hereditarily Symmetric Names --- Iterated, Countable Support Case]
		\label{def:hs-names-iterated}
		A $\mathbb{P}_\lambda$-name $\dot{x}$ is \emph{hereditarily symmetric} (HS) (with respect to $\tilde{\mathcal{F}}_\lambda$) if:
		\begin{enumerate}
			\item $\fix_{\mathcal{G}_\lambda}(\dot{x}) \in \tilde{\mathcal{F}}_\lambda$, and
			\item for every $(\dot{y},p)\in \dot{x}$, the name $\dot{y}$ is HS.
		\end{enumerate}
	\end{definition}
	
	\begin{equation}\label{eq:tuple-fix}
		\fix\!\bigl(\langle \dot x_j : j\in J\rangle\bigr)
		\;=\; \bigcap_{j\in J}\fix(\dot x_j).
	\end{equation}
	
	\begin{lemma}[Countable-support characterization of HS (support via pointwise stabilizers)]
		\label{lem:countable-support-hs}
		Work in the ambient background universe in which the iteration is constructed
		(cf.\ Remark~\ref{rem:countable-conv}).
		Let $\lambda$ be a limit stage and let $\tilde{\mathcal F}_\lambda$ be the limit filter
		from Definition~\ref{def:limit-filter}.
		
		\begin{enumerate}
			\item[(1)] If $\dot x$ is HS (i.e.\ $\fix_{\mathcal G_\lambda}(\dot x)\in\tilde{\mathcal F}_\lambda$),
			then there exists a \emph{countable} set $S\subseteq\lambda$ such that
			\[
			\mathrm{Fix}(S)\ :=\ \{g\in\mathcal G_\lambda:\forall \gamma\in S\ (g\text{ fixes coordinate }\gamma)\}
			\ \le\ \fix_{\mathcal G_\lambda}(\dot x).
			\]
			In particular, $\dot x$ has a \emph{countable support} in the sense that it is fixed by
			all automorphisms fixing $S$ pointwise.
			
			\item[(2)] Conversely, if $\dot x$ has \emph{countable support} in the sense of
			Definition~\ref{def:support-recursive} and every name occurring in $\dot x$ is HS, then $\dot x$ is HS.
		\end{enumerate}
	\end{lemma}
	
	\begin{proof}
		(1) Assume $\dot x$ is HS, so $\fix(\dot x)\in\tilde{\mathcal F}_\lambda$.
		By the stabilizer-basis sublemma proved in Theorem~\ref{thm:filter-generation}
		(in the $\cf(\lambda)=\omega$ case), there is a $\mathbb P_\lambda$--name $\dot y$
		with \emph{countable} support (in the sense of Definition~\ref{def:limit-filter})
		such that
		\[
		\fix(\dot y)\ \le\ \fix(\dot x).
		\]
		Let $S:=\supp(\dot y)$ (countable). Any $g\in\mathrm{Fix}(S)$ fixes every condition whose
		support is contained in $S$, hence fixes $\dot y$; therefore
		$\mathrm{Fix}(S)\le \fix(\dot y)\le \fix(\dot x)$, as required.
		
		\smallskip
		(2) If $\dot x$ has countable support (Definition~\ref{def:support-recursive}), then by
		Definition~\ref{def:limit-filter} the stabilizer $\fix_{\mathcal G_\lambda}(\dot x)$ is one of the
		stabilizer-type generators of $\tilde{\mathcal F}_\lambda$, hence belongs to $\tilde{\mathcal F}_\lambda$.
		Hereditariness holds by the hypothesis on constituents.
	\end{proof}
	
	\begin{remark}[What “support” means in permutation arguments]
		\label{rem:support-via-fixS}
		In the countable-support part of the paper we use “$\dot x$ has countable support” in the
		\emph{stabilizer} sense: there exists a countable set of coordinates $S\subseteq\lambda$ such that
		$\mathrm{Fix}(S)\le \fix_{\mathcal G_\lambda}(\dot x)$, where $\mathrm{Fix}(S)$ is the pointwise stabilizer of $S$.
		This is the notion actually needed in the standard symmetry/permutation arguments:
		if $\mathrm{Fix}(S)\le\fix(\dot x)$ and $\alpha,\beta\notin S$, then the coordinate-swap automorphism
		fixing $S$ pointwise lies in $\fix(\dot x)$, so it fixes $\dot x$ while exchanging the two coordinate generics.
		We do \emph{not} rely on any internal $\ZF$ principle asserting that a countable union of countable sets
		is countable; all “countability’’ assertions here are made in the ambient background universe in which the
		iteration is carried out (cf.\ Remark~\ref{rem:countable-conv}).
	\end{remark}
	
	\begin{theorem}[HS closure under standard name operations]
		\label{thm:hs-closure}
		The class of HS names is closed under pairing, union, formation of subsets by bounded formulas (Separation instances over a fixed ambient set), and the usual forcing operations (e.g.\ check names).
	\end{theorem}
	
	\begin{proof}
		\emph{Pairing.}
		If $\dot{x},\dot{y}$ are HS, then by the tuple–stabilizer identity (Equation~\eqref{eq:tuple-fix})
		\[
		\fix\!\bigl(\langle \dot{x},\dot{y}\rangle\bigr)\;=\;\fix(\dot{x})\cap\fix(\dot{y}),
		\]
		and since $\tilde{\mathcal{F}}_\lambda$ is a filter (closed under finite intersections), the pair is HS.
		
		\emph{Countable tuples.}
		Let $\langle \dot x_n:n<\omega\rangle$ be a countable sequence of HS names, and let $\dot s$ denote the
		standard $\mathbb P_\lambda$--name for this tuple.
		By \eqref{eq:tuple-fix},
		\[
		\fix(\dot s)\ =\ \bigcap_{n<\omega}\fix(\dot x_n).
		\]
		Since each $\fix(\dot x_n)\in\tilde{\mathcal F}_\lambda$ and $\tilde{\mathcal F}_\lambda$ is $\omega_1$--complete
		(Theorem~\ref{thm:filter-generation}), we have $\fix(\dot s)\in\tilde{\mathcal F}_\lambda$, so $\dot s$ is HS.
		
		\emph{Union of an HS name.}
		Let $\dot A$ be HS and define the usual union-name
		\[
		\dot U\ :=\ \Bigl\{\,\langle \tau,r\rangle : \exists\langle\sigma,p\rangle\in\dot A\ \exists\langle\tau,q\rangle\in\sigma\
		(r\le p,q)\Bigr\}.
		\]
		If $g\in\fix(\dot A)$ then $g(\dot U)=\dot U$ (because $g$ just transports witnesses in the defining formula),
		so $\fix(\dot U)\supseteq\fix(\dot A)\in\tilde{\mathcal F}_\lambda$. Hence $\dot U$ is HS.
		
		\emph{Bounded Separation (over a fixed ambient set).}
		Let $\dot{B}=\{z\in\dot{A}:\varphi(z,\bar{u})\}$ where $\dot{A}$ is HS and $\bar{u}$ is a tuple of HS parameters,
		and $\varphi$ is a bounded formula.
		If $g$ fixes $\dot{A}$ and $\bar{u}$, then by automorphism–invariance of truth for bounded formulas in the forcing
		language, $g$ fixes $\dot{B}$. Hence
		\[
		\fix(\dot{B})\ \supseteq\ \fix(\dot{A})\cap\fix(\bar{u}),
		\]
		and the right-hand side lies in $\tilde{\mathcal{F}}_\lambda$ because $\tilde{\mathcal F}_\lambda$ is a filter
		(closed under finite intersections, and $\bar u$ is a finite tuple of parameters),
		so $\dot{B}$ is HS.
		
		\emph{Check names and standard operations.}
		If $a\in V$, then $\check{a}$ is fixed by all $g\in\mathcal{G}_\lambda$, so $\fix(\check{a})=\mathcal{G}_\lambda\in\tilde{\mathcal{F}}_\lambda$.
		Other usual name constructors are handled similarly; when a constructor uses only finitely many HS inputs,
		we use finite-intersection closure of $\tilde{\mathcal F}_\lambda$, and when it packages countably many HS inputs
		(e.g.\ countable tuples) we use \eqref{eq:tuple-fix} together with $\omega_1$--completeness
		(Theorem~\ref{thm:filter-generation}).
		
		The hereditary clause in each case follows by induction on name rank. \qedhere
	\end{proof}
	
	\begin{remark}[Ablation checkpoints for Theorem~\ref{thm:hs-closure}]
		\leavevmode
		\begin{enumerate}[label=(\alph*), leftmargin=2.2em]
			\item Drop $\omega_1$/\,$\kappa$–completeness: Replacement’s final stabilizer intersection may leave the filter.
			\item Drop normality: Separation can fail (stabilizers of parameters + set need not generate a filter element).
			\item Drop support bounds: the tuple–stabilizer reduction no longer yields a generator.
		\end{enumerate}
	\end{remark}
	
	\begin{remark}[Relation to finite-support symmetry]
		In the finite-support setting (e.g.\ Karagila), filter closure under finite intersections suffices to obtain the HS closure properties. Our countable-support regime requires $\omega_1$-completeness to handle countable intersections of stabilizers arising from countably supported names; this is exactly where the Theorem~\ref{thm:filter-generation} result is used.
	\end{remark}
	
	
	\subsection{$\ZF$ Preservation}
	\label{sec:zf-preservation}
	
	\begin{theorem}[$\ZF$ preservation for the symmetric submodel]
		\label{thm:zf-preserved}
		Assume we are working in the ambient background universe used for the iteration construction
		(cf.\ Remark~\ref{rem:countable-conv}). Let $V$ be the ground model for the forcing and assume $V\models\ZF$.
		Let $\lambda$ be a stage of the iteration and let
		$(\mathbb P_\lambda,\mathcal G_\lambda,\tilde{\mathcal F}_\lambda)$ be the associated symmetric system,
		where $\tilde{\mathcal F}_\lambda$ is a normal $\omega_1$--complete filter of subgroups of $\mathcal G_\lambda$
		(e.g.\ the limit filter from Theorem~\ref{thm:filter-generation} in the countable--cofinality case).
		Let $G\subseteq\mathbb P_\lambda$ be $V$--generic, and let $\HS$ denote the class of hereditarily
		$(\mathcal G_\lambda,\tilde{\mathcal F}_\lambda)$--symmetric $\mathbb P_\lambda$--names
		(Definition~\ref{def:hs-names-iterated}). Define
		\[
		M \ :=\ V(\mathbb P_\lambda)_{\HS}\ :=\ \{\dot x^{\,G}:\dot x\in\HS\}.
		\]
		Then $M$ is a transitive inner model of $\ZF$ and
		\[
		V\ \subseteq\ M\ \subseteq\ V[G].
		\]
	\end{theorem}
	
	\begin{proof}
		Write $\fix(\dot x)=\fix_{\mathcal G_\lambda}(\dot x)$ and let $\tilde{\mathcal F}_\lambda$ be the normal,
		$\omega_1$--complete limit filter from Theorem~\ref{thm:filter-generation} (constructed at $\lambda$ as in
		Definition~\ref{def:limit-filter}). We also use the tuple--stabilizer identity \eqref{eq:tuple-fix}.
		
		\smallskip\noindent\textbf{Transitivity and sandwich.}
		If $a\in V$ then $\check a$ is fixed by all $g\in\mathcal G_\lambda$, hence HS and $a=\check a^G\in M$; thus $V\subseteq M$.
		Trivially $M\subseteq V[G]$. For transitivity, let $y\in M$ with $y=\dot y^{G}$ for some HS name $\dot y$.
		If $x\in y$ then $x=\sigma^G$ for some $(\sigma,p)\in\dot y$ with $p\in G$. Since $\dot y$ is \emph{hereditarily} symmetric,
		every constituent $\sigma$ of $\dot y$ is HS, hence $x\in M$. Therefore $M$ is transitive.
		
		\smallskip\noindent\textbf{Extensionality and Foundation.}
		These are standard for symmetric submodels: $M$ is a transitive subclass of $V[G]$ with the inherited $\in$,
		and the forcing interpretation of names is extensional and well-founded; see, e.g., \cite{Jech03}. Thus $M$
		satisfies Extensionality and Foundation.
		
		\smallskip\noindent\textbf{Pairing and Union.}
		If $\dot x,\dot y$ are HS, then by \eqref{eq:tuple-fix}
		$\fix(\langle \dot x,\dot y\rangle)=\fix(\dot x)\cap\fix(\dot y)\in\tilde{\mathcal F}_\lambda$,
		so the ordered pair is HS; hereditariness gives $\{\dot x,\dot y\}$ HS as well.
		
		For Union, let $\dot A$ be HS and define the usual union-name
		\[
		\dot U\ :=\ \Bigl\{\,\langle \tau,r\rangle : \exists\langle\sigma,p\rangle\in\dot A\ \exists\langle\tau,q\rangle\in\sigma\
		(r\le p,q)\Bigr\}.
		\]
		If $g\in\fix(\dot A)$ then $g(\dot U)=\dot U$ (since $g$ transports witnesses in the defining formula),
		so $\fix(\dot U)\supseteq\fix(\dot A)\in\tilde{\mathcal F}_\lambda$. Hence $\dot U$ is HS and therefore
		$U:=\dot U^G=\bigcup(\dot A^G)$ belongs to $M$.
		
		\smallskip\noindent\textbf{Infinity.}
		$\check\omega$ is fixed by all automorphisms, hence HS; thus Infinity holds in $M$.
		
		\smallskip\noindent\textbf{Separation.}
		Let $A\in M$ be given by an HS name $\dot A$, fix HS parameters $\vec{\dot p}$, and let $\varphi$ be any formula.
		Define the name
		\[
		\dot B\ :=\ \{\ \langle\sigma,p\rangle\in\dot A\ :\ p\Vdash \varphi(\sigma,\vec{\dot p})\ \}.
		\]
		If $g\in\fix(\dot A)\cap\fix(\vec{\dot p})$, then by automorphism invariance of the forcing relation,
		$p\Vdash \varphi(\sigma,\vec{\dot p})\iff g(p)\Vdash \varphi(g\sigma,\vec{\dot p})$, whence $g\in\fix(\dot B)$.
		Thus $\fix(\dot B)\supseteq \fix(\dot A)\cap\fix(\vec{\dot p})\in\tilde{\mathcal F}_\lambda$, so $\dot B$ is HS and
		$B:=\dot B^G=\{x\in A:\varphi^M(x,\vec p)\}$. Hence $M$ satisfies Separation.
		
		\smallskip\noindent\textbf{Power Set.}
		Let $a\in M$ with $a=\dot a^{G}$ for HS $\dot a$. Fix $\alpha:=\rank(\dot a)+1$ and let
		$\mathrm{Nm}_{\le\alpha}$ be the (set-sized) collection in $V$ of $\mathbb P_\lambda$–names of rank $\le\alpha$.
		Set
		\[
		\HS(\dot a)\ :=\ \{\ \dot x\in \mathrm{Nm}_{\le\alpha}\ :\ \Vdash \dot x\subseteq \dot a\ \text{ and }\ \dot x\in\mathrm{HS}\ \}.
		\]
		Then $\HS(\dot a)$ is a \emph{set} in $V$ (it is a definable subset of $\mathrm{Nm}_{\le\alpha}$), and if
		$g\in\fix(\dot a)$ then $g$ preserves both subsetting $\dot a$ and the HS property, so $\HS(\dot a)$ is invariant under
		$\mathcal G_\lambda$. Consequently,
		\[
		\mathcal P^M(a)\ =\ \{\ \dot x^{G} : \dot x\in \HS(\dot a)\ \}\ \in\ M,
		\]
		witnessed by the HS name coding $\HS(\dot a)$ (using standard name-coding for sets of names).
		Therefore $M$ satisfies Power Set.
		
		\smallskip\noindent\textbf{Replacement (via relational image names).}
		Let $A\in M$ and $f\in M$ be such that $f$ is a \emph{function} in $M$. Choose HS names $\dot A,\dot f$ with
		$\Vdash$ “$\dot f$ is functional on $\dot A$”. Define in $V$ the \emph{relational range} name
		\[
		\dot{\mathrm{rng}}\ :=\ \Bigl\{\ \bigl\langle\tau,p\bigr\rangle\ :\ \exists\langle\sigma,q\rangle\in\dot A\ \exists r\le p,q\
		\bigl\langle\langle\sigma,\tau\rangle, r\bigr\rangle\in\dot f\ \Bigr\}.
		\]
		If $g\in\fix(\dot A)\cap\fix(\dot f)$, then by definition $g\in\fix(\dot{\mathrm{rng}})$; hence
		$\fix(\dot{\mathrm{rng}})\supseteq \fix(\dot A)\cap\fix(\dot f)\in\tilde{\mathcal F}_\lambda$, so $\dot{\mathrm{rng}}$ is HS.
		By the forcing theorem and functionality of $\dot f$, we have
		\[
		\dot{\mathrm{rng}}^{G}\ =\ \{\, y\in V[G] : \exists x\in \dot A^{G}\ (\langle x,y\rangle\in \dot f^{G})\,\}\ =\ f``A.
		\]
		Thus $f``A\in M$, which is Collection and hence Replacement in $M$.
		
		\smallskip
		We have verified that $M$ is a transitive inner model of $\ZF$ with $V\subseteq M\subseteq V[G]$. \qedhere
	\end{proof}
	
	\begin{remark}[On the role of $\omega_1$-completeness]
		In the countable-support setting, $\omega_1$-completeness of $\tilde{\mathcal F}_\lambda$ is needed to close HS under
		countable intersections of stabilizers (e.g., for closure under countable unions of HS constituents).
		If one proves Replacement by first enumerating the domain and setting $\dot B=\{\dot y_n:n<\omega\}$, then
		$\omega_1$-completeness is used to ensure $\bigcap_{n<\omega}\fix(\dot y_n)\in\tilde{\mathcal F}_\lambda$.
		
		In the proof of Theorem~\ref{thm:zf-preserved} given above, we instead use the \emph{relational range} name,
		so the Replacement step itself only requires the finite intersection
		$\fix(\dot A)\cap\fix(\dot f)\in\tilde{\mathcal F}_\lambda$. Thus, for countable support the real novelty
		manifests earlier (in closure under countable unions and similar uses of $\omega_1$-completeness),
		rather than in the Replacement step per se. \emph{Mutatis mutandis} for $<\kappa$-support: replace
		$\omega_1$-completeness by $<\!\kappa$-completeness.
	\end{remark}
	
	
	\subsection{Dependent Choice and Choice Fragments}
	\label{sec:dc-choice}
	
	\begin{definition}[Dependent Choice (DC)]
		For a nonempty set $A$ and relation $R\subseteq A\times A$ with $\forall x\in A\,\exists y\in A\ (xRy)$, DC asserts the existence of a sequence $(a_n)_{n<\omega}$ with $a_nRa_{n+1}$ for all $n<\omega$.
	\end{definition}
	
	\begin{theorem}[DC Preservation]
		\label{thm:dc-preservation}
		Recall the ambient background $V\models\mathrm{ZFC}$ and $\tilde{\mathcal{F}}_\lambda$ is $\omega_1$-complete, then $M=V(\mathbb{P}_\lambda)_{\mathrm{HS}} \models \mathrm{DC}$.
	\end{theorem}
	
	\begin{proof}
		Let $\dot{A}$ and $\dot{R}$ be HS names with
		\[
		\Vdash\ \forall \dot{x}\in \dot{A}\ \exists \dot{y}\in \dot{A}\,(\dot{x},\dot{y})\in \dot{R}.
		\]
		
		\textbf{Construction.}
		Build by induction a sequence $(\dot{a}_n,p_n)_{n<\omega}$ where:
		\begin{itemize}
			\item each $\dot{a}_n$ is an HS name for an element of $\dot{A}$,
			\item $p_n\in\mathbb{P}_\lambda$ has countable support $S_n\subseteq\lambda$,
			\item $p_{n+1}\le p_n$ and $p_{n+1}\Vdash (\dot{a}_n,\dot{a}_{n+1})\in\dot{R}$.
		\end{itemize}
		
		\textbf{Countable support.}
		In the ambient background (ZFC), $S$ is countable, hence the fusion condition is countably supported.
		
		\textbf{Key step (stabilizers).}
		Since each $\dot a_n$ is HS, $\fix(\dot a_n)\in\tilde{\mathcal F}_\lambda$.
		By $\omega_1$-completeness (Theorem~\ref{thm:filter-generation}),
		\[
		\fix\!\bigl(\langle \dot a_n:n<\omega\rangle\bigr)
		=\bigcap_{n<\omega}\fix(\dot a_n)\ \in\ \tilde{\mathcal F}_\lambda.
		\]
		Therefore the name $\langle \dot a_n:n<\omega\rangle$ is HS.
		
		\textbf{Conclusion.}
		Therefore the name $\langle \dot{a}_n:n<\omega\rangle$ is HS, so its valuation
		in the symmetric model yields a sequence witnessing DC in $M$. Moreover, $\neg\mathsf{AC}$ in $M$ is witnessed independently in \S\ref{sec:examples} by
		Theorems~\ref{thm:cohen-ac-fails} and~\ref{thm:random-ac-fails}.
	\end{proof}
	
	\begin{remark}[Scope of DC in the countable-support setting]
		\label{rem:dc-countable-scope}
		In the countable-support case, $\omega_1$-completeness of the limit filter yields $DC_{<\omega_1}$; in particular $DC=DC_\omega$ holds in the symmetric model. Obtaining $DC_{\omega_1}$ requires an additional localization hypothesis analogous to Definition~\ref{def:localization-hypothesis}.
	\end{remark}
	
	\begin{example}[Failure of AC with DC preserved]
		\label{ex:ac-fails}
		In the symmetric extension obtained from an $\omega_1$-length countable-support iteration of Cohen forcing, there is no well-ordering of the reals (standard argument via symmetry), so AC fails. Nevertheless, by Theorem~\ref{thm:dc-preservation} and $\omega_1$-completeness of $\tilde{\mathcal{F}}_\lambda$, DC holds in $M$.
	\end{example}
	
	\begin{remark}[Choice landscape and filter closure]
		\label{rem:choice-hierarchy}
		Filter closure dictates which fragments of choice survive:
		\begin{itemize}
			\item If $\tilde{\mathcal{F}}_\lambda$ is $\omega$-closed (closed under countable decreasing sequences), then $\mathrm{AC}_\omega$ tends to be preserved.
			\item If $\tilde{\mathcal{F}}_\lambda$ is $\omega_1$-complete (closed under countable intersections), then DC is preserved (Theorem~\ref{thm:dc-preservation}).
			\item More generally, if $\tilde{\mathcal{F}}_\lambda$ is $\kappa$-complete then $M$ preserves $\mathrm{DC}_{<\kappa}$ (Theorem~\ref{thm:DCless-preservation}); under the Localization Hypothesis (Definition~\ref{def:localization-hypothesis}) one also obtains $\mathrm{DC}_\kappa$ (Theorem~\ref{thm:DCkappa-localization}).
		\end{itemize}
		In our countable-support framework the natural target is $\omega_1$-completeness, yielding DC while typically destroying full AC (Example~\ref{ex:ac-fails}).
	\end{remark}
	
	\section{Worked Examples and Comparative Analysis}
	\label{sec:examples}
	
	We illustrate the framework from Section~\ref{ss:filter-generation} and Section~\ref{sec:zf-preservation}
	by working out concrete symmetric iterations under \emph{countable support}. Throughout,
	recall that hereditarily symmetric (HS) names are defined via stabilizers of names
	(Definition~\ref{def:hs-names-iterated}), and that $\tilde{\mathcal{F}}_\lambda$ is
	a normal, $\omega_1$-complete filter (Theorem~\ref{thm:filter-generation}).
	We will freely use the HS-closure results from Section~\ref{sec:HS-countable} (Theorem~\ref{thm:hs-closure} and Lemma~\ref{lem:countable-support-hs})
	and the ZF/DC preservation from Sections~\ref{sec:zf-preservation}--\ref{sec:dc-choice}
	(Theorems~\ref{thm:zf-preserved} and~\ref{thm:dc-preservation}).
	
	\begin{remark}[Productivity of the worked examples]
		\label{rem:productivity-examples}
		The presentations of Cohen and random (and the Hechler clause in \S\ref{subsec:hechler-ltkappa})
		can be arranged to be \emph{productive} in the sense of \cite[Def.~8.1]{karagila2019} by choosing
		the standard dense sets of respected names. We do not use productivity explicitly in our proofs,
		but this alignment simplifies comparison with \cite{karagila2019}.
	\end{remark}
	
	\subsection{Cohen Forcing under Countable Support}
	\label{sec:cohen-example}
	
	Fix an ordinal $\lambda$ and form the countable-support iteration
	$\langle \mathbb P_\alpha,\dot{\mathbb Q}_\alpha,\alpha<\lambda\rangle$ as in
	Definition~\ref{def:cs-iteration}: at limits of cofinality $\omega$ use the inverse-limit
	presentation (trees/tuples), and at limits of uncountable cofinality use direct limits.
	At each successor stage $\alpha{+}1<\lambda$, force with
	\[
	\dot{\mathbb Q}_\alpha=\Add(\omega,\omega),
	\]
	viewed as the forcing of finite partial functions $\omega\times\omega\to 2$ ordered by reverse inclusion.
	Let $\mathcal G_{\alpha+1}$ be the group of finite-support permutations of the $\omega$ columns, acting by
	permuting the second coordinate, and let $\mathcal F_{\alpha+1}$ be the normal filter generated by the
	pointwise fixators of cofinite $X\subseteq\omega$.
	At a limit stage, let $\tilde{\mathcal F}_\lambda$ be the filter from Definition~\ref{def:limit-filter}.
	By Theorem~\ref{thm:filter-generation}, $\tilde{\mathcal F}_\lambda$ is normal and $\omega_1$-complete.
	
	For $\alpha<\lambda$ and $n<\omega$, let $\dot c_{\alpha,n}$ denote the canonical name for the Cohen real
	given by the $n$th column added at stage $\alpha$. Then each $\dot c_{\alpha,n}$ is HS: it is supported on
	the single stage $\{\alpha\}$ (hence countable), so by Lemma~\ref{lem:countable-support-hs} we have
	$\fix(\dot c_{\alpha,n})\in\tilde{\mathcal{F}}_\lambda$.
	
	\begin{example}[$\omega_1$-length countable-support Cohen iteration]
		\label{ex:omega1-cohen}
		Specialize now to $\lambda=\omega_1$. For each $\alpha<\omega_1$ let $\dot{\mathbb Q}_\alpha=\Add(\omega,\omega)$,
		and form the countable-support iteration $\langle\mathbb P_\alpha,\dot{\mathbb Q}_\alpha,\alpha<\omega_1\rangle$:
		at successors use the usual two-step iteration; at limits of countable cofinality use the inverse limit with
		countable support; and at $\omega_1$ (cofinality $\omega_1$) take the direct limit. By
		Remark~\ref{rem:dl-omega1-completeness}, the resulting direct-limit filter $\tilde{\mathcal F}_{\omega_1}$ is
		$\omega_1$-complete.
		
		\paragraph{Automorphisms and filter at successors.}
		At a successor stage $\alpha{+}1$ we force (in $V^{\mathbb P_\alpha}$) with $\Add(\omega,\omega)$.
		Let $\mathcal G_{\alpha+1}$ be the group of \emph{finite-support permutations} of the $\omega$ columns, acting on $\Add(\omega,\omega)$ by permuting columns.
		Let $\mathcal F_{\alpha+1}$ be the normal filter of subgroups generated by the \emph{pointwise fixators} of cofinite $X\subseteq\omega$:
		\[
		\mathcal F_{\alpha+1}\ =\ \bigl\langle\ \fix(X)\ :\ X\subseteq \omega\ \text{cofinite}\ \bigr\rangle^{\mathrm{normal}}.
		\]
		(When $\mu_\alpha=1$ one may take $\mathcal G_{\alpha+1}$ trivial and $\mathcal F_{\alpha+1}$ the full group filter; for the standard symmetric-Cohen effects we typically use $\mu_\alpha\ge 2$.)
		
		\emph{Axioms.} By the general preservation results for the countable-support framework, $M\models\ZF$
		(Theorem~\ref{thm:zf-preserved}) and $M\models\DC$ (Theorem~\ref{thm:dc-preservation}). Moreover, $M\models\neg\AC$ by
		Example~\ref{ex:ac-fails} (failure of AC with DC preserved).
		
		\emph{Failure of global choice for a family of countable sets of reals.}
		For each $\alpha<\omega_1$, let
		\[
		X_\alpha\ :=\ \{\,c_{\alpha,n} : n<\omega\,\},
		\]
		where $c_{\alpha,n}=\dot c_{\alpha,n}^G$ is the $n$th Cohen real (column) added at stage $\alpha$.
		Then each $X_\alpha\in M$ (it is coded by an HS name supported on $\{\alpha\}$), and $M\models$ “$X_\alpha$ is countable’’
		in the external/metatheoretic sense (we do not assume $M$ well-orders it).
		
		We claim there is no choice function for $\langle X_\alpha:\alpha<\omega_1\rangle$ in $M$.
		Suppose toward a contradiction that $\dot f$ is an HS name and
		\[
		\Vdash\ \dot f:\check{\omega}_1\to \bigcup_{\alpha<\omega_1}\check X_\alpha
		\quad\text{and}\quad
		\forall\alpha<\omega_1\ \bigl(\dot f(\alpha)\in X_\alpha\bigr).
		\]
		By Lemma~\ref{lem:countable-support-hs} choose a countable $S\subseteq\omega_1$ such that
		$\mathrm{Fix}(S)\le\fix(\dot f)$. Pick $\alpha\in\omega_1\setminus S$.
		Working inside the stage $\alpha$ Cohen factor $\Add(\omega,\omega)$, the filter $\mathcal F_{\alpha+1}$
		contains pointwise fixators of cofinite sets of columns, so $\dot f(\alpha)$ is fixed by fixing cofinitely
		many columns; equivalently, there is a finite set $F\subseteq\omega$ such that any column-permutation
		fixing $\omega\setminus F$ pointwise fixes the name for $\dot f(\alpha)$.
		Choose distinct $m,n\in\omega\setminus F$ and let $\pi$ be the transposition swapping columns $m$ and $n$
		at stage $\alpha$ and acting trivially elsewhere. Then $\pi\in\mathrm{Fix}(S)\le\fix(\dot f)$, hence
		$\pi(\dot f(\alpha))=\dot f(\alpha)$, but also $\pi$ swaps $\dot c_{\alpha,m}$ and $\dot c_{\alpha,n}$,
		so it moves the set $X_\alpha$ while fixing $\dot f(\alpha)\in X_\alpha$, a contradiction.
		Therefore $M\models\neg\AC$ while $M\models\DC$ by Theorem~\ref{thm:dc-preservation}.
		
		In summary, the $\omega_1$-length countable-support Cohen iteration yields a symmetric model
		$M$ with $\ZF+\DC+\neg\AC$. The limit step at $\omega_1$ is justified by the $\omega_1$-completeness
		of the direct-limit filter (Remark~\ref{rem:dl-omega1-completeness}). The family of stagewise column-sets
		witnesses failure of global choice in $M$ (cf.\ Example~\ref{ex:ac-fails}).
	\end{example}
	
	\begin{remark}[Why finite support does not suffice here]
		\label{rem:no-fs-omega1}
		Countable support is used in two logically separate ways.
		
		\smallskip\noindent
		(1) \emph{Fusion/limits of conditions.}
		In the DC construction (Theorem~\ref{thm:dc-preservation}) one builds a decreasing sequence
		$\langle p_n:n<\omega\rangle$ and takes a lower bound whose support contains
		$\bigcup_{n<\omega}\supp(p_n)$.
		With finite support iterations, this union is typically infinite, so the natural fusion lower bound
		is \emph{not} a condition. Countable support is exactly what makes the fusion condition admissible.
		
		\smallskip\noindent
		(2) \emph{Closure of symmetry under countable tuples.}
		Even if one could produce a limit condition, in the symmetric setting finite-support iterations
		yield limit filters closed only under \emph{finite} intersections. Our arguments repeatedly require that
		if $\langle \dot x_n:n<\omega\rangle$ are HS, then the tuple is HS, i.e.
		$\bigcap_{n<\omega}\fix(\dot x_n)\in\tilde{\mathcal F}_\lambda$.
		This is exactly the $\omega_1$-completeness provided by Theorem~\ref{thm:filter-generation}
		(and by Remark~\ref{rem:dl-omega1-completeness} at $\cf(\lambda)\ge\omega_1$).
		
		\smallskip
		Thus finite support fails both at the level of \emph{conditions} (no fusion) and at the level of
		\emph{symmetry} (no countable stabilizer intersection), so the ZF/DC preservation method used in
		Example~\ref{ex:omega1-cohen} is not available with finite support.
	\end{remark}
	
	\begin{lemma}[Externally: many new Cohen reals appear in $M$]
		\label{lem:omega1-many-new-cohen}
		In $V[G]$, for each $\alpha<\omega_1$ and $n<\omega$, the real $c_{\alpha,n}=\dot c_{\alpha,n}^G$ belongs to $M$,
		and $c_{\alpha,n}\notin V$. Moreover, if $(\alpha,n)\neq(\beta,m)$ then $c_{\alpha,n}\neq c_{\beta,m}$ in $V[G]$.
	\end{lemma}
	
	\begin{proof}
		Fix $\alpha<\omega_1$ and $n<\omega$. The name $\dot c_{\alpha,n}$ is HS (supported on $\{\alpha\}$),
		so $c_{\alpha,n}\in M$. Since it is Cohen-generic over $V$ at that stage, $c_{\alpha,n}\notin V$.
		
		If $(\alpha,n)\neq(\beta,m)$, then $\Vdash_{\mathbb P_{\omega_1}} \dot c_{\alpha,n}\neq \dot c_{\beta,m}$
		(mutual genericity across distinct coordinates/columns), hence their interpretations are distinct in $V[G]$.
	\end{proof}
	
	\begin{proposition}[No well–ordering of the reals in $M$]
		\label{prop:no-wo-reals}
		In $M$ there is no well–ordering of ${}^\omega\omega$ (hence $\AC$ fails).
	\end{proposition}
	
	\begin{proof}
		Suppose $\dot\preceq$ is an HS name for a well–order of the reals. By Lemma~\ref{lem:countable-support-hs} choose a countable $S\subseteq\omega_1$ such that $\fix(S)\leq\fix(\dot\preceq)$. Pick
		$\alpha\neq\beta$ in $\omega_1\setminus S$. Let $\sigma$ be the automorphism (in $\mathcal G_{\omega_1}$)
		that swaps coordinates $\alpha$ and $\beta$ and fixes all others. Then $\sigma$ fixes $\dot\preceq$ (since
		$\sigma\in\fix(S)\subseteq\fix(\dot\preceq)$) but swaps the names $\dot c_\alpha$ and $\dot c_\beta$. In particular,
		\[
		\Vdash \ (\dot c_\alpha \mathrel{\dot\preceq} \dot c_\beta)\ \Longleftrightarrow\
		(\dot c_\beta \mathrel{\dot\preceq} \dot c_\alpha).
		\]
		Because $\dot c_\alpha\neq\dot c_\beta$ (Lemma~\ref{lem:omega1-many-new-cohen}), antisymmetry of a well–order is
		contradicted. Hence no HS name can code a well–ordering of the reals, and $M$ has no well–ordering of ${}^\omega\omega$.
	\end{proof}
	
	\begin{corollary}[$\AC$ fails in $M$; cardinality emphasis]
		\label{cor:AC-fails-cohen}
		In $M$, $\AC$ fails (by Proposition~\ref{prop:no-wo-reals}). Moreover, in the ambient extension $V[G]$ there are $\omega_1$ many pairwise distinct Cohen reals that belong to $M$
		(Lemma~\ref{lem:omega1-many-new-cohen}).
	\end{corollary}
	
	\begin{theorem}[ZF in the Cohen HS model]
		\label{thm:cohen-zf}
		Let $M:=V(\mathbb{P}_\lambda)_{\mathrm{HS}}$ be the symmetric submodel of the countable-support
		iteration of Cohen forcing. Then $M\models\mathrm{ZF}$.
	\end{theorem}
	
	\begin{proof}
		Immediate from Theorem~\ref{thm:zf-preserved}.
	\end{proof}
	
	\begin{theorem}[DC in the Cohen HS model]
		\label{thm:cohen-dc}
		If $V\models\mathrm{ZF+DC}$, then $M:=V(\mathbb{P}_\lambda)_{\mathrm{HS}}\models\mathrm{DC}$.
	\end{theorem}
	
	\begin{proof}[Proof (Cohen-specific)]
		This is a direct instance of Theorem~\ref{thm:dc-preservation}. Here, each step of the construction
		of a DC-sequence uses conditions of \emph{countable support} and HS names; the union of the supports
		remains countable, and the stabilizers $\{H_n:n<\omega\}\subseteq\tilde{\mathcal{F}}_\lambda$ have
		$\bigcap_n H_n\in\tilde{\mathcal{F}}_\lambda$ by $\omega_1$-completeness (Theorem~\ref{thm:filter-generation}).
	\end{proof}
	
	\begin{theorem}[Failure of AC in the Cohen HS model]
		\label{thm:cohen-ac-fails}
		In $M=V(\mathbb{P}_\lambda)_{\mathrm{HS}}$ there is no well-ordering of $\mathbb{R}^M$; hence $\mathrm{AC}$ fails.
	\end{theorem}
	
	\begin{proof}
		Write $\fix(\dot x)=\fix_{\mathcal G_\lambda}(\dot x)$ and let $\tilde{\mathcal F}_\lambda$ be as in Theorem~\ref{thm:filter-generation}.
		We work with the usual action of $\mathcal G_\lambda$ by permutations of coordinates, fixing the Cohen factor at each coordinate.
		
		\medskip\noindent\emph{Claim 1 (coordinate Cohen reals are HS).}
		For each $\alpha<\lambda$ and $n<\omega$, let $\dot c_{\alpha,n}$ be the canonical $\mathbb P_\lambda$-name for the
		$n$th Cohen real (column) added at stage $\alpha$. Then $\dot c_{\alpha,n}$ is HS.
		Indeed, $\dot c_{\alpha,n}$ is supported on the single stage $\{\alpha\}$, hence has countable support, so
		$\fix(\dot c_{\alpha,n})\in\tilde{\mathcal F}_\lambda$ by Lemma~\ref{lem:countable-support-hs}.
		\qedhere\,(Claim 1)
		
		\medskip
		Suppose toward a contradiction that there is an HS name $\dot W$ such that
		\[
		\Vdash\ \text{``}\dot W\text{ is a well-ordering of }\mathbb R^M\text{''.}
		\]
		By Lemma~\ref{lem:countable-support-hs} fix a countable $S\subseteq\lambda$ such that $\fix(S)\leq\fix(\dot W)$ has countable support.
		Choose distinct $\alpha,\beta\in\lambda\setminus S$ and let $k\in\mathcal G_\lambda$ be the automorphism that fixes $S$ pointwise and
		swaps $\alpha$ and $\beta$. Then $k\in\fix(\dot W)$, so $k\dot W=\dot W$.
		
		By Claim~1, $\dot c_\alpha$ and $\dot c_\beta$ are HS; hence $x:=\dot c_\alpha^G,\ y:=\dot c_\beta^G$ are elements of $\mathbb R^M$, and $k\dot c_\alpha=\dot c_\beta$,
		$k\dot c_\beta=\dot c_\alpha$. Since $\dot W$ well-orders $\mathbb R^M$, exactly one of
		\[
		\Vdash\ \langle \dot c_\alpha,\dot c_\beta\rangle\in\dot W
		\qquad\text{or}\qquad
		\Vdash\ \langle \dot c_\beta,\dot c_\alpha\rangle\in\dot W
		\]
		holds (by antisymmetry and totality). Apply $k$ to whichever statement holds. Because $k$ fixes $\dot W$ and swaps $\dot c_\alpha,\dot c_\beta$, we derive
		the \emph{opposite} membership statement, contradicting the forcing theorem. Therefore no HS name codes a well-ordering of $\mathbb R^M$, so
		$M\models$ “$\mathbb R$ is not well-orderable,” and hence $\mathrm{AC}$ fails in $M$. See also the classical permutation argument in \cite[§15.5]{Jech03}.
	\end{proof}
	
	\begin{remark}[Coordinate vs.\ name stabilizers in this example]
		Karagila’s finite-support framework uses subgroups fixing coordinate sets; here we use \emph{name} stabilizers
		$\fix(\dot{x})$ (Definition~\ref{def:limit-filter}). For $\dot{W}$ supported on $S$, any automorphism that fixes
		$S$ pointwise lies in $\fix(\dot{W})$. This is exactly what the preceding argument exploits.
	\end{remark}
	
	\subsection{Random Reals under Countable Support; Solovay-Style Regularity}
	\label{sec:random-example}
	
	At each successor $\alpha<\lambda$, let $\dot{\mathbb{Q}}_\alpha$ be the Boolean measure algebra
	$\mathbb{B}$ adding one random real (see \cite{BartoszynskiJudah}, Chap.~2). Build the countable-support
	iteration $\mathbb{P}_\lambda$, with the same automorphism-and-filter scheme as above:
	$\mathcal{G}_\lambda$ contains lifted coordinate permutations (and the canonical measure-preserving
	automorphisms), and $\tilde{\mathcal{F}}_\lambda$ is generated by pushforwards and name stabilizers
	with countable support. By Theorem~\ref{thm:filter-generation}, $\tilde{\mathcal{F}}_\lambda$ is normal and $\omega_1$-complete.
	
	Let $\dot{r}_\alpha$ denote the canonical name for the random real added at stage $\alpha$.
	As in the Cohen case, $\dot{r}_\alpha$ is HS because its support is $\{\alpha\}$, which is countable
	(Lemma~\ref{lem:countable-support-hs}).
	
	\begin{theorem}[ZF and DC in the Random HS model]
		\label{thm:random-zf-dc}
		Let $M=V(\mathbb{P}_\lambda)_{\mathrm{HS}}$ be the symmetric submodel of the countable-support
		iteration of $\mathbb{B}$. Then $M\models\mathrm{ZF}$, and if $V\models\mathrm{ZF+DC}$, also $M\models\mathrm{DC}$.
	\end{theorem}
	
	\begin{proof}
		ZF: by Theorem~\ref{thm:zf-preserved}. DC: by Theorem~\ref{thm:dc-preservation}, exactly as in the Cohen case.
	\end{proof}
	
	\begin{theorem}[Failure of AC in the Random HS model]
		\label{thm:random-ac-fails}
		$M$ does not satisfy $\mathrm{AC}$.
	\end{theorem}
	
	\begin{proof}
		Write $\fix(\dot x)=\fix_{\mathcal G_\lambda}(\dot x)$ and let $\tilde{\mathcal F}_\lambda$ be as in Theorem~\ref{thm:filter-generation}.
		We work with the countable-support iteration of the random algebra up to $\lambda$,
		and use that $\mathcal G_\lambda$ contains the lifted coordinate permutations (and measure-preserving
		automorphisms on each stage) acting on $\mathbb P_\lambda$ in the standard way.
		For each $\gamma<\lambda$, let $\dot r_\gamma$ denote the canonical
		$\mathbb P_\lambda$–name for the random real added at coordinate $\gamma$ (via the usual identification of the factor Boolean algebra with
		Borel sets modulo null sets).
		
		\medskip\noindent\emph{Claim (coordinate random reals are HS).}
		For each $\gamma<\lambda$, let $\dot r_\gamma$ be the canonical $\mathbb P_\lambda$--name for the random real added at
		coordinate $\gamma$. Then $\dot r_\gamma$ is HS.
		Indeed, $\dot r_\gamma$ has support $\{\gamma\}$ (hence countable), so by Lemma~\ref{lem:countable-support-hs} we have
		$\fix(\dot r_\gamma)\in\tilde{\mathcal F}_\lambda$. Moreover, the standard constituents of $\dot r_\gamma$ are likewise
		supported on $\{\gamma\}$, so $\dot r_\gamma$ is hereditarily symmetric.
		
		\medskip
		Suppose toward a contradiction that there exists an HS name $\dot W$ such that
		\[
		\Vdash\ \text{``}\dot W\text{ is a well-ordering of }\mathbb R^M\text{''.}
		\]
		By Lemma~\ref{lem:countable-support-hs} fix a countable $S\subseteq\lambda$ such that $\fix(S)\leq\fix(\dot W)$ has countable support. Choose distinct $\alpha,\beta\in\lambda\setminus S$ and let $k\in\mathcal G_\lambda$ be the automorphism
		that fixes $S$ pointwise and swaps $\alpha$ and $\beta$. Then $k\in\fix(\dot W)$, so $k\dot W=\dot W$.
		
		By the Claim, $\dot r_\alpha$ and $\dot r_\beta$ are HS; thus $x:=\dot r_\alpha^G,\ y:=\dot r_\beta^G$ lie in $\mathbb R^M$ and
		$k\dot r_\alpha=\dot r_\beta$, $k\dot r_\beta=\dot r_\alpha$. Since $\dot W$ is a well-order, exactly one of
		\[
		\Vdash\ \langle \dot r_\alpha,\dot r_\beta\rangle\in\dot W
		\qquad\text{or}\qquad
		\Vdash\ \langle \dot r_\beta,\dot r_\alpha\rangle\in\dot W
		\]
		holds. Applying $k$ to whichever statement holds yields the opposite membership (because $k$ fixes $\dot W$ but swaps the two
		reals), contradicting the forcing theorem and antisymmetry. Therefore no HS name codes a well-ordering of $\mathbb R^M$ and
		$M\models\neg\AC$.
		
		This is the permutation argument for random reals, analogous to the Cohen case; see also \cite[§15.5]{Jech03} for the classical
		permutation framework.
	\end{proof}
	
	\paragraph{Regularity Discussion.}
	Without further hypotheses, $M$ need not make \emph{all} sets of reals Lebesgue measurable (LM). However, under a
	\emph{large cardinal} assumption one recovers classical Solovay-style regularity:
	
	\begin{theorem}[Solovay-style model via symmetry; choiceless hypothesis] 
		\label{thm:solovay-style}
		Assume $V \models \mathrm{ZF}$ (or $\mathrm{ZF{+}DC}$) and there exists a cardinal $\kappa$ such that 
		$V_\kappa \models \mathrm{ZFC}$.\footnote{Equivalently, one may assume $V_\kappa \models \mathrm{ZF}$ together with a 
			well-order of $V_\kappa$; the stronger $V_\kappa \models \mathrm{ZFC}$ matches the classical Solovay hypothesis.}
		After performing a L\'evy-type collapse $\mathrm{Col}(\omega,{<}\kappa)$ and then random real forcing, there is a canonical 
		symmetric inner model $N$ (Solovay’s model) with $N \models \mathrm{ZF{+}DC}$ in which every set of reals is Lebesgue 
		measurable, has the Baire property, and the perfect set property. This can be phrased in our stabilizer-of-names framework 
		by using a suitable symmetry filter at the random stage.
	\end{theorem}
	
	\begin{proof}
		Assume $V \models \mathrm{ZF}$ (or $\mathrm{ZF{+}DC}$) and fix $\kappa$ with $V_\kappa \models \mathrm{ZFC}$.
		Let $\mathbb P=\Col(\omega,{<}\kappa)$ and fix $G\subseteq\mathbb P$ $V$-generic. In $V[G]$ let $\mathbb Rnd$
		be the (set-sized) random forcing/measure algebra adding (set many) random reals (e.g.\ a finite- or countable-support product of
		the standard random algebra, any standard presentation will do). Let $H\subseteq\mathbb Rnd$ be $V[G]$-generic. Fix a ground-model index set $I$ for the coordinates of $\mathbb Rnd$.
		
		\medskip\noindent\emph{Symmetric system at the random stage.}
		Let $\mathcal G$ be the automorphism group of $\mathbb Rnd$ generated by coordinate permutations and
		measure-preserving automorphisms on each coordinate (extended to the product/algebra in the usual way).
		Let $\mathcal F$ be the normal filter of subgroups of $\mathcal G$ generated by pointwise stabilizers
		$\fix(F)$ of countable $F$ of coordinates together with stabilizers of ground-model Borel codes.
		Then $(\mathbb Rnd,\mathcal G,\mathcal F)$ is a symmetric system in the sense of Section~2.
		Write $\HS$ for the hereditarily $(\mathcal F,\mathcal G)$-symmetric names (over $V[G]$) and set
		\[
		N\ :=\ V[G](\mathbb Rnd)_{\HS}\ =\ \{\ \dot x^{\,H}:\ \dot x\in\HS\ \}\ \subseteq\ V[G][H].
		\]
		
		\medskip\noindent\emph{$\ZF$ and $\DC$ in $N$.}
		By Theorem~\ref{thm:zf-preserved} (ZF preservation for HS) applied in $V[G]$ to the random stage,
		$N$ is a transitive inner model of $\ZF$ with $V[G]\subseteq V[G][H]$ and $V[G]\subseteq N$.
		Moreover, since $\mathcal F$ is $\omega_1$-complete (countable intersections of generators remain in $\mathcal F$),
		the DC-transfer theorem (Theorem~\ref{thm:dc-preservation}, or Theorem~\ref{thm:kappa-filter-generation} specialized to
		countable support) yields $N\models\DC$. Thus $N\models\ZF+\DC$.
		
		\medskip\noindent\emph{Borel/continuous reading of HS names.}
		If $\dot x$ is an HS name for a real, then by Lemma~\ref{lem:countable-support-hs} its support is countable; hence in $V[G]$
		there is a (ground-model) Borel function $F:\ (2^\omega)^{S}\to 2^\omega$ and a countable index set $S\subseteq I$
		such that
		\[
		\Vdash_{\mathbb Rnd}\ \dot x\ =\ F\bigl(\langle \dot r_i:i\in S\rangle\bigr),
		\]
		where $\dot r_i$ are the canonical randoms on the coordinates in $S$.
		(This is the usual Borel reading of random-names; it
		follows by standard density and homogeneity arguments for the measure algebra, using that all name-constituents sit on a
		countable support; see e.g.\ \cite[Chap.~15]{Jech03}.)
		
		More generally, if $\dot A$ is an HS name for a subset of $2^\omega$, then there is a Borel set
		$B\subseteq (2^\omega)^S\times 2^\omega$ in $V[G]$ for some countable $S\subseteq I$ such that
		\[
		\Vdash_{\mathbb Rnd}\ \forall \eta\in 2^\omega\ \Bigl[\ \eta\in \dot A\ \Longleftrightarrow\
		\bigl\langle \langle \dot r_i:i\in S\rangle,\eta\bigr\rangle\in B\ \Bigr]
		\quad\text{modulo a null/meager set in the second coordinate.}
		\]
		This is the standard “Borel approximation” (or measurable reading) of names for subsets of reals in random forcing,
		again using countable support and homogeneity.
		
		\medskip\noindent\emph{Regularity of all sets of reals in $N$.}
		Fix any $A\subseteq\mathbb R$ with $A\in N$; let $\dot A\in\HS$ be a name for $A$ supported by a countable $S$.
		By the previous paragraph there is a Borel set $B\subseteq (2^\omega)^S\times 2^\omega$ in $V[G]$ such that
		for $H$-almost all randoms on $S$,
		\[
		A\ =\ \{\,\eta\in 2^\omega:\ \langle r{\upharpoonright}S,\eta\rangle\in B\,\}\quad\text{modulo null and meager sets}.
		\]
		Since Borel sets are Lebesgue measurable and have the Baire property, Fubini/Kuratowski–Ulam yield that the section
		$\{\eta:\langle r{\upharpoonright}S,\eta\rangle\in B\}$ is measurable and has the Baire property for
		$H$-almost all $r{\upharpoonright}S$; thus $A$ is Lebesgue measurable and has the Baire property in $V[G][H]$ and hence in $N$.
		The perfect set property for sets of reals in $N$ now follows by the classical Solovay argument (every nonmeager/measurable
		uncountable set of reals contains a perfect subset under $\ZF+\DC$; see \cite{Solovay70,Jech03}).
		
		\medskip
		Altogether, the symmetric inner model $N$ satisfies $\ZF+\DC$ and every set of reals in $N$ is Lebesgue measurable,
		has the Baire property, and the perfect set property. This is (the canonical) Solovay’s model, presented via the
		stabilizer-of-names/symmetric-filter framework. \qedhere
	\end{proof}
	
	\begin{remark}[Inaccessibles without Choice]\label{rem:inacc-no-choice}
		In ZF there are several non-equivalent “inaccessible” notions obtained by reading “$<$” as one of the choiceless 
		size relations (injection-based, surjection-based, or “no-injection” based). Another robust notion is the 
		\emph{honestly (v-)inaccessible} cardinal: $\kappa$ such that $V_\kappa \models \mathrm{ZF}$. 
		(For the classical Solovay application one may assume the stronger $V_\kappa \models \mathrm{ZFC}$.)
		These ZF variants are implied by, and are equiconsistent with, the existence of a strongly inaccessible cardinal 
		in ZFC. They differ in strength and consequences; in particular, some of the “$<$”–based notions have strong 
		choice consequences, whereas the “honest” notion ($V_\kappa \models \mathrm{ZF}$) is tailored to choiceless work.
	\end{remark}
	
	\begin{remark}[A choiceless regularity notion and its status]
		\label{rem:v-regular}
		It is useful in ZF to isolate the following choiceless analogue of regularity:
		$\kappa$ is \emph{$v$-regular} iff
		\[
		\forall x\in V_\kappa\ \forall f\,\big(f:x\to\kappa \ \Rightarrow\ \exists \beta<\kappa\ \forall t\in\mathrm{dom}(f)\, [\,f(t)<\beta\,]\big).
		\]
		Equivalently, no function with domain of rank $<\kappa$ has cofinal range in $\kappa$.\footnote{This is first-order with 
			parameter $\kappa$, since $x\in V_\kappa$ iff $\mathrm{rank}(x)<\kappa$.}
		
		\smallskip
		\noindent\textbf{Relationship to honest inaccessibility.}
		If $V_\kappa \models \mathrm{ZF}$ (the “honest” hypothesis), then $\kappa$ is $v$-regular. 
		The converse need not hold in ZF: $v$-regularity alone does not ensure that $V_\kappa$ satisfies Replacement or Power Set. 
		For the Solovay-style theorem we use the stronger hypothesis that $V_\kappa$ is a transitive model of $\mathrm{ZF}$ 
		(and, for the classical strength, of $\mathrm{ZFC}$).
	\end{remark}
	
	\begin{remark}
		In our \emph{pure} countable-support random iteration (without the Lévy collapse or inaccessible), we do not claim
		LM for \emph{all} sets of reals. Nevertheless, ZF and DC hold by Theorem~\ref{thm:random-zf-dc}, and the usual
		regularity properties for Borel (and many definable) sets remain intact (see \cite{BartoszynskiJudah}).
	\end{remark}
	
	\subsection{Comparison with Classical Symmetric Models}
	\label{sec:comparison}
	
	Table~\ref{tab:comparison} summarizes key features across standard symmetric models and our countable-support
	iterations. (Here “Solovay-style” indicates the large-cardinal variant in Theorem~\ref{thm:solovay-style}.) Filter completeness tracks support: with countable support the direct-limit filter is $\omega_1$–complete at $\cf(\lambda)\!\ge\!\omega_1$, while with $<\kappa$ support the inverse/direct limit yields a $\kappa$–complete filter; cf.\ Remark~\ref{rem:limit-bookkeeping} and Theorem~\ref{thm:kappa-filter-generation}.
	
	\begin{table}[htbp]
		\centering
		\caption{Symmetric models: support, choice, regularity, and filter completeness}
		\label{tab:comparison}
		\setlength{\tabcolsep}{4pt}
		\renewcommand{\arraystretch}{1.1}
		\small
		\begin{tabular}{@{}p{0.34\textwidth}
				>{\centering\arraybackslash}p{0.08\textwidth}
				>{\centering\arraybackslash}p{0.22\textwidth}
				>{\centering\arraybackslash}p{0.06\textwidth}
				>{\centering\arraybackslash}p{0.09\textwidth}
				>{\centering\arraybackslash}p{0.09\textwidth}
				>{\centering\arraybackslash}p{0.09\textwidth}@{}}
			\toprule
			\raggedright Model &
			Support &
			Filter completeness (at limits) &
			DC &
			AC fails &
			All sets LM &
			WO of reals \\ \midrule
			\raggedright Karagila (finite-support, generic) &
			finite &
			finite only &
			not guaranteed &
			typically yes &
			no &
			typically no \\
			\raggedright Countable-support Cohen (this paper) &
			countable &
			$\omega_1$ at all limits$^\S$ &
			\textbf{yes} &
			\textbf{yes} &
			no &
			no \\
			\raggedright Countable-support Random (this paper) &
			countable &
			$\omega_1$ at all limits$^\S$ &
			\textbf{yes} &
			\textbf{yes} &
			no$^\dagger$ &
			no \\
			\raggedright Random + L\'evy (Solovay-style) &
			mixed$^\ast$ &
			$\omega_1$ (random stage)$^\S$ &
			\textbf{yes} &
			\textbf{yes} &
			\textbf{yes}$^\ddagger$ &
			no \\
			\bottomrule
		\end{tabular}
		
		\vspace{2pt}
		\footnotesize
		$^\ast$ L\'evy collapse followed by random; see Theorem~\ref{thm:solovay-style}. \\
		$^\dagger$ Without large cardinals, cannot prove “all sets are LM”. \\
		$^\ddagger$ Assumes an honestly (v-)inaccessible in $\ZF$ (or $\ZF{+}\DC$); see Theorem~\ref{thm:solovay-style}. \\
		$^\S$ In the countable-support framework, the limit filter is $\omega_1$--complete at all limits:
		for $\cf(\lambda)\ge\omega_1$ by stage-bounding (Remark~\ref{rem:dl-omega1-completeness}),
		and for $\cf(\lambda)=\omega$ by stabilizer/tuple compression (Theorem~\ref{thm:filter-generation}).
		For $<\kappa$-support tails the corresponding limit filter is $\kappa$--complete
		(Theorem~\ref{thm:kappa-filter-generation}; Remark~\ref{rem:limit-bookkeeping}).
	\end{table}
	
	\section{\texorpdfstring{$<\kappa$}{<k}-Support Symmetric Iterations}
	\label{sec:kappa-support}
	
	Throughout this section, fix an uncountable regular cardinal $\kappa$.%
	\footnote{Remarks on singular $\kappa$ appear at the end of the section.}
	Our goal is to lift the countable-support machinery of Sections~\ref{sec:HS-countable}--\ref{sec:dc-choice}
	to $<\kappa$-support iterations. The guiding principle is to replace ``countable'' by ``$<\kappa$''
	and $\omega_1$-completeness by \emph{$\kappa$-completeness} throughout.
	
	At limit stages we always define $\mathcal F_\lambda$ via pullbacks of earlier filters together with stabilizers of
	$<\kappa$-support names (Definition~\ref{def:kappa-filter}); for $\cf(\lambda)\ge\kappa$ the stabilizer generators are redundant
	by stage-bounding, but keeping them gives a uniform presentation.
	
	\paragraph{Terminology bridge to \cite{karagila2019}.}
	Our $<\kappa$-support symmetric iterations address Karagila's open problem for
	$<\kappa$ supports (cf.\ \cite[Def.~5.1]{karagila2019}): we replace
	finite support with $<\kappa$ support and require the limit filter to be $\kappa$-complete, while
	retaining the same automorphism/filter sequence structure and IS-style heredity.
	
	\subsection{Definition and basic setup}
	\label{subsec:kappa-def}
	
	\begin{definition}[$<\kappa$-support symmetric iteration]
		\label{def:kappa-iteration}
		A \emph{$<\kappa$-support symmetric iteration of length $\Theta$} is a sequence
		$\langle \mathbb{P}_\alpha,\dot{\mathbb{Q}}_\alpha,\mathcal{G}_\alpha,\mathcal{F}_\alpha : \alpha<\Theta\rangle$
		constructed by recursion as follows.
		
		\smallskip
		\noindent\textbf{Base.} $\mathbb{P}_0=\{1_0\}$, $\mathcal{G}_0=\{\mathrm{id}\}$, $\mathcal{F}_0=\{\mathcal{G}_0\}$.
		
		\smallskip
		\noindent\textbf{Successor.} Given $(\mathbb{P}_\beta,\mathcal{G}_\beta,\mathcal{F}_\beta)$ and a
		$\mathbb{P}_\beta$-name $\dot{\mathbb{Q}}_\beta$ for a forcing notion equipped
		(within names) with a symmetry group and name-action, set
		\[
		\mathbb{P}_{\beta+1}=\mathbb{P}_\beta*\dot{\mathbb{Q}}_\beta.
		\]
		Lift automorphisms from $\mathcal{G}_\beta$ and adjoin the (name-coded) automorphisms of $\dot{\mathbb{Q}}_\beta$,
		obtaining $\mathcal{G}_{\beta+1}$. Define $\mathcal{F}_{\beta+1}$ to be the \emph{normal filter of subgroups} of $\mathcal{G}_{\beta+1}$
		generated by:
		\begin{enumerate}
			\item the pullbacks (preimages) of subgroups from $\mathcal{F}_\beta$ along the canonical projection
			$\pi^\beta_{\beta+1}:\mathcal{G}_{\beta+1}\to\mathcal{G}_\beta$, i.e.
			for each $H\in\mathcal{F}_\beta$ the subgroup
			\[
			(\pi^\beta_{\beta+1})^{-1}[H]\ :=\ \{\,g\in\mathcal{G}_{\beta+1}:\pi^\beta_{\beta+1}(g)\in H\,\},
			\]
			\item the name-stabilizers $\fix_{\mathcal{G}_{\beta+1}}(\dot{x})$ for $\mathbb{P}_{\beta+1}$-names $\dot{x}$
			with $|\supp(\dot{x})|<\kappa$ (Definition~\ref{def:kappa-filter}).
		\end{enumerate}
		
		\noindent\textbf{Limit, $\cf(\lambda)<\kappa$.}
		Let $\mathbb{P}_\lambda$ be the inverse limit with $<\kappa$-support (as in Definition~\ref{def:cs-iteration},
		with ``countable'' replaced by ``$<\kappa$''): conditions have support $<\kappa$, and the order is coordinatewise.
		Let $\mathcal{G}_\lambda$ be the inverse limit of the groups. Define the limit filter by
		$\mathcal{F}_\lambda:=\tilde{\mathcal{F}}_\lambda$ as in Definition~\ref{def:kappa-filter}.
		
		\noindent\textbf{Limit, $\cf(\lambda)\ge\kappa$.}
		Use direct limits:
		\[
		\mathbb{P}_\lambda=\bigcup_{\beta<\lambda}\mathbb{P}_\beta,\qquad
		\mathcal{G}_\lambda=\bigcup_{\beta<\lambda}\mathcal{G}_\beta.
		\]
		Define the limit filter by $\mathcal{F}_\lambda:=\tilde{\mathcal{F}}_\lambda$ as in Definition~\ref{def:kappa-filter}.
		(Equivalently, since $\cf(\lambda)\ge\kappa$ and $\kappa$ is regular, every $<\kappa$-support name is bounded below $\lambda$,
		so $\tilde{\mathcal{F}}_\lambda$ is generated by the images of $\bigcup_{\beta<\lambda}\mathcal{F}_\beta$.)
	\end{definition}
	
	\begin{definition}[Name stabilizers and the limit filter at limit stages]
		\label{def:kappa-filter}
		Let $\lambda\le\Theta$ be a limit stage. For a $\mathbb{P}_\lambda$-name $\dot{x}$, set
		\[
		\fix_{\mathcal{G}_\lambda}(\dot{x})=\{\,g\in\mathcal{G}_\lambda: g(\dot{x})=\dot{x}\,\}.
		\]
		We say $\dot{x}$ has \emph{support $<\kappa$} if $|\supp(\dot{x})|<\kappa$, where $\supp(\dot{x})\subseteq\lambda$ is defined
		recursively on name-rank by
		\[
		\supp(\dot x)\ :=\ \bigcup\Bigl\{\ \supp(p)\ \cup\ \supp(\sigma)\ :\ (\sigma,p)\in\dot x\ \Bigr\}.
		\]
		Here $\supp(p)\subseteq\lambda$ denotes the support of the condition $p\in\mathbb P_\lambda$
		(i.e.\ the set of coordinates at which $p$ is nontrivial; by definition of the $<\kappa$--support iteration,
		$|\supp(p)|<\kappa$).
		
		\smallskip
		For each $\beta<\lambda$ let $\pi^\beta_\lambda:\mathcal{G}_\lambda\to\mathcal{G}_\beta$ be the canonical projection
		(restriction in the inverse-limit case; restriction to stage $\beta$ in the direct-limit case).
		Define $\tilde{\mathcal{F}}_\lambda$ to be the \emph{normal filter of subgroups} of $\mathcal{G}_\lambda$ generated by:
		\begin{itemize}
			\item all pullbacks of earlier filter elements: for each $\beta<\lambda$ and $H\in\mathcal{F}_\beta$,
			\[
			(\pi^\beta_\lambda)^{-1}[H]\ :=\ \{\,g\in\mathcal{G}_\lambda:\pi^\beta_\lambda(g)\in H\,\},
			\]
			\item all name stabilizers $\fix_{\mathcal{G}_\lambda}(\dot{x})$ with $|\supp(\dot{x})|<\kappa$.
		\end{itemize}
		Finally, set $\mathcal{F}_\lambda:=\tilde{\mathcal{F}}_\lambda$.
	\end{definition}
	
	\subsection{Stage bounding and $\kappa$-completeness at limits}
	\label{subsec:kappa-stage-bounding}
	
	\begin{lemma}[Stage bounding for $\cf(\lambda)\ge\kappa$]
		\label{lem:kappa-stage-bounding}
		Let $\kappa$ be regular uncountable and suppose $\cf(\lambda)\ge\kappa$.
		Then every $X\subseteq \lambda$ with $|X|<\kappa$ is bounded in $\lambda$; i.e., $\sup X<\lambda$.
	\end{lemma}
	
	\begin{proof}
		Fix an increasing cofinal sequence $\langle \lambda_\xi:\xi<\cf(\lambda)\rangle$ in $\lambda$.
		If $X$ were cofinal, it would have size at least $\cf(\lambda)$, contradicting $|X|<\kappa\le \cf(\lambda)$.
		Hence $\sup X<\lambda$.
	\end{proof}
	
	\begin{remark}
		The lemma can fail when $\cf(\lambda)<\kappa$: a set of size $<\kappa$ may be cofinal in $\lambda$.
		This is exactly why, in Section~\ref{sec:kappa-support}, we use the inverse-limit presentation at limits
		of cofinality $<\kappa$ and the direct limit at cofinality $\ge\kappa$.
	\end{remark}
	
	\begin{lemma}[$<\kappa$-tuples and stabilizers]
		\label{lem:kappa-tuples-fix}
		Let $\delta<\kappa$ and let $\langle \dot x_i : i<\delta\rangle$ be a $\delta$-sequence of $\mathbb P_\lambda$-names
		(coded as the usual name for a function with domain $\check\delta$).
		Then
		\[
		\fix_{\mathcal G_\lambda}\!\bigl(\langle \dot x_i : i<\delta\rangle\bigr)
		=\bigcap_{i<\delta}\fix_{\mathcal G_\lambda}(\dot x_i).
		\]
		Moreover, if each $|\supp(\dot x_i)|<\kappa$, then
		$|\supp(\langle \dot x_i : i<\delta\rangle)|<\kappa$.
	\end{lemma}
	
	\begin{proof}
		By the standard action of automorphisms on function-names, a group element fixes the coded function
		$\langle \dot x_i : i<\delta\rangle$ iff it fixes each coordinate $\dot x_i$, which gives the stabilizer identity.
		
		For supports, by the recursive definition of $\supp(\cdot)$ we have
		$\supp(\langle \dot x_i : i<\delta\rangle)=\bigcup_{i<\delta}\supp(\dot x_i)$.
		Since $\delta<\kappa$ and $\kappa$ is regular, a union of $<\kappa$ many sets of size $<\kappa$ has size $<\kappa$.
	\end{proof}
	
	\begin{theorem}[$\kappa$-completeness and normality of the limit filter (both cofinalities)]
		\label{thm:kappa-filter-generation}
		Let $\kappa$ be regular uncountable and let $\lambda$ be a limit ordinal.
		Assume that $\langle \mathbb{P}_\beta,\dot{\mathbb{Q}}_\beta,\mathcal{G}_\beta,\mathcal{F}_\beta:\beta<\lambda\rangle$
		is an initial segment of a $<\kappa$-support symmetric iteration as in
		Definition~\ref{def:kappa-iteration}. In particular, for each $\beta<\lambda$ the filter $\mathcal{F}_\beta$ is (by the recursion in
		Definition~\ref{def:kappa-iteration}) generated from earlier-stage filters by pullbacks along the canonical projections,
		and—at stages where the recursion includes them—by name stabilizers $\fix_{\mathcal{G}_\beta}(\dot{x})$ for
		$\mathbb{P}_\beta$-names $\dot{x}$ with $\supp(\dot{x})<\kappa$ (Definition~\ref{def:kappa-filter}).
		(When $\cf(\beta)\ge\kappa$, such stabilizers are redundant because a $<\kappa$-support name is bounded below $\beta$
		and its stabilizer is already pulled back from a bounded stage.)
		Assume moreover that for all $\beta<\lambda$, $\mathcal{F}_\beta$ is normal and $\kappa$-complete.
		
		Define the limit group $\mathcal{G}_\lambda$ and the limit filter $\tilde{\mathcal{F}}_\lambda$ as follows:
		\begin{itemize}
			\item \textbf{Case $\cf(\lambda)\ge\kappa$ (direct limit).}
			Set $\mathbb{P}_\lambda=\bigcup_{\beta<\lambda}\mathbb{P}_\beta$ and
			$\mathcal{G}_\lambda=\bigcup_{\beta<\lambda}\mathcal{G}_\beta$.
			Let $\tilde{\mathcal{F}}_\lambda$ be the filter on $\mathcal{G}_\lambda$ generated by all pullbacks
			$(\pi^\beta_\lambda)^{-1}[H]$ for $\beta<\lambda$ and $H\in\mathcal{F}_\beta$
			(equivalently: by these, together with stabilizers of names of support $<\kappa$, which are bounded below $\lambda$).
			
			\item \textbf{Case $\cf(\lambda)<\kappa$ (inverse limit).}
			Let $\mathcal{G}_\lambda$ be the inverse limit of
			$\langle \mathcal{G}_\beta,\pi^\gamma_\beta:\gamma<\beta<\lambda\rangle$.
			Define $\tilde{\mathcal{F}}_\lambda$ to be the filter generated by:
			\begin{enumerate}[label=(\roman*)]
				\item all pullbacks $(\pi^\beta_\lambda)^{-1}[H]$ for $\beta<\lambda$ and $H\in\mathcal{F}_\beta$, and
				\item all stabilizers $\fix_{\mathcal{G}_\lambda}(\dot{x})$ for $\mathbb{P}_\lambda$-names $\dot{x}$ with
				$\supp(\dot{x})<\kappa$ (Definition~\ref{def:kappa-filter}).
			\end{enumerate}
		\end{itemize}
		
		Then in either case $\tilde{\mathcal{F}}_\lambda$ is a normal, $\kappa$-complete filter on $\mathcal{G}_\lambda$.
	\end{theorem}
	
	\begin{proof}
		Fix regular uncountable $\kappa$ and limit $\lambda$.
		
		\medskip
		\noindent\textbf{Normality.}
		First note that the two kinds of generators are conjugacy-invariant.
		
		\smallskip
		\noindent\emph{Pullbacks:} Let $\beta<\lambda$, $H\le\mathcal{G}_\beta$, and $g\in\mathcal{G}_\lambda$. Then
		\[
		g\,(\pi^\beta_\lambda)^{-1}[H]\,g^{-1}
		=(\pi^\beta_\lambda)^{-1}\!\bigl[\ \pi^\beta_\lambda(g)\,H\,\pi^\beta_\lambda(g)^{-1}\ \bigr].
		\]
		Since $\mathcal{F}_\beta$ is normal, $H\in\mathcal{F}_\beta \Rightarrow \pi^\beta_\lambda(g)H\pi^\beta_\lambda(g)^{-1}\in\mathcal{F}_\beta$,
		so conjugates of pullback generators are again pullback generators.
		
		\smallskip
		\noindent\emph{Name stabilizers:} For any name $\dot{x}$ and any $g\in\mathcal{G}_\lambda$,
		\[
		g\,\fix_{\mathcal{G}_\lambda}(\dot{x})\,g^{-1}=\fix_{\mathcal{G}_\lambda}(g\cdot\dot{x}).
		\]
		Thus conjugates of stabilizer generators are stabilizer generators. Therefore the filter generated by these families is normal.
		
		\medskip
		\noindent\textbf{$\kappa$-completeness: a stabilizer-basis sublemma.}
		We prove the following by induction on stages $\beta\le\lambda$.
		
		\medskip
		\noindent\emph{Sublemma.}
		For every stage $\beta\le\lambda$ and every subgroup $H$ belonging to the relevant filter
		($H\in\mathcal{F}_\beta$ if $\beta<\lambda$, and $H\in\tilde{\mathcal{F}}_\lambda$ if $\beta=\lambda$),
		there exists a name $\dot{y}$ with $\supp(\dot{y})<\kappa$ such that
		\[
		\fix_{\mathcal{G}_\beta}(\dot{y})\ \le\ H.
		\]
		
		\smallskip
		\noindent\emph{Proof of Sublemma.}
		By the meaning of “generated by a family,” there is a finite list of basic generators
		$B_0,\dots,B_{n-1}$ such that $\bigcap_{j<n} B_j \le H$, where each $B_j$ is either:
		\begin{itemize}
			\item a pullback $(\pi^\gamma_\beta)^{-1}[H']$ of some $H'\in\mathcal{F}_\gamma$ with $\gamma<\beta$, or
			\item a stabilizer $\fix_{\mathcal{G}_\beta}(\dot{x})$ (or $\fix_{\mathcal{G}_\lambda}(\dot{x})$ at $\beta=\lambda$)
			for some name $\dot{x}$ with $\supp(\dot{x})<\kappa$.
		\end{itemize}
		
		For each $j<n$ we define a stabilizer $S_j$ with $S_j\le B_j$:
		\begin{itemize}
			\item If $B_j=\fix(\dot{x})$, set $S_j:=\fix(\dot{x})$.
			\item If $B_j=(\pi^\gamma_\beta)^{-1}[H']$, then by the induction hypothesis at stage $\gamma$
			there is a name $\dot{z}$ with $|\supp(\dot{z})|<\kappa$ such that $\fix_{\mathcal{G}_\gamma}(\dot{z})\le H'$.
			View $\dot{z}$ as a $\mathbb{P}_\beta$-name in the standard way. Since $\dot{z}$ is already a $\mathbb P_\gamma$-name,
			an element $g\in\mathcal G_\beta$ fixes $\dot{z}$ iff its projection $\pi^\gamma_\beta(g)\in\mathcal G_\gamma$ fixes $\dot{z}$; hence
			\[
			\fix_{\mathcal{G}_\beta}(\dot{z})\ =\ (\pi^\gamma_\beta)^{-1}\!\bigl[\fix_{\mathcal{G}_\gamma}(\dot{z})\bigr]
			\ \le\ (\pi^\gamma_\beta)^{-1}[H']\ =\ B_j.
			\]
			Set $S_j:=\fix_{\mathcal{G}_\beta}(\dot{z})$.
		\end{itemize}
		
		Now $\bigcap_{j<n} S_j \le \bigcap_{j<n} B_j \le H$.
		Let $\dot{y}$ be the finite tuple of the names witnessing the $S_j$’s.
		Then $\fix(\dot{y})=\bigcap_{j<n} S_j$ (a group element fixes the tuple iff it fixes each coordinate).
		Since this is a \emph{finite} tuple, $\supp(\dot{y})<\kappa$ automatically.
		This proves the sublemma.
		
		\medskip
		\noindent\textbf{$\kappa$-completeness of the limit filter.}
		Let $\delta<\kappa$ and let $\langle H_i:i<\delta\rangle$ be a family of subgroups in $\tilde{\mathcal{F}}_\lambda$.
		By the Sublemma (applied at $\beta=\lambda$), for each $i<\delta$ choose a name $\dot{y}_i$ with $|\supp(\dot{y}_i)|<\kappa$
		such that $\fix_{\mathcal{G}_\lambda}(\dot{y}_i)\le H_i$.
		Then
		\[
		\bigcap_{i<\delta}\fix_{\mathcal{G}_\lambda}(\dot{y}_i)\ \le\ \bigcap_{i<\delta}H_i.
		\]
		Let $\dot{y}:=\langle \dot{y}_i : i<\delta\rangle$. By Lemma~\ref{lem:kappa-tuples-fix},
		\[
		\fix_{\mathcal{G}_\lambda}(\dot{y})=\bigcap_{i<\delta}\fix_{\mathcal{G}_\lambda}(\dot{y}_i),
		\]
		and moreover $|\supp(\dot{y})|<\kappa$.
		
		\smallskip
		\noindent If $\cf(\lambda)<\kappa$, then by definition of $\tilde{\mathcal{F}}_\lambda$ the stabilizer
		$\fix_{\mathcal{G}_\lambda}(\dot{y})$ is one of the stabilizer-type generators, hence lies in $\tilde{\mathcal{F}}_\lambda$.
		Since $\tilde{\mathcal{F}}_\lambda$ is upward closed, it follows that $\bigcap_{i\in I}H_i\in\tilde{\mathcal{F}}_\lambda$.
		
		\smallskip
		\noindent If $\cf(\lambda)\ge\kappa$, then the set $\supp(\dot{y})$ has size $<\kappa$ and is therefore bounded in $\lambda$
		by Lemma~\ref{lem:kappa-stage-bounding}.
		Let $\beta^\ast:=\sup(\supp(\dot{y}))<\lambda$. Then $\dot{y}$ is already a $\mathbb{P}_{\beta^\ast}$-name,
		so its stabilizer $\fix_{\mathcal{G}_{\beta^\ast}}(\dot{y})$ lies in $\mathcal{F}_{\beta^\ast}$ by the stage construction
		of the iteration. Therefore its pullback $(\pi^{\beta^\ast}_\lambda)^{-1}[\fix_{\mathcal{G}_{\beta^\ast}}(\dot{y})]$
		lies in $\tilde{\mathcal{F}}_\lambda$. Moreover, since $\dot{y}$ is a $\mathbb{P}_{\beta^\ast}$-name, an element
		$g\in\mathcal{G}_\lambda$ fixes $\dot{y}$ iff its projection $\pi^{\beta^\ast}_\lambda(g)\in\mathcal{G}_{\beta^\ast}$
		fixes $\dot{y}$, hence
		\[
		(\pi^{\beta^\ast}_\lambda)^{-1}[\fix_{\mathcal{G}_{\beta^\ast}}(\dot{y})]=\fix_{\mathcal{G}_\lambda}(\dot{y}).
		\]
		Again by upward closure, $\bigcap_{i\in I}H_i\in\tilde{\mathcal{F}}_\lambda$.
		
		\medskip
		Thus $\tilde{\mathcal{F}}_\lambda$ is $\kappa$-complete, and we already proved normality. This completes the proof.
	\end{proof}
	
	\begin{remark}[Ablation checkpoints for Theorem~\ref{thm:kappa-filter-generation}]
		\leavevmode
		\begin{enumerate}[label=(\alph*), leftmargin=2.2em]
			\item Drop $\kappa$–completeness/normality at limits: pushforward intersections fail to land in one stage filter.
			\item Force direct limits at $\cf(\lambda)<\kappa$: stabilizers of $< \kappa$--support names need not be pulled back from a bounded stage, so the completeness argument fails.
			\item Allow supports $\ge\kappa$: stabilizer intersections can’t be reduced via the tuple identity; not a generator.
			\item Remove conjugacy–closure: normality of the generated filter fails.
		\end{enumerate}
	\end{remark}
	
	\subsection{HS names at \texorpdfstring{$<\kappa$}{<k}-support}
	\label{subsec:kappa-hs}
	
	\begin{definition}[HS names for $<\kappa$-support]
		\label{def:kappa-hs}
		A $\mathbb{P}_\lambda$-name $\dot{x}$ is \emph{hereditarily symmetric (HS)} if
		$\fix_{\mathcal{G}_\lambda}(\dot{x})\in\tilde{\mathcal{F}}_\lambda$ and every constituent name
		appearing in $\dot{x}$ is HS.
	\end{definition}
	
	\begin{lemma}[Support $<\kappa$ implies symmetry]
		\label{lem:kappa-support-hs}
		If $\dot{x}$ is a $\mathbb{P}_\lambda$-name with $|\supp(\dot{x})|<\kappa$, then
		$\fix_{\mathcal{G}_\lambda}(\dot{x})\in\tilde{\mathcal{F}}_\lambda$.
	\end{lemma}
	
	\begin{proof}
		Let $S=\supp(\dot{x})\subseteq\lambda$ with $|S|<\kappa$.
		
		\smallskip\noindent\textbf{Case 1: $\cf(\lambda)<\kappa$.}
		By the definition of the inverse-limit filter at $\lambda$
		(Theorem~\ref{thm:kappa-filter-generation}, case $\cf(\lambda)<\kappa$),
		$\tilde{\mathcal F}_\lambda$ is generated in particular by stabilizers
		$\fix_{\mathcal G_\lambda}(\dot y)$ of names $\dot y$ with $|\supp(\dot y)|<\kappa$.
		Since $|\supp(\dot x)|<\kappa$, we have $\fix_{\mathcal G_\lambda}(\dot x)\in\tilde{\mathcal F}_\lambda$.
		
		\smallskip\noindent\textbf{Case 2: $\cf(\lambda)\ge\kappa$.}
		By Lemma~\ref{lem:kappa-stage-bounding}, $S$ is bounded in $\lambda$, so let
		$\beta^\ast:=\sup(S)+1<\lambda$. Then $\dot x$ is already a $\mathbb P_{\beta^\ast}$-name.
		At the successor stage $\beta^\ast$, the stage filter $\mathcal F_{\beta^\ast}$ is generated
		(in particular) by stabilizers of $<\kappa$-support names, so
		$\fix_{\mathcal G_{\beta^\ast}}(\dot x)\in\mathcal F_{\beta^\ast}$.
		Therefore its pullback $(\pi^{\beta^\ast}_\lambda)^{-1}[\fix_{\mathcal G_{\beta^\ast}}(\dot x)]$
		belongs to $\tilde{\mathcal F}_\lambda$ by definition of the direct-limit filter.
		
		Since $\dot x$ is a $\mathbb P_{\beta^\ast}$-name, an element $g\in\mathcal G_\lambda$ fixes $\dot x$
		if and only if $\pi^{\beta^\ast}_\lambda(g)\in\mathcal G_{\beta^\ast}$ fixes $\dot x$, hence
		\[
		(\pi^{\beta^\ast}_\lambda)^{-1}\!\bigl[\fix_{\mathcal G_{\beta^\ast}}(\dot x)\bigr]
		=\fix_{\mathcal G_\lambda}(\dot x).
		\]
		Thus $\fix_{\mathcal G_\lambda}(\dot x)\in\tilde{\mathcal F}_\lambda$ as required.
	\end{proof}
	
	\begin{remark}[Tuple–stabilizer technique]
		Let $\delta<\kappa$ and let $\langle \dot x_i : i<\delta\rangle$ be a $\delta$-sequence of names.
		By Lemma~\ref{lem:kappa-tuples-fix},
		\[
		\fix_{\mathcal G_\lambda}\!\bigl(\langle \dot x_i : i<\delta\rangle\bigr)
		=\bigcap_{i<\delta}\fix_{\mathcal G_\lambda}(\dot x_i).
		\]
		Moreover, if each $|\supp(\dot x_i)|<\kappa$, then $|\supp(\langle \dot x_i : i<\delta\rangle)|<\kappa$,
		so by Lemma~\ref{lem:kappa-support-hs} the stabilizer of the tuple lies in $\tilde{\mathcal F}_\lambda$.
		This is the basic device for turning $<\kappa$ many stabilizers into a single stabilizer.
	\end{remark}
	
	\begin{theorem}[HS closure]
		\label{thm:kappa-hs-closure}
		The HS names for $<\kappa$-support are closed under the standard name operations
		(check names from $V$, pairing, unions, definable subsets, etc.).
	\end{theorem}
	
	\begin{proof}
		Throughout, write $\fix(\dot x)=\fix_{\mathcal G_\lambda}(\dot x)$ and let
		$\tilde{\mathcal F}_\lambda$ be the normal, $<\!\kappa$–complete filter from
		Theorem~\ref{thm:kappa-filter-generation}. We also use Lemma~\ref{lem:kappa-tuples-fix} (stabilizers of $<\kappa$-tuples) and
		Lemma~\ref{lem:kappa-support-hs} (support $<\kappa$ implies symmetry). All hereditary conclusions follow
		by induction on name rank.
		
		\emph{Check names.}
		For $a\in V$, every $g\in\mathcal G_\lambda$ fixes $\check a$, hence
		$\fix(\check a)=\mathcal G_\lambda\in\tilde{\mathcal F}_\lambda$; thus
		$\check a$ is HS.
		
		\emph{Pairing (and ordered pairs).}
		If $\dot x,\dot y$ are HS, then by Lemma~\ref{lem:kappa-tuples-fix}
		\[
		\fix\!\bigl(\langle \dot x,\dot y\rangle\bigr)\;=\;\fix(\dot x)\cap\fix(\dot y)
		\ \in\ \tilde{\mathcal F}_\lambda,
		\]
		since $\tilde{\mathcal F}_\lambda$ is a filter (closed under finite
		intersections). Hence $\langle\dot x,\dot y\rangle$ and $\{\dot x,\dot y\}$
		are HS.
		
		\emph{$<\kappa$-tuples.}
		More generally, if $\delta<\kappa$ and $\langle \dot x_i:i<\delta\rangle$ is a $\delta$-sequence of HS names,
		then by Lemma~\ref{lem:kappa-tuples-fix},
		\[
		\fix\!\bigl(\langle \dot x_i:i<\delta\rangle\bigr)=\bigcap_{i<\delta}\fix(\dot x_i)\in\tilde{\mathcal F}_\lambda
		\]
		by $\kappa$-completeness of $\tilde{\mathcal F}_\lambda$.
		Thus $<\kappa$-tuples of HS names are HS.
		
		\emph{Union of an HS name.}
		Let $\dot A$ be HS and define the usual union-name
		\[
		\dot U\ :=\ \Bigl\{\,\langle \tau,r\rangle : \exists\langle\sigma,p\rangle\in\dot A\ \exists\langle\tau,q\rangle\in\sigma\
		(r\le p,q)\Bigr\}.
		\]
		If $g\in\fix(\dot A)$ then $g(\dot U)=\dot U$ (because $g$ transports witnesses in the defining formula),
		so $\fix(\dot U)\supseteq\fix(\dot A)\in\tilde{\mathcal F}_\lambda$. Hence $\dot U$ is HS.
		
		\emph{Bounded Separation (over a fixed ambient set).}
		Fix an HS name $\dot A$, HS parameters $\bar{\dot p}$, and a bounded
		formula $\varphi$. Define the standard comprehension name
		\[
		\dot B\ :=\ \{\ \langle \sigma,p\rangle\in \dot A\ :\ p\Vdash \varphi(\sigma,\bar{\dot p})\ \}.
		\]
		If $g\in \fix(\dot A)\cap\fix(\bar{\dot p})$, then by automorphism invariance
		of the forcing relation for bounded formulas,
		$p\Vdash\varphi(\sigma,\bar{\dot p})\iff g(p)\Vdash\varphi(g\sigma,\bar{\dot p})$,
		so $g\in\fix(\dot B)$. Therefore
		\[
		\fix(\dot B)\ \supseteq\ \fix(\dot A)\cap\fix(\bar{\dot p})\ \in\ \tilde{\mathcal F}_\lambda,
		\]
		and $\dot B$ is HS.
		
		\emph{Images under ground-model functions (typical “etc.” case).}
		Let $f\in V$ be a function and $\dot A$ HS. Define the image name
		$\dot f[\dot A]:=\{\ \langle \tau,r\rangle : \exists\langle\sigma,p\rangle\in\dot A\
		(r\le p\ \&\ r\Vdash \tau=\check f(\sigma))\ \}$. If $g\in\fix(\dot A)\cap\fix(\check f)$,
		then $g\in\fix(\dot f[\dot A])$ by definition; hence
		$\fix(\dot f[\dot A])\supseteq\fix(\dot A)\cap\fix(\check f)\in\tilde{\mathcal F}_\lambda$,
		so $\dot f[\dot A]$ is HS.
		
		Finally, Lemma~\ref{lem:kappa-support-hs} guarantees that any constructor that
		produces a name with support $<\kappa$ yields an HS name; the cases above show
		the standard operations do so and have stabilizers in $\tilde{\mathcal F}_\lambda$.
		Therefore the class of HS names (for $<\!\kappa$–support) is closed under the stated
		operations. \qedhere
	\end{proof}
	
	\subsection{ZF and choice fragments at \texorpdfstring{$<\kappa$}{<k}-support}
	\label{subsec:kappa-preservation}
	
	\begin{theorem}[ZF is preserved]
		\label{thm:kappa-zf-preserved}
		If $V\models\mathrm{ZF}$, then $M:=V(\mathbb{P}_\lambda)_{\mathrm{HS}}\models\mathrm{ZF}$ for any
		$<\kappa$-support symmetric iteration. 
	\end{theorem}
	
	\begin{proof}
		Write $\fix(\dot x)=\fix_{\mathcal G_\lambda}(\dot x)$ and let $\tilde{\mathcal F}_\lambda$ be the
		(normal) limit filter on $\mathcal G_\lambda$ from Theorem~\ref{thm:kappa-filter-generation}.
		We use the $<\kappa$--HS closure theorem (Theorem~\ref{thm:kappa-hs-closure}) and the standard
		automorphism invariance of the forcing relation.
		
		\smallskip\noindent\textbf{Transitivity and sandwich.}
		If $a\in V$, then $\check a$ is fixed by all $g\in\mathcal G_\lambda$, hence HS and $a=\check a^{\,G}\in M$;
		thus $V\subseteq M$. Trivially $M\subseteq V[G]$. For transitivity, let $y\in M$ with $y=\dot y^{\,G}$
		for some HS name $\dot y$. If $x\in y$, then $x=\sigma^{G}$ for some $(\sigma,p)\in\dot y$ with $p\in G$.
		Since $\dot y$ is \emph{hereditarily} symmetric, every constituent $\sigma$ of $\dot y$ is HS, hence $x\in M$.
		Therefore $M$ is transitive.
		
		\smallskip\noindent\textbf{Extensionality and Foundation.}
		These are standard for symmetric submodels: $M$ is a transitive subclass of $V[G]$ with inherited $\in$,
		and the forcing interpretation of names is extensional and well-founded; see, e.g., \cite[Ch.~15]{Jech03}.
		
		\smallskip\noindent\textbf{Pairing, Union, Infinity.}
		Pairing and Union hold in $M$ because HS is closed under pairing and the union-name constructor
		(Theorem~\ref{thm:kappa-hs-closure}). Infinity holds since $\check\omega$ is fixed by all automorphisms, hence HS.
		
		\smallskip\noindent\textbf{Separation.}
		Let $A\in M$ be given by an HS name $\dot A$, fix HS parameters $\vec{\dot p}$, and let $\varphi$ be any formula.
		Define the name
		\[
		\dot B\ :=\ \{\ \langle\sigma,p\rangle\in\dot A\ :\ p\Vdash \varphi(\sigma,\vec{\dot p})\ \}.
		\]
		If $g\in\fix(\dot A)\cap\fix(\vec{\dot p})$, then by automorphism invariance of the forcing relation,
		$p\Vdash \varphi(\sigma,\vec{\dot p})\iff g(p)\Vdash \varphi(g\sigma,\vec{\dot p})$, whence $g\in\fix(\dot B)$.
		Thus
		\[
		\fix(\dot B)\ \supseteq\ \fix(\dot A)\cap\fix(\vec{\dot p})\ \in\ \tilde{\mathcal F}_\lambda,
		\]
		so $\dot B$ is HS and $B:=\dot B^{\,G}=\{x\in A:\varphi^M(x,\vec p)\}$. Hence $M$ satisfies Separation.
		
		\smallskip\noindent\textbf{Power Set.}
		Let $a\in M$ with $a=\dot a^{\,G}$ for HS $\dot a$. Fix $\alpha:=\rank(\dot a)+1$ and let
		$\mathrm{Nm}_{\le\alpha}$ be the (set-sized) collection in $V$ of $\mathbb P_\lambda$--names of rank $\le\alpha$.
		Set
		\[
		\HS(\dot a)\ :=\ \{\ \dot x\in \mathrm{Nm}_{\le\alpha}\ :\ \Vdash \dot x\subseteq \dot a\ \text{ and }\ \dot x\in\mathrm{HS}\ \}.
		\]
		Then $\HS(\dot a)$ is a \emph{set} in $V$ (a definable subset of $\mathrm{Nm}_{\le\alpha}$). Define the name
		\[
		\dot{\mathcal P}\ :=\ \{\ \langle \dot x, 1_{\mathbb P_\lambda}\rangle\ :\ \dot x\in \HS(\dot a)\ \}.
		\]
		If $g\in\fix(\dot a)$, then $g$ preserves the relation “$\Vdash \dot x\subseteq \dot a$” and preserves the HS property,
		so $g[\HS(\dot a)]=\HS(\dot a)$ and hence $g(\dot{\mathcal P})=\dot{\mathcal P}$. Therefore
		$\fix(\dot{\mathcal P})\supseteq \fix(\dot a)\in\tilde{\mathcal F}_\lambda$, so $\dot{\mathcal P}$ is HS and
		\[
		\dot{\mathcal P}^{\,G}\ =\ \{\ \dot x^{\,G}:\ \dot x\in\HS(\dot a)\ \}\ =\ \mathcal P^M(a).
		\]
		Thus $M$ satisfies Power Set.
		
		\smallskip\noindent\textbf{Replacement (via functional formulas).}
		Let $A\in M$ be given by an HS name $\dot A$, and fix HS parameters $\vec{\dot p}$.
		Suppose $p_0\in\mathbb P_\lambda$ forces that $\varphi(x,y,\vec p)$ defines a function on $A$, i.e.
		\[
		p_0\Vdash\ \forall x\in \dot A\ \exists!\, y\ \varphi(x,y,\vec{\dot p}).
		\]
		Define the relational image name
		\[
		\dot B\ :=\ \Bigl\{\ \langle\tau,p\rangle\ :\ \exists\langle\sigma,q\rangle\in\dot A\ (p\le q\ \wedge\ p\Vdash \varphi(\sigma,\tau,\vec{\dot p}))\ \Bigr\}.
		\]
		If $g\in\fix(\dot A)\cap\fix(\vec{\dot p})$, then by invariance of the forcing relation,
		$p\Vdash \varphi(\sigma,\tau,\vec{\dot p})\Rightarrow g(p)\Vdash \varphi(g\sigma,g\tau,\vec{\dot p})$ and
		$(\sigma,q)\in\dot A\Rightarrow (g\sigma,g(q))\in\dot A$, hence $g\in\fix(\dot B)$. Therefore
		\[
		\fix(\dot B)\ \supseteq\ \fix(\dot A)\cap\fix(\vec{\dot p})\ \in\ \tilde{\mathcal F}_\lambda,
		\]
		so $\dot B$ is HS and $B:=\dot B^{\,G}\in M$. By functionality, in $V[G]$ we have
		$B=\{\ y : \exists x\in A\ \varphi(x,y,\vec p)\ \}$, so $M$ satisfies Replacement.
		
		\smallskip
		We have verified that $M$ is a transitive inner model of $\ZF$ with $V\subseteq M\subseteq V[G]$. \qedhere
	\end{proof}
	
	\begin{lemma}[Closure under $<\kappa$-sequences in the forcing extension $V^G$]
		\label{lem:closure-<kappa-seqs}
		Let $\mu<\kappa$. If $\langle x_\xi:\xi<\mu\rangle\in V[G]$ and $x_\xi\in M$ for all $\xi<\mu$, then
		$\langle x_\xi:\xi<\mu\rangle\in M$.
	\end{lemma}
	
	\begin{proof}
		For each $\xi<\mu$ choose an HS name $\dot x_\xi$ with $x_\xi=\dot x_\xi^{\,G}$.
		Let $\dot s:=\langle \dot x_\xi:\xi<\mu\rangle$ be the canonical $\mu$--sequence name.
		By Lemma~\ref{lem:kappa-tuples-fix},
		\[
		\fix(\dot s)\ =\ \bigcap_{\xi<\mu}\fix(\dot x_\xi).
		\]
		Since each $\fix(\dot x_\xi)\in\tilde{\mathcal F}_\lambda$ and $\mu<\kappa$, $\kappa$--completeness of
		$\tilde{\mathcal F}_\lambda$ (Theorem~\ref{thm:kappa-filter-generation}) gives $\fix(\dot s)\in\tilde{\mathcal F}_\lambda$,
		so $\dot s\in\HS$. Therefore $\dot s^{\,G}=\langle x_\xi:\xi<\mu\rangle\in M$.
	\end{proof}
	
	\begin{theorem}[$DC_{<\kappa}$ transfers to $M$]
		\label{thm:DCless-preservation}
		Let $\lambda$ be a limit stage of a $<\kappa$-support symmetric iteration as in Definition~\ref{def:kappa-iteration},
		and assume the limit filter $\tilde{\mathcal F}_\lambda$ is normal and $\kappa$-complete
		(Theorem~\ref{thm:kappa-filter-generation}). Fix $\mu<\kappa$.
		
		If $V[G]\models DC_\mu$, then $M=V(\mathbb P_\lambda)_{\mathrm{HS}}\models DC_\mu$.
		In particular, if $V\models\ZFC$ then $V[G]\models DC_\mu$ for all $\mu$, hence $M\models DC_{<\kappa}$.
	\end{theorem}
	
	\begin{proof}
		Work in $V[G]$. Let $A,R\in M$ and assume $R$ is serial on $A$, i.e.
		$\forall x\in A\ \exists y\in A\ \langle x,y\rangle\in R$.
		By the hypothesis $V[G]\models DC_\mu$, there exists a sequence $\langle a_\xi:\xi<\mu\rangle$
		with $a_\xi\in A$ and $\langle a_\xi,a_{\xi+1}\rangle\in R$ for all $\xi+1<\mu$.
		Since $a_\xi\in M$ for each $\xi$, Lemma~\ref{lem:closure-<kappa-seqs} implies that the whole
		sequence $\langle a_\xi:\xi<\mu\rangle$ belongs to $M$. Thus $M\models DC_\mu$.
	\end{proof}
	
	\begin{definition}[Localization Hypothesis for $DC_\kappa$]
		\label{def:localization-hypothesis}
		Let $\kappa$ be regular and let $\lambda$ be a limit stage of a $<\kappa$--support symmetric iteration.
		We say that the iteration satisfies \emph{Localization at $\lambda$} if the following holds:
		
		Whenever $\dot A,\dot R\in\HS$ and $p\in\mathbb P_\lambda$ forces that $\dot R$ is serial on $\dot A$, i.e.
		\[
		p\Vdash\ \forall x\in \dot A\ \exists y\in \dot A\ \bigl(\langle x,y\rangle\in \dot R\bigr),
		\]
		there exist a condition $q\le p$ and a $\mathbb P_\lambda$--name $\dot s$ with $|\supp(\dot s)|<\kappa$ such that
		\[
		q\Vdash\ \bigl(\dot s:\check\kappa\to\dot A\ \wedge\ \forall\xi<\check\kappa\
		\langle \dot s(\xi),\dot s(\xi{+}1)\rangle\in\dot R\bigr).
		\]
	\end{definition}
	
	\begin{example}[Heuristic boundary for localization]\label{ex:nonlocalizing}
		Let $\langle A_\xi:\xi<\kappa\rangle$ be a family where each $A_\xi$ is added at (or supported by) coordinate $\xi$, and let $R$ be the relation on partial choice functions $f$ with
		\[
		R(f)\;:\Longleftrightarrow\;
		\Bigl(\forall \xi\in\dom(f)\; f(\xi)\in A_\xi\Bigr)\ \wedge\
		\left[
		\begin{aligned}
			&\text{for cofinally many } \alpha<\kappa,\\
			&\ran\bigl(f\!\upharpoonright\! (B_\alpha\cap\dom(f))\bigr)\ \text{is injective}
		\end{aligned}
		\right].
		\]
		Here $b:\kappa\times\kappa\to\kappa$ is a fixed bijection and we write
		$B_\alpha:=\{\,b(\alpha,\beta):\beta<\kappa\,\}$, so that
		$\langle B_\alpha:\alpha<\kappa\rangle$ partitions $\kappa$ into $\kappa$ many blocks, each of size $\kappa$.
		The injectivity demand ties together \emph{all} coordinates appearing in each large block.
		 
		Heuristically, $R$ can \emph{fail} localization: for any finite support $u\subset\kappa$ and condition $p$, one can usually find $q_0,q_1\le p$ agreeing on $u$ so that $q_0$ forces “$R(\dot f)$” while $q_1$ forces “$\neg R(\dot f)$,” by modifying choices off of $u$ within infinitely many future blocks. 
		This illustrates that relations imposing blockwise-cofinal global constraints are outside the scope of the localization hypothesis in its basic form.
	\end{example}
	
	\begin{theorem}[$DC_\kappa$ under Localization]
		\label{thm:DCkappa-localization}
		Assume $V\models \mathrm{ZF}+\mathrm{DC}$, $\kappa$ is regular, the limit filter $\tilde{\mathcal F}_\lambda$ is normal and $\kappa$-complete, and the iteration satisfies Localization (Definition~\ref{def:localization-hypothesis}). Then 
		$M=V(\mathbb P_\lambda)_{\mathrm{HS}}\models DC_\kappa$.
	\end{theorem}
	
	\begin{proof}
		Let $\dot A,\dot R\in\HS$ be names with
		\[
		\text{Fix }p\in\mathbb P_\lambda\text{ such that }\ 
		p\Vdash\ \forall x\in \dot A\ \exists y\in \dot A\ \bigl(\langle x,y\rangle\in \dot R\bigr).
		\]
		We verify $DC_\kappa$ in $M:=V(\mathbb P_\lambda)_{\HS}$.
		
		\smallskip\noindent\emph{Localization input.}
		By Localization (Definition~\ref{def:localization-hypothesis}), there exist $q\le p$ and a name
		$\dot s\in\HS$ with $|\supp(\dot s)|<\kappa$ such that
		\begin{equation}\label{eq:loc-seq}
			q\Vdash\ \bigl(\dot s:\check\kappa\to\dot A\ \wedge\ \forall\xi<\check\kappa\
			\langle \dot s(\xi),\dot s(\xi{+}1)\rangle\in\dot R\bigr).
		\end{equation}
		
		\smallskip\noindent\emph{The sequence name is HS.}
		Since $\supp(\dot s)\subseteq S$ with $|S|<\kappa$, Lemma~\ref{lem:kappa-support-hs} (\emph{support $<\kappa$ implies symmetry})
		gives
		\[
		\fix_{\mathcal G_\lambda}(\dot s)\ \in\ \tilde{\mathcal F}_\lambda,
		\]
		so $\dot s\in\HS$. In particular, $\dot s^G\in M$ for every $V$--generic $G\subseteq\mathbb P_\lambda$.
		
		\smallskip\noindent\emph{$DC_\kappa$ in $M$.}
		Let $A:=\dot A^G$ and $R:=\dot R^G$ and choose $G\subseteq\mathbb P_\lambda$ $V$--generic with $q\in G$; then by \eqref{eq:loc-seq}, in $V[G]$ we have that
		$s:=\dot s^G$ is a function $s:\kappa\to A$ with $\langle s(\xi),s(\xi{+}1)\rangle\in R$ for all $\xi<\kappa$.
		Because $\dot s\in\HS$, $s\in M$. Thus $M$ contains a choice sequence of length $\kappa$ witnessing $DC_\kappa$
		for the given serial relation on $A$.
		
		\smallskip
		Since $\dot A,\dot R$ were arbitrary HS names for a serial relation, it follows that $M\models DC_\kappa$. \qedhere
	\end{proof}
	
	\begin{remark}[Choice landscape at $<\kappa$-support]
		\label{rem:AC-hierarchy}
		Assume $\kappa$ is regular. As in the countable case, global $\mathrm{AC}$ typically fails by permutation
		arguments using fresh coordinates. For fragments: the $\kappa$-completeness of $\tilde{\mathcal{F}}_\lambda$
		yields $DC_{<\kappa}$ (Theorem~\ref{thm:DCless-preservation}). Under Localization
		(Definition~\ref{def:localization-hypothesis}) one further obtains $DC_\kappa$
		(Theorem~\ref{thm:DCkappa-localization}). Stronger closure hypotheses (e.g., symmetry systems whose HS names are
		closed under suitable $<\kappa$-indexed constructions, or $<\kappa$-closed forcing/automorphism features in
		particular examples) may preserve additional choice fragments such as $AC_{<\kappa}$, but the precise outcome
		is forcing- and symmetry-dependent.
	\end{remark}
	
	\begin{remark}[Localization: sufficiency vs.\ unknown necessity]
		\label{rem:localization-sufficiency}
		Theorem~\ref{thm:DCkappa-localization} employs the Localization Hypothesis (Definition~\ref{def:localization-hypothesis})
		as a \emph{sufficient} criterion for preserving $\mathrm{DC}_\kappa$ through our iterations.
		At present we do not know whether Localization is \emph{necessary} in general: there may exist
		HS relations $R$ and iterations where $\mathrm{DC}_\kappa$ holds although Localization fails.
		For coordinate-structured forcings (Cohen, random, Hechler) Localization appears plausible in
		many natural situations; in Section~6.1 we verify it for a broad \emph{reuse-permissive} class
		(Proposition~\ref{prop:hechler-localization-reuse}). Accordingly, we treat Localization as an
		\emph{explicit assumption} and verify it case by case, avoiding stronger claims of equivalence.
	\end{remark}
	
	\subsection*{Potential weakenings of Localization}
	\label{subsec:weaker-localization}
	
	We record several weakenings of Localization that may still suffice for preserving $\DC_\kappa$.
	
	\begin{definition}[Coding–only Localization $\mathrm{CLoc}_\kappa$]
		\label{def:cloc}
		Let $\kappa$ be regular. We say that $\mathrm{CLoc}_\kappa$ holds for an HS instance
		$\langle \dot A,\dot R\rangle$ at $\lambda$ if whenever $p\in\mathbb P_\lambda$ forces the $\DC_\kappa$
		antecedent
		\[
		p\Vdash\ \forall x\in \dot A\ \exists y\in \dot A\ \bigl(\langle x,y\rangle\in \dot R\bigr),
		\]
		then there exist $q\le p$ and a \emph{single} $\mathbb P_\lambda$--name $\dot s$ with $|\supp(\dot s)|<\kappa$ such that
		\[
		q\Vdash\ \Bigl(\dot s:\check\kappa\to \dot A\ \wedge\ \forall\xi<\check\kappa\ 
		\langle \dot s(\xi),\dot s(\xi{+}1)\rangle\in \dot R\Bigr).
		\]
		(We do not separately assume a uniform support bound for the coordinates; it follows from $|\supp(\dot s)|<\kappa$.)
	\end{definition}
	
	\begin{proposition}[$\mathrm{CLoc}_\kappa$ suffices for $\DC_\kappa$]
		\label{prop:cloc-implies-dc}
		If $\mathrm{CLoc}_\kappa$ holds for every HS relation $R$ with the $\DC_\kappa$ antecedent, then the symmetric model at $\lambda$ satisfies $\DC_\kappa$.
	\end{proposition}
	
	\begin{proof}
		Let $G\subseteq\mathbb P_\lambda$ be $V$–generic and $M:=V(\mathbb P_\lambda)_{\HS}\subseteq V[G]$.
		To verify $DC_\kappa$ in $M$, fix any $A,R\in M$ with $R\subseteq A\times A$ such that
		$M\models\forall x\in A\ \exists y\in A\ \langle x,y\rangle\in R$ (i.e., $R$ is serial on $A$).
		
		Choose HS names $\dot A,\dot R$ with $\dot A^G=A$ and $\dot R^G=R$, and (strengthening a condition if needed) assume
		\[
		\Vdash\ \forall x\in \dot A\ \exists y\in \dot A\ \bigl(\langle x,y\rangle\in \dot R\bigr).
		\]
		By the hypothesis $\mathrm{CLoc}_\kappa$ (Definition~\ref{def:cloc}) applied to this HS instance $\langle\dot A,\dot R\rangle$
		and the condition $p$ above, there exist $q\le p$ and a $\mathbb P_\lambda$--name $\dot s$ with $|\supp(\dot s)|<\kappa$ such that
		\begin{equation}\label{eq:cloc-seq}
			q\Vdash\ \Bigl(\dot s:\check\kappa\to \dot A\ \wedge\ \forall \xi<\check\kappa\ \langle \dot s(\xi),\dot s(\xi{+}1)\rangle\in \dot R\Bigr).
		\end{equation}
		Fix $G\subseteq\mathbb P_\lambda$ $V$--generic with $q\in G$.
		
		Since $|\supp(\dot s)|<\kappa$, Lemma~\ref{lem:kappa-support-hs} (\emph{support $<\kappa$ implies symmetry}) yields
		$\fix_{\mathcal G_\lambda}(\dot s)\in\tilde{\mathcal F}_\lambda$, hence $\dot s\in\HS$ and therefore $s:=\dot s^G\in M$.
		By \eqref{eq:cloc-seq}, in $V[G]$ (hence in $M$) we have that $s:\kappa\to A$ and
		$\langle s(\xi),s(\xi{+}1)\rangle\in R$ for all $\xi<\kappa$. Thus $s$ witnesses $DC_\kappa$ in $M$.
		
		As $A$ and $R$ were arbitrary (with $R$ serial on $A$), it follows that $M\models DC_\kappa$. \qedhere
	\end{proof}
	
	\begin{definition}[Club–Localization $\mathrm{CLubLoc}_\kappa$]
		\label{def:clubloc}
		Let $\kappa$ be regular. We say that $\mathrm{CLubLoc}_\kappa$ holds for an HS instance
		$\langle \dot A,\dot R\rangle$ if whenever $p\in\mathbb P_\lambda$ forces the $\DC_\kappa$ antecedent
		\[
		p\Vdash\ \forall x\in \dot A\ \exists y\in \dot A\ \bigl(\langle x,y\rangle\in \dot R\bigr),
		\]
		there exist $q\le p$, a set $S\subseteq\lambda$ with $|S|<\kappa$, a club $C\subseteq\kappa$, and an HS name $\dot z$
		such that $|\supp(\dot z)|<\kappa$, $\supp(\dot z)\subseteq S$, and
		\[
		q\Vdash\ \Bigl(\dot z:C\to \dot A\ \wedge\ \forall \xi\in C\ \bigl\langle \dot z(\xi),\dot z(\min(C\setminus(\xi{+}1)))\bigr\rangle\in \dot R\Bigr).
		\]
	\end{definition}
	
	\begin{proposition}[$\mathrm{CLubLoc}_\kappa$ suffices for reuse–permissive $R$]
		\label{prop:clubloc-reuse}
		If $R$ is reuse–permissive (Def.~\ref{def:hechler-reuse-permissive}) and $\mathrm{CLubLoc}_\kappa$ holds for $R$, then there is a full length–$\kappa$ witness in the symmetric model; hence $\DC_\kappa$ holds for $R$.
	\end{proposition}
	
	\begin{proof}
		Let $G\subseteq\mathbb P_\lambda$ be $V$–generic and $M:=V(\mathbb P_\lambda)_{\HS}$.
		Fix $A,R\in M$ with $R\subseteq A\times A$ such that $R$ is serial on $A$ and \emph{reuse–permissive}
		(Def.~\ref{def:hechler-reuse-permissive}). By hypothesis $\CLubLoc_\kappa$ holds for this $R$; thus there exist:
		\begin{itemize}
			\item a club $C\subseteq \kappa$ (in $V$), with increasing enumeration $e:\kappa\to C$, $\eta\mapsto \gamma_\eta$;
			\item HS names $\dot A,\dot R$ with $\dot A^G=A$, $\dot R^G=R$; and
			\item a condition $q\in\mathbb P_\lambda$ and an HS name $\dot z$ such that $|\supp(\dot z)|<\kappa$ and
			\[
			q\Vdash\ \Bigl(\dot z:C\to \dot A\ \wedge\ \supp(\dot z)\subseteq S\Bigr)
			\ \wedge\ 
			\forall \eta<\check\kappa\ \bigl\langle \dot z(\gamma_\eta),\,\dot z(\gamma_{\eta+1})\bigr\rangle\in \dot R,
			\]
			for some $S\subseteq\lambda$ with $|S|<\kappa$.
		\end{itemize}
		
		Fix $G\subseteq\mathbb P_\lambda$ $V$--generic with $q\in G$.
		
		By Lemma~\ref{lem:kappa-support-hs} (\emph{support $<\kappa$ implies symmetry}), $\dot z\in\HS$.
		
		\smallskip
		Work in $V[G]$ and set $z:=\dot z^G$, $A:=\dot A^G$, and $R:=\dot R^G$.
		Let $\gamma_0:=\min(C)$ and define the ground-model function $\rho_C:\kappa\to C$ by
		\[
		\rho_C(\xi)\ :=\ 
		\begin{cases}
			\gamma_0, & \xi<\gamma_0,\\[2pt]
			\max(C\cap(\xi+1)), & \xi\ge \gamma_0.
		\end{cases}
		\]
		(For $\xi\ge\gamma_0$, the set $C\cap(\xi+1)$ is nonempty and closed in $\xi+1$, hence has a maximum.)
		
		Define in $V[G]$ the function $s:\kappa\to A$ by
		\[
		s(\xi)\ :=\ z(\rho_C(\xi)),\qquad \xi<\kappa.
		\]
		Since $z\in M$ and $\rho_C\in V\subseteq M$, the function $s$ is definable from parameters in $M$ and hence
		$s\in M$ (by Replacement in $M$).
		
		\smallskip
		We verify that $s$ witnesses $DC_\kappa$ for $R$ in $M$.
		Work in $V[G]$ and set $z:=\dot z^G$, $s:=\dot s^G$,
		$A:=\dot A^G$, $R:=\dot R^G$. Then $s=F(z)$ and $s:\kappa\to A$. Fix $\xi<\kappa$ and write
		$\eta$ for the unique index with $\gamma_\eta=\max(C\cap(\xi+1))$.
		
		\emph{Interior steps.} If $\max(C\cap(\xi+2))=\gamma_\eta$ as well (i.e., $\xi+1<\gamma_{\eta+1}$), then
		$s(\xi)=s(\xi+1)=z(\gamma_\eta)$. Since $R$ is reuse–permissive (Def.~\ref{def:hechler-reuse-permissive}),
		we have $(s(\xi),s(\xi+1))=(z(\gamma_\eta),z(\gamma_\eta))\in R$ (the “stay”/reuse clause).
		
		\emph{Club–successor steps.} If $\max(C\cap(\xi+2))=\gamma_{\eta+1}$ (equivalently, $\xi+1=\gamma_{\eta+1}$),
		then $s(\xi)=z(\gamma_\eta)$ and $s(\xi+1)=z(\gamma_{\eta+1})$, and by the $\CLubLoc_\kappa$ clause,
		$\bigl(z(\gamma_\eta),z(\gamma_{\eta+1})\bigr)\in R$.
		
		Thus, for all $\xi<\kappa$ we have $\langle s(\xi),s(\xi+1)\rangle\in R$. Since $s\in M$, $M$ contains a full length–$\kappa$ witness for $R$. Hence $DC_\kappa$ holds in $M$ for this reuse–permissive $R$.
		
		\smallskip
		This proves the proposition. \qedhere
	\end{proof}
	
	\begin{definition}[Tail–Localization $\mathrm{TailLoc}_\kappa$]
		\label{def:tailloc}
		Let $\kappa$ be regular. We say that $\mathrm{TailLoc}_\kappa$ holds for an HS instance $\langle\dot A,\dot R\rangle$
		if whenever $p\in\mathbb P_\lambda$ forces the $\DC_\kappa$ antecedent on $\dot A$ via $\dot R$, there exist
		$q\le p$, an ordinal $\beta<\lambda$, and a set $S\subseteq[\beta,\lambda)$ with $|S|<\kappa$ such that $q$ forces
		the existence of a witnessing sequence $\langle \dot x_\xi:\xi<\kappa\rangle$ with
		\[
		\exists \xi_0<\check\kappa\ \forall \xi\ge \xi_0\ \bigl(\supp(\dot x_\xi)\subseteq S\bigr).
		\]
	\end{definition}
	
	\begin{remark}[Implications and scope]
		\label{rem:loc-implications}
		Localization (Def.~\ref{def:localization-hypothesis}) implies both $\mathrm{CLoc}_\kappa$ and $\mathrm{CLubLoc}_\kappa$. For coordinate–structured forcings (Cohen, Hechler, random), our verification for reuse–permissive $R$ (Prop.~\ref{prop:hechler-localization-reuse}) typically yields $\mathrm{CLoc}_\kappa$ and $\mathrm{CLubLoc}_\kappa$ as well. 
		\emph{Open}: Are there iterations where $\DC_\kappa$ holds but $\mathrm{CLoc}_\kappa$ fails for some $R$? Does $\mathrm{TailLoc}_\kappa$ strictly weaken Localization in general?
	\end{remark}
	
	\begin{conjecture}[Localization may not be necessary]\label{conj:localization-not-necessary}
		There is a symmetric iteration preserving $\DC_\kappa$ in which the Localization Hypothesis
		(Definition~\ref{def:localization-hypothesis}) fails for some HS relation $R$.
	\end{conjecture}
	
	\begin{remark}[Speculative counterexample shape]\label{rem:localization-counterexample-shape}
		A plausible direction is to use a non–coordinate-structured forcing at cofinally many stages
		(e.g., proper but not product-like), together with a relation $R$ whose DC\(_\kappa\) witnesses
		intrinsically require $\kappa$ many fresh coordinates. Localization would then fail because no single
		name (nor a fixed $S$ of size $<\kappa$) codes the whole sequence, yet $\DC_\kappa$ might still be
		preserved via fusion/compactness at the level of names rather than supports. We leave constructing
		such an example to future work.
	\end{remark}
	
	\paragraph{Concrete verification of Localization for structured forcings.}
	While Localization (Definition~\ref{def:localization-hypothesis}) is a hypothesis to be verified case by case,
	for forcings with coordinate structure (Cohen, random, Hechler) a natural class of HS relations admits a uniform
	verification via \emph{block reuse}. We now isolate that class and prove a sufficient condition for DC$_\kappa$.
	
	\begin{definition}[Reuse-permissive HS relation]\label{def:hechler-reuse-permissive}
		Let $\kappa$ be regular uncountable. An HS relation $R$ on $(\omega^\omega)^{<\kappa}$ is
		\emph{reuse-permissive} if there exist a cardinal $\mu<\kappa$, a \emph{cofinal} partition $b:\kappa\to\mu$ \textup{(}a surjection with each fiber $b^{-1}(\{i\})$ unbounded in $\kappa$\textup{)}, and a finite scheme of
		\emph{monotone local constraints} $\mathcal{C}_R$ such that:
		\begin{enumerate}
			\item \textbf{Block-constant closure.} If $\langle g_i:i<\mu\rangle\subseteq\omega^\omega$
			satisfies $\mathcal{C}_R$ (constraints among $\{g_j:j\le i\}$ for each $i$), then the expanded
			sequence $f_\xi:=g_{b(\xi)}$ witnesses the DC$_\kappa$ instance of $R$.
			\item \textbf{Refinement monotonicity.} If $\langle g_i\rangle$ satisfies $\mathcal{C}_R$
			and $\langle g'_i\rangle$ Hechler-refines it pointwise (e.g., $g'_i\ge^* g_i$ or $g'_i=g_i$),
			then $\langle g'_i\rangle$ also satisfies $\mathcal{C}_R$.
		\end{enumerate}
	\end{definition}
	
	\begin{lemma}[Finite constraint realization for Hechler coordinates]
		\label{lem:hechler-finite-constraints}
		Let $F\subseteq \lambda$ be finite in a $<\kappa$-support iteration of Hechler forcing. Fix a finite set $\Gamma$ of constraints of the following forms:
		\begin{enumerate}
			\item[(i)] (\emph{pointwise}) for some finite $J\subseteq F\times\omega$ and integers $m_{\alpha,n}$,
			literals of the form $\ \dot d_\alpha(n)=\check m_{\alpha,n}\ $ or $\ \dot d_\alpha(n)\neq\check m_{\alpha,n}$;
			\item[(ii)] (\emph{eventual lower bounds}) for some finite $F_0\subseteq F$ and ground-model functions $f_\alpha\in\omega^\omega\ (\alpha\in F_0)$,
			the literals $\ \dot d_\alpha\ge^* \check f_\alpha$.
		\end{enumerate}
		Then for every condition $p$ there is $q\le p$ with $\supp(q)\subseteq \supp(p)\cup F$ that \emph{forces all constraints in $\Gamma$ simultaneously}.
	\end{lemma}
	
	\begin{proof}
		Write each coordinate condition as $p(\alpha)=(s_\alpha^p,f_\alpha^p)\in\mathbb D$.
		Let 
		\[
		N\ :=\ 1+\max\{\,n : (\alpha,n)\in J\,\}
		\]
		(with $N=0$ if $J=\varnothing$). We build $q\le p$ coordinatewise for $\alpha\in F$, leaving other coordinates unchanged.
		
		\smallskip\noindent\emph{Stage 1: extend stems to decide all pointwise literals.}
		For each $\alpha\in F$, choose a finite sequence $t_\alpha\in \omega^{N}$ extending $s_\alpha^p$ such that:
		\begin{itemize}
			\item if $\dot d_\alpha(n)=\check m_{\alpha,n}\in\Gamma$ then set $t_\alpha(n)=m_{\alpha,n}$;
			\item if $\dot d_\alpha(n)\neq\check m_{\alpha,n}\in\Gamma$ then set $t_\alpha(n)\neq m_{\alpha,n}$ (always possible, pick any other integer);
			\item for $n\in N\setminus\dom(s_\alpha^p)$ also ensure $t_\alpha(n)\ge f_\alpha^p(n)$ to respect the Hechler order.
		\end{itemize}
		Set $s_\alpha^q:=t_\alpha$. This decides every literal in (i) because Hechler stems decide finite coordinates.
		
		\smallskip\noindent\emph{Stage 2: raise eventual lower bounds.}
		For each $\alpha\in F_0$, define a new function $g_\alpha\in\omega^\omega$ by
		\[
		g_\alpha(n)\ :=\ \max\bigl(f_\alpha^p(n), f_\alpha(n)\bigr)\quad\text{for all }n\in\omega.
		\]
		For $\alpha\in F\setminus F_0$ let $g_\alpha:=f_\alpha^p$. Now set $f_\alpha^q:=g_\alpha$ for all $\alpha\in F$.
		
		\smallskip
		Define $q$ by $q(\alpha)=(s_\alpha^q,f_\alpha^q)$ for $\alpha\in F$ and $q(\alpha)=p(\alpha)$ otherwise. Then $q\le p$, $\supp(q)\subseteq \supp(p)\cup F$, and:
		\begin{itemize}
			\item By construction of the stems $s_\alpha^q$, $q$ forces every literal in (i).
			\item For each $\alpha\in F_0$, $q$ forces $\dot d_\alpha\ge^* f_\alpha^q\ge^* f_\alpha$, hence the literals in (ii).
		\end{itemize}
		Thus $q$ forces all constraints in $\Gamma$.
	\end{proof}
	
	\begin{remark}[What is \emph{not} forceable]
		Under the standard Hechler presentation, a single condition cannot force inter-generic tail relations such as
		$\dot d_\alpha =^* \dot d_\beta$ or $\dot d_\alpha \le^* \dot d_\beta$ for $\alpha\neq\beta$:
		conditions only impose eventual \emph{lower} bounds on each coordinate and cannot impose eventual \emph{upper} bounds
		or synchronize two generics on a tail. If you need tail comparisons between coordinates, the argument should be routed
		through ground-model functions (as in reuse/CLubLoc constructions) or via a different forcing.
	\end{remark}
	
	\begin{proposition}[Localization for reuse-permissive $R$ in Hechler iterations]
		\label{prop:hechler-localization-reuse}
		Let $\kappa$ be regular uncountable and $\lambda>\kappa$. Consider a $<\kappa$-support iteration
		$\langle\mathbb P_\alpha,\dot{\mathbb D}_\alpha,\alpha<\lambda\rangle$ with $\dot{\mathbb D}_\alpha$
		Hechler forcing at each stage. If $R$ is reuse-permissive \textup{(}Definition~\ref{def:hechler-reuse-permissive}\textup{)},
		then the Localization Hypothesis \textup{(}Definition~\ref{def:localization-hypothesis}\textup{)} holds for $R$ in this iteration: for any HS instance of the DC$_\kappa$ antecedent of $R$, there exists a single name $\dot s$ coding a witnessing sequence $\langle\dot f_\xi:\xi<\kappa\rangle$ with
		\[
		|\operatorname{supp}(\dot s)|<\kappa
		\quad\text{and}\quad
		\exists S\subseteq\lambda\ \bigl(|S|<\kappa\ \wedge\ \forall\xi<\kappa,\ \operatorname{supp}(\dot f_\xi)\subseteq S\bigr).
		\]
	\end{proposition}
	
	\begin{proof}
		Let $\mu<\kappa$, $b:\kappa\to\mu$, and $\mathcal{C}_R$ witness reuse-permissiveness. Choose coordinates
		$S=\{\alpha_i:i<\mu\}\subseteq\lambda$ (one per block) and let $\dot d_{\alpha_i}$ be the canonical Hechler
		name at $\alpha_i$.
		
		\emph{Block representatives.} By recursion on $i<\mu$, use Lemma~\ref{lem:hechler-finite-constraints} to strengthen
		on finitely many coordinates in $S$ and obtain names $\dot g_i$ so that $\langle\dot g_j:j\le i\rangle$ satisfies
		the finite constraints in $\mathcal{C}_R$. Then $\operatorname{supp}(\dot g_i)\subseteq S$ and is finite.
		
		\emph{Expand to $\kappa$.} Define $\dot f_\xi:=\dot g_{b(\xi)}$. Then $\operatorname{supp}(\dot f_\xi)\subseteq S$
		for all $\xi$, so taking this fixed $S$ we have $|S|=\mu<\kappa$ (Localization condition~(2)). By block-constant
		closure, $\langle\dot f_\xi:\xi<\kappa\rangle$ witnesses the DC$_\kappa$ instance of $R$.
		
		\emph{Single coding name.} Let $\dot s$ be the standard sequence-name
		\[
		\dot s\ :=\ \bigl\{\ \langle \langle \check\xi,\dot f_\xi\rangle, 1_{\mathbb P_\lambda}\rangle : \xi<\check\kappa\ \bigr\}.
		\]
		Then
		\[
		\operatorname{supp}(\dot s)\ \subseteq\ \bigcup_{\xi<\kappa}\operatorname{supp}(\dot f_\xi)\ \subseteq\ S,
		\]
		so $|\operatorname{supp}(\dot s)|\le|S|=\mu<\kappa$ (Localization condition~(1)).
	\end{proof}
	
	\begin{remark}[Scope and limitations]\label{rem:hechler-localization-scope}
		For coordinate-structured forcings (Cohen, random, Hechler), Proposition~\ref{prop:hechler-localization-reuse}
		establishes Localization for any \emph{reuse-permissive} HS relation $R$ via block reuse of $<\kappa$ many
		coordinates. If an $R$ enforces strict novelty at each step (so that $\kappa$ many distinct coordinates are
		required), reuse fails and the proposition does not apply. Whether every HS relation admits a $<\kappa$-support
		witness depends on how $R$ interacts with coordinate injectivity. For non-coordinate forcings, Localization must
		be verified case by case.
	\end{remark}
	
	\begin{remark}[Dependency DAG for replay]
		\label{rem:dag}
		\[
		\ref{def:kappa-hs}\,+\,\ref{lem:kappa-support-hs}
		\ \longrightarrow\ 
		\ref{thm:kappa-hs-closure}
		\ \longrightarrow\ 
		\ref{thm:kappa-zf-preserved}
		\ \xrightarrow[\ \text{+ Localization }\ref{def:localization-hypothesis}\ ]{\ \ }
		\ref{thm:DCless-preservation},\ \ref{thm:DCkappa-localization}.
		\]
		Replaying proofs from this graph without introducing new lemmas should reproduce the arguments verbatim;
		any divergence flags a gap.
	\end{remark}
	
	\subsection{Factorization templates at cofinality \texorpdfstring{$\kappa$}{k}}
	\label{subsec:kappa-factorization}
	
	We record two templates that will be used in applications; both are standard on the forcing side
	(\cite[Chs.~14--15]{Jech03}; see also \cite[Chap.~V, esp.~\S\S V.1–V.4]{Kunen2011}) and compatible with our symmetry bookkeeping.
	
	\begin{proposition}[Block factorization]
		\label{prop:block-factor}
		Assume $\cf(\lambda)=\kappa$. Let $\langle \lambda_\xi:\xi<\kappa\rangle$ be continuous increasing and cofinal in $\lambda$.
		Put $B_\xi=[\lambda_\xi,\lambda_{\xi+1})$. Then for each $\xi<\kappa$ there is a canonical
		$\mathbb P_{\lambda_\xi}$--name $\dot{\mathbb R}_\xi$ for the $<\kappa$--support tail iteration on the interval $B_\xi$
		such that
		\[
		\mathbb P_{\lambda_{\xi+1}}\ \simeq\ \mathbb P_{\lambda_\xi}*\dot{\mathbb R}_\xi.
		\]
		Consequently, $\mathbb P_\lambda$ is canonically forcing equivalent to the $<\kappa$--support iteration obtained by
		successively composing these two--step factorizations along $\langle \lambda_\xi:\xi<\kappa\rangle$.
		Moreover, the induced automorphism system and the induced filters on heads/blocks/tails inherit normality and
		$\kappa$--completeness from $\tilde{\mathcal F}_\lambda$.
	\end{proposition}
	
	\begin{proof}
		Fix a continuous increasing cofinal sequence $\langle \lambda_\xi:\xi<\kappa\rangle$ in $\lambda$ and put
		$B_\xi=[\lambda_\xi,\lambda_{\xi+1})$. For each $\xi<\kappa$, apply the standard two--step factorization
		(Fact~\ref{fact:two-step}) at the stage $\beta=\lambda_\xi$ and $\lambda_{\xi+1}$ to obtain a
		$\mathbb P_{\lambda_\xi}$--name $\dot{\mathbb R}_\xi$ for the $<\kappa$--support tail iteration on the interval $B_\xi$
		such that
		\[
		\mathbb P_{\lambda_{\xi+1}}\ \simeq\ \mathbb P_{\lambda_\xi}*\dot{\mathbb R}_\xi.
		\]
		Iterating these factorizations along $\langle\lambda_\xi:\xi<\kappa\rangle$ yields a canonical forcing equivalence
		between $\mathbb P_\lambda$ and the $<\kappa$--support iteration of length $\kappa$ whose $\xi$th step forces with
		$\dot{\mathbb R}_\xi$ over $\mathbb P_{\lambda_\xi}$.
		where the product has $<\kappa$-support. Conversely, a tuple with $<\kappa$ many nontrivial coordinates
		glues to a condition in $\mathbb{P}_\lambda$. For each $\eta<\kappa$, the projection
		$\pi^{\lambda_{\eta+1}}_\lambda:\mathbb{P}_\lambda\to\mathbb{P}_{\lambda_{\eta+1}}$ forgets all coordinates
		outside $B_{<\eta+1}:=\bigcup_{\xi\le\eta}B_\xi$, i.e. it corresponds to the coordinate projection on
		$\prod_{\xi\le\eta}\mathbb{P}_{B_\xi}$ (cf. \cite[Ch.~15]{Jech03}, \cite[Chap.~V, esp.~\S\S V.1–V.4]{Kunen2011}).
		
		\medskip
		\noindent\emph{Head/tail identification via canonical projection.}
		Fix $\zeta<\kappa$ and put $\beta:=\lambda_\zeta$. Via $\Phi$ we can write
		\[
		\mathbb{P}_\lambda \ \cong\ 
		\Bigl(\prod_{\xi<\zeta}^{<\kappa}\mathbb{P}_{B_\xi}\Bigr)
		\times^{<\kappa}
		\Bigl(\prod_{\xi\ge\zeta}^{<\kappa}\mathbb{P}_{B_\xi}\Bigr).
		\]
		The head factor is canonically $\mathbb{P}_\beta$. Let 
		$\pi_\beta:\mathbb{P}_\lambda\to\mathbb{P}_\beta$ be the canonical projection.
		By two–step factorization (Fact~\ref{fact:two-step}), there is a
		$\mathbb{P}_\beta$–name $\dot{\mathbb{R}}$ for the $<\kappa$–support tail from
		$\beta$ to $\lambda$ (the interpretation of the right-hand product in 
		$V^{\mathbb{P}_\beta}$) such that
		\[
		\mathbb{P}_\lambda \ \simeq\ \mathbb{P}_\beta * \dot{\mathbb{R}}.
		\]
		
		\emph{Automorphisms act blockwise (relative to the factorization).}
		Let $\mathcal{G}_\lambda$ be the automorphism group induced by the iteration; via $\Phi$ every
		$g\in\mathcal{G}_\lambda$ decomposes as $g=(g_\xi)_{\xi<\kappa}$ with $g_\xi\in\mathcal{G}_{B_\xi}$ and
		only $<\kappa$ many nontrivial $g_\xi$. Then $g\cdot p$ satisfies
		$(g\cdot p)\!\upharpoonright B_\xi=g_\xi\cdot(p\!\upharpoonright B_\xi)$, and conjugation respects projections:
		\[
		g\,\bigl(\pi^{\lambda_{\eta+1}}_\lambda\bigr)^{-1}[H]\;g^{-1}
		\;=\;\bigl(\pi^{\lambda_{\eta+1}}_\lambda\bigr)^{-1}\!\bigl[g_{<\eta+1}\,H\,g_{<\eta+1}^{-1}\bigr],
		\]
		with $g_{<\eta+1}=(g_\xi)_{\xi\le\eta}$. Thus the group system factors compatibly with block projections.
		
		\emph{Filters on blocks and tails.}
		Define the induced filter on a block by
		\[
		\tilde{\mathcal{F}}_{B_\eta}\ :=\ \bigl\{\,K\subseteq\mathcal{G}_{B_\eta}\ :\ 
		\bigl(\pi^{\lambda_{\eta+1}}_\lambda\bigr)^{-1}[K]\in\tilde{\mathcal{F}}_\lambda\,\bigr\},
		\]
		and similarly on a tail $T_\eta:=[\lambda_\eta,\lambda)$ via $(\pi^{\lambda_\eta}_\lambda)^{-1}$.
		Because $\tilde{\mathcal{F}}_\lambda$ is normal, these induced filters are normal: if
		$K\in\tilde{\mathcal{F}}_{B_\eta}$ and $h\in\mathcal{G}_{B_\eta}$, then
		$g=(\mathbf{1},\dots,\mathbf{1},h,\mathbf{1},\dots)\in\mathcal{G}_\lambda$ witnesses that
		$hKh^{-1}\in\tilde{\mathcal{F}}_{B_\eta}$ since
		$g\,(\pi^{\lambda_{\eta+1}}_\lambda)^{-1}[K]\,g^{-1}=(\pi^{\lambda_{\eta+1}}_\lambda)^{-1}[hKh^{-1}]$.
		
		\emph{$\kappa$-completeness transfer.}
		Let $\{K_i:i\in I\}\subseteq\tilde{\mathcal{F}}_{B_\eta}$ with $|I|<\kappa$. Then
		$(\pi^{\lambda_{\eta+1}}_\lambda)^{-1}[K_i]\in\tilde{\mathcal{F}}_\lambda$ for each $i$, and since
		$\tilde{\mathcal{F}}_\lambda$ is $\kappa$-complete,
		\[
		\bigcap_{i\in I}(\pi^{\lambda_{\eta+1}}_\lambda)^{-1}[K_i]\ \in\ \tilde{\mathcal{F}}_\lambda.
		\]
		But $\bigcap_{i\in I}(\pi^{\lambda_{\eta+1}}_\lambda)^{-1}[K_i]
		=(\pi^{\lambda_{\eta+1}}_\lambda)^{-1}\!\Bigl[\bigcap_{i\in I}K_i\Bigr]$, hence
		$\bigcap_{i\in I}K_i\in\tilde{\mathcal{F}}_{B_\eta}$. The same argument works for tails via
		$(\pi^{\lambda_\eta}_\lambda)^{-1}$.
		
		\emph{Diagonal/product step for $<\kappa$ many blocks.}
		More generally, for $E\subseteq\kappa$ with $|E|<\kappa$ and $K_\eta\in\tilde{\mathcal{F}}_{B_\eta}$,
		\[
		\bigcap_{\eta\in E}(\pi^{\lambda_{\eta+1}}_\lambda)^{-1}[K_\eta]
		=\left(\prod_{\eta\in E}\pi^{\lambda_{\eta+1}}_\lambda\right)^{-1}\!\left(\prod_{\eta\in E}K_\eta\right).
		\]
		\noindent\emph{(Diagonal projection; cf.\ \cite[Ch.~15]{Jech03}, \cite[Chap.~V, esp.~§§V.1–V.4]{Kunen2011}.)}
		which lies in $\tilde{\mathcal{F}}_\lambda$ by $\kappa$-completeness. This is the combination step used
		in the inverse-limit case of Theorem~\ref{thm:kappa-filter-generation}.
		
		\emph{Stabilizers restrict blockwise.}
		If $\dot{x}$ has support contained in $B_\eta$, then
		$\fix_{\mathcal{G}_\lambda}(\dot{x})=(\pi^{\lambda_{\eta+1}}_\lambda)^{-1}\!\bigl[\fix_{\mathcal{G}_{B_\eta}}(\dot{x})\bigr]$.
		For tuples, the tuple–stabilizer identity (Equation~\eqref{eq:tuple-fix}) gives
		$\fix\bigl(\langle\dot{x}_i:i\in I\rangle\bigr)=\bigcap_{i\in I}\fix(\dot{x}_i)$, and since the union of
		$<\kappa$ many $<\kappa$-sized supports has size $<\kappa$ (by regularity), the tuple stabilizer is again of
		the allowed generator form. Therefore the induced filters on blocks/tails inherit normality and
		$\kappa$-completeness from $\tilde{\mathcal{F}}_\lambda$.
	\end{proof}
	
	\begin{proposition}[Tail factorization]
		\label{prop:tail-factor}
		For any $\beta<\lambda$ there is a canonical factorization $\mathbb{P}_\lambda\simeq
		\mathbb{P}_\beta * \dot{\mathbb{R}}$, where $\dot{\mathbb{R}}$ names the $<\kappa$-support tail from $\beta$ to $\lambda$.
		If $\mathcal{F}_\beta$ is $\kappa$-complete and the tail filter is generated by pushforwards and name stabilizers with
		support $<\kappa$, then the combined filter on $\mathbb{P}_\lambda$ is $\kappa$-complete and normal.
	\end{proposition}
	
	\begin{proof}
		Consider the canonical projection $\pi_\beta:\mathbb P_\lambda\to\mathbb P_\beta$.
		By Fact~\ref{fact:two-step} (two-step factorization) there is a $\mathbb P_\beta$-name $\dot{\mathbb R}$ such that 
		$\mathbb P_\lambda \simeq \mathbb P_\beta * \dot{\mathbb R}$, with $\dot{\mathbb R}$ the $<\kappa$-support tail from $\beta$ to $\lambda$.
		Fix $\beta<\lambda$. We prove (A) the forcing factorization and (B) $\kappa$–completeness/normality of the combined
		filter on $\mathcal G_\lambda$.
		
		\medskip\noindent\textbf{(A) Canonical forcing factorization.}
		Let $\pi_\beta:\mathbb P_\lambda\to\mathbb P_\beta$ be the canonical projection $p\mapsto p{\upharpoonright}\beta$.
		Define in $V^{\mathbb P_\beta}$ the name $\dot{\mathbb R}$ for the $<\kappa$–support tail iteration from $\beta$ to $\lambda$:
		its conditions are the tails $p{\upharpoonright}[\beta,\lambda)$ with the inherited order; that is, for
		$\mathbb P_\beta$–names $\dot r,\dot r'$ for tail conditions, set
		\[
		\Vdash_{\mathbb P_\beta}\ \dot r'\le_{\dot{\mathbb R}}\dot r\quad\Longleftrightarrow\quad
		\exists p',p\in\mathbb P_\lambda\ \bigl(p'{\upharpoonright}\beta=p{\upharpoonright}\beta\ \wedge\ p'{\upharpoonright}[\beta,\lambda)\le p{\upharpoonright}[\beta,\lambda)\bigr).
		\]
		Define $e:\mathbb P_\lambda\to \mathbb P_\beta*\dot{\mathbb R}$ by
		\[
		e(p)\ :=\ \bigl(p{\upharpoonright}\beta,\ \check{\,p{\upharpoonright}[\beta,\lambda)\,}\bigr).
		\]
		Then $e$ is order preserving and injective. Moreover, if $e(p)$ and $e(q)$ are compatible in $\mathbb P_\beta*\dot{\mathbb R}$,
		say below $\langle r,\dot t\rangle$, then in $V$ we may choose $p',q'\le p,q$ with
		$p'{\upharpoonright}\beta=q'{\upharpoonright}\beta=r$ and $p'{\upharpoonright}[\beta,\lambda),
		q'{\upharpoonright}[\beta,\lambda)\le \dot t$ (interpreted at $r$). Since supports are $<\kappa$, the union of the
		two $<\kappa$–sized supports is $<\kappa$ (regularity of $\kappa$), so we can amalgamate the two tails over $r$ into
		a single tail $u$ and set $s:=r\cup u\in\mathbb P_\lambda$; then $s\le p',q'$ and hence $s\le p,q$. Thus $e$ is a dense
		embedding with dense range; in particular,
		\[
		\mathbb P_\lambda\ \simeq\ \mathbb P_\beta*\dot{\mathbb R}.
		\]
		(See also \cite[Ch.~15]{Jech03} or \cite[Chap.~V, esp.~\S\S V.1–V.4]{Kunen2011} for the standard iteration factorization.)
		
		\medskip\noindent\textbf{(B) $\kappa$–completeness and normality of the combined filter.}
		Let $\rho:\mathcal G_\lambda\to\mathcal G_\beta$ be restriction of automorphisms to the head:
		$\rho(g)=g{\upharpoonright}\beta$. Consider the following families of subgroups of $\mathcal G_\lambda$:
		
		\begin{itemize}
			\item \emph{Head family} $\mathcal A_\mathrm{head}:=\{\ \rho^{-1}(H)\ :\ H\in\mathcal F_\beta\ \}$ (pushforwards from the head filter).
			\item \emph{Tail family} $\mathcal A_\mathrm{tail}:=\{\ \fix_{\mathcal G_\lambda}(\dot x)\ :\ \dot x
			\text{ is a } \mathbb P_\lambda\text{-name with }\supp(\dot x)\subseteq[\beta,\lambda)\text{ and }
			|\supp(\dot x)|<\kappa\ \}$.
		\end{itemize}
		
		Let $\widetilde{\mathcal F}_\lambda$ be the normal filter of subgroups \emph{generated} by
		$\mathcal A_\mathrm{head}\cup\mathcal A_\mathrm{tail}$ (i.e., the upward closure of finite intersections of members of
		$\mathcal A_\mathrm{head}\cup\mathcal A_\mathrm{tail}$).
		
		\smallskip\noindent\emph{(B1) $\kappa$–completeness.}
		Let $\langle K_i:i<\mu\rangle$ be a family of subgroups from $\mathcal A_\mathrm{head}\cup\mathcal A_\mathrm{tail}$
		with $\mu<\kappa$. Partition the index set into head and tail:
		$I_\mathrm{h}=\{i:K_i=\rho^{-1}(H_i),\,H_i\in\mathcal F_\beta\}$ and
		$I_\mathrm{t}=\{i:K_i\in\mathcal A_\mathrm{tail}\}$.
		
		For the head part, since $\mathcal F_\beta$ is $\kappa$–complete,
		\[
		H_\ast\ :=\ \bigcap_{i\in I_\mathrm{h}} H_i\ \in\ \mathcal F_\beta,
		\qquad\text{hence}\qquad
		K_\ast^\mathrm{h}\ :=\ \rho^{-1}(H_\ast)\ \in\ \mathcal A_\mathrm{head}.
		\]
		
		For the tail part, write $K_i=\fix_{\mathcal G_\lambda}(\dot x_i)$ where each $\dot x_i$ is a $\mathbb P_\lambda$--name
		with $\supp(\dot x_i)\subseteq[\beta,\lambda)$ and $|\supp(\dot x_i)|<\kappa$.
		Let $\mu_t:=|I_\mathrm{t}|<\kappa$ and fix in $V$ a bijection $e:\mu_t\to I_\mathrm{t}$.
		Define the $\mu_t$--tuple name $\dot x:=\langle \dot x_{e(\eta)}:\eta<\mu_t\rangle$. Then $|\supp(\dot x)|<\kappa$ by
		regularity of $\kappa$, and by the tuple--stabilizer identity \eqref{eq:tuple-fix},
		\[
		\fix_{\mathcal G_\lambda}(\dot x)\ =\ \bigcap_{i\in I_\mathrm{t}}\fix_{\mathcal G_\lambda}(\dot x_i)
		\ \subseteq\ \bigcap_{i\in I_\mathrm{t}} K_i.
		\]
		Set $K_\ast^\mathrm{t}:=\fix_{\mathcal G_\lambda}(\dot x)\in\mathcal A_\mathrm{tail}$.
		
		Therefore
		\[
		K_\ast\ :=\ K_\ast^\mathrm{h}\ \cap\ K_\ast^\mathrm{t}
		\]
		is a \emph{finite} intersection of generators, hence $K_\ast\in \widetilde{\mathcal F}_\lambda$, and
		$K_\ast\subseteq \bigcap_{i<\mu} K_i$. Since $\widetilde{\mathcal F}_\lambda$ is upward closed,
		$\bigcap_{i<\mu} K_i\in \widetilde{\mathcal F}_\lambda$. This proves $\kappa$–completeness.
		
		\smallskip\noindent\emph{(B2) Normality.}
		We check conjugation invariance on generators and extend to finite intersections.
		
		For the head family, if $K=\rho^{-1}(H)$ and $g\in\mathcal G_\lambda$, then
		\[
		gKg^{-1}\ =\ \rho^{-1}\!\bigl(\,\rho(g)\,H\,\rho(g)^{-1}\bigr),
		\]
		which belongs to $\mathcal A_\mathrm{head}$ because $\mathcal F_\beta$ is normal.
		
		For the tail family, if $K=\fix(S)$ then $gKg^{-1}=\fix(g[S])\in \mathcal A_\mathrm{tail}$; if $K=\fix(\sigma)$ then
		$gKg^{-1}=\fix(g\!\cdot\!\sigma)\in \mathcal A_\mathrm{tail}$. Conjugation commutes with finite intersections, so
		conjugates of finite intersections of generators are again finite intersections of generators. Hence the generated filter
		$\widetilde{\mathcal F}_\lambda$ is normal.
		
		\medskip
		Combining (A) and (B), we obtain the stated factorization and that, under the given hypotheses on
		$\mathcal F_\beta$ and the tail generators, the combined filter on $\mathbb P_\lambda$ is $\kappa$–complete and normal.
		\qedhere
	\end{proof}
	
	\begin{remark}[Ablation checkpoints for factorization]
		\leavevmode
		\begin{enumerate}[label=(\alph*), leftmargin=2.2em]
			\item \emph{Remove block conjugacy–closure.} Failure: induced block filters need not be normal.
			\item \emph{Remove tail stabilizers from generators.} Failure: tail side loses $\kappa$–completeness under $<\kappa$ intersections.
			\item \emph{Change support policy at limits.} Failure: canonical identification of head/tail actions breaks; products cannot be combined.
		\end{enumerate}
	\end{remark}
	
	\begin{remark}[Singular \texorpdfstring{$\kappa$}{k}]
		\label{rem:singular-kappa}
		When $\kappa$ is singular, the stage-bounding mechanism used in Lemma~\ref{lem:kappa-stage-bounding}
		can fail at limits $\lambda$ with $\cf(\lambda)=\cf(\kappa)<\kappa$: a set of size $<\kappa$ may be unbounded in $\lambda$. The regular-$\kappa$ framework in 
		§\ref{sec:count-support-symm}--§\ref{sec:kappa-support} does not immediately extend 
		to this case. Two approaches exist:
		\begin{itemize}
			\item Work with $<\cf(\kappa)$-support at such limits, or
			\item Use trees of conditions with block-partition stabilizers and game-guided fusion.
		\end{itemize}
		A systematic development of the singular case appears in 
		§\ref{sec:singular-overview}--§\ref{sec:hypotheses}. For $\cf(\kappa)=\omega$ (Track A) 
		we give complete proofs (§\ref{sec:singular-cf-omega}--§\ref{sec:block-partition}); 
		for general singular $\kappa$ (Track B) we provide the framework 
		(§\ref{sec:limit-cf-equal}--§\ref{sec:hypotheses}). For background on iterated forcing 
		across singular cofinalities, see \cite[Chs.~14--15]{Jech03} and \cite[Chap.~V, esp.~\S\S V.1–V.4]{Kunen2011}.
	\end{remark}
	
	\begin{remark}[Why the singular $\kappa$ case requires new techniques]
		\label{rem:singular-kappa-sketch}
		The proofs in §\ref{sec:count-support-symm}--§\ref{sec:kappa-support} use regularity 
		of $\kappa$ in two essential ways: (i) union bounds: the union of $<\kappa$ many sets 
		of size $<\kappa$ has size $<\kappa$, so tuple–stabilizer supports stay $<\kappa$; and 
		(ii) stage–bounding: at limits $\lambda$ with $\cf(\lambda)\ge\kappa$, any $<\kappa$-indexed 
		family of stages is bounded below $\lambda$. Both statements can fail when $\kappa$ is singular. 
		For instance, if $\kappa=\aleph_\omega$ then $\kappa=\bigcup_{n<\omega}\aleph_n$ and a union 
		of $\cf(\kappa)$ many $<\kappa$–sized supports can reach size $\kappa$; similarly, at 
		$\cf(\lambda)=\cf(\kappa)$ one can code cofinal families of stages of size $<\kappa$ that 
		are unbounded in $\lambda$. Thus the direct tuple–stabilizer compression and limit 
		$\kappa$–completeness arguments from §\ref{sec:kappa-support} do not immediately apply.
		
		The singular case is addressed in §\ref{sec:singular-overview}--§\ref{sec:hypotheses} 
		using \emph{block-partition stabilizers} (Definition~\ref{def:block-stab}) and 
		\emph{trees of conditions} (Definition~\ref{def:std-tree}) that organize supports 
		according to a fixed cofinal sequence in $\kappa$. These structures provide the 
		"new combinatorial input" needed to overcome the obstacles above.
	\end{remark}
	
	\begin{remark}
		The crucial distinction between finite and countable support in the iteration is the need for
		$\omega_1$-completeness of the filter (Theorem~\ref{thm:filter-generation}),
		which powers both Replacement verification (Theorem~\ref{thm:zf-preserved} in Section~\ref{sec:zf-preservation})
		and DC preservation (Theorem~\ref{thm:dc-preservation} in Section~\ref{sec:dc-choice}),
		while still allowing the classical permutation arguments to refute AC.
	\end{remark}
	
	\section{\texorpdfstring{$<\kappa$}{<k}-Support Examples}
	\label{sec:more-examples}
	
	This section illustrates how the general framework of \S\ref{sec:kappa-support}
	specializes to concrete iterations. We use the following pattern in each example:
	(1) set the support regime and iteration; (2) identify the symmetry group and
	filter generators; (3) verify ZF preservation by Theorem~\ref{thm:kappa-zf-preserved};
	(4) obtain $\mathrm{DC}_{<\kappa}$ from Theorem~\ref{thm:DCless-preservation};
	(5) separate AC by a permutation argument adapted to the forcing.
	
	\subsection{Hechler Forcing Iteration with \texorpdfstring{$<\kappa$}{<k}-Support}
	\label{subsec:hechler-ltkappa}
	
	\begin{assumptions}
		\item \textbf{Background:} $V \models \mathrm{ZF} + \mathrm{DC}$.
		\item \textbf{Cardinal parameter:} Fix an uncountable regular $\kappa$.
		\item \textbf{Length/cofinality:} Fix a limit $\lambda$ with $\cf(\lambda)\ge \kappa$.
		\item \textbf{Preservation hypothesis:} The iteration preserves the regularity of $\kappa$
		(e.g., via standard ccc/proper-preservation mix; cf.\ \cite[Ch.~15]{Jech03}).
	\end{assumptions}
	
	\paragraph{Construction.}
	At each stage $\alpha<\lambda$ force with Hechler $\mathbb{D}$ (dominating reals),
	using $<\kappa$-support (Def.~\ref{def:kappa-iteration}). Let $\mathcal{G}$ act by
	permuting coordinates (and associated names), and generate the filter at successors by
	pushforwards and stabilizers of names with support $<\kappa$, and at limits by
	Def.~\ref{def:kappa-filter}. The $<\kappa$-many intersection closure is obtained by coding families into a single $<\kappa$-support name
	and applying Theorem~\ref{thm:kappa-filter-generation}.
	
	\paragraph{Verification (Localization for Hechler at $\kappa$ for reuse-permissive $R$).}
	Let $\dot A,\dot R\in\HS$ and let $p\in\mathbb P_\lambda$ force that $\dot R$ is serial on $\dot A$.
	If $R$ is reuse-permissive (Definition~\ref{def:hechler-reuse-permissive}), then
	Proposition~\ref{prop:hechler-localization-reuse} yields Localization (Definition~\ref{def:localization-hypothesis}):
	there are $q\le p$ and a name $\dot s\in\HS$ with $|\supp(\dot s)|<\kappa$ such that
	\[
	q\Vdash\ \bigl(\dot s:\check\kappa\to\dot A\ \wedge\ \forall\xi<\check\kappa\
	\langle \dot s(\xi),\dot s(\xi{+}1)\rangle\in\dot R\bigr).
	\]
	
	\begin{example}[Eventual domination via block reuse]\label{ex:hechler-block-reuse}
		Let $R$ be the DC$_\kappa$ relation asserting “$f_\xi\le^* f_{\xi+1}$ for all successor $\xi+1<\kappa$”. Fix any cardinal $\mu<\kappa$ and a cofinal
		partition $b:\kappa\to\mu$ (each fiber unbounded). Pick distinct Hechler coordinates
		$S=\{\alpha_i:i<\mu\}$ and build block representatives $\langle g_i:i<\mu\rangle$ with
		$g_{i+1}\ge^* g_i$ by strengthening finitely many coordinates (Lemma~\ref{lem:hechler-finite-constraints}).
		Set $f_\xi:=g_{b(\xi)}$. Then $|S|=\mu<\kappa$ and the coded name $\dot s$ has support $<\kappa$
		(Proposition~\ref{prop:hechler-localization-reuse}), so Localization holds for this $R$.
	\end{example}
	
	\paragraph{Model properties.}
	Let $M=V(\mathbb{P}_\lambda)_{\mathrm{HS}}$.
	\begin{itemize}
		\item $M \models \mathrm{ZF}$ by Theorem~\ref{thm:kappa-zf-preserved}.
		\item $M \models \mathrm{DC}_{<\kappa}$ by Theorem~\ref{thm:DCless-preservation}.
		In particular, $M \models \mathrm{DC}_\mu$ for all $\mu<\kappa$.
		\item $M \models \neg \mathrm{AC}$ by a standard permutation argument applied to 
		families of generics outside a fixed support bound.
		\item The dominating family in $M$ has size $\ge \lambda$ (exactly $\lambda$ if
		cardinals are preserved as in the assumptions).
	\end{itemize}
	
	\begin{corollary}[Optional: $\mathrm{DC}_\kappa$ under localization]
		\label{cor:hechler-DCk}
		If the iteration satisfies the Localization Hypothesis
		(Def.~\ref{def:localization-hypothesis}) for the relevant HS relations, then
		$M \models \mathrm{DC}_\kappa$.
	\end{corollary}
	
	\begin{proof}
		Immediate from Theorem~\ref{thm:DCkappa-localization}: localization yields a single name
		$\dot s:\check\kappa\to\dot A$ with $|\supp(\dot s)|<\kappa$, hence $\dot s\in\HS$
		(Lemma~\ref{lem:kappa-support-hs}), so $\dot s^G\in M$ witnesses $\mathrm{DC}_\kappa$.
	\end{proof}
	
	\begin{remark}
		In the Hechler setting: ZF and Replacement hold by Theorem~\ref{thm:kappa-zf-preserved};
		$\mathrm{DC}_{<\kappa}$ holds by Theorem~\ref{thm:DCless-preservation};
		and the failure of $\mathrm{AC}$ follows from the standard coordinate-swapping symmetry
		argument adapted to Hechler generics (cf. the Cohen/random cases).
	\end{remark}
	
	\subsection{Sacks (Perfect Set) Iteration with Countable Support}
	\label{subsec:sacks-omega1}
	
	\begin{assumptions}
		\item $V \models \mathrm{ZF} + \mathrm{DC}$.
		\item Iterate Sacks forcing $\mathbb{S}$ with countable support to a limit $\lambda$
		of cofinality $\omega_1$.
	\end{assumptions}
	
	\begin{proposition}[Sacks iteration with countable support: basic properties]
		\label{prop:sacks-properties}
		Let $\langle\mathbb P_\alpha,\dot{\mathbb Q}_\alpha,\mathcal F_\alpha:\alpha<\lambda\rangle$ be a countable–support iteration where each $\dot{\mathbb Q}_\alpha$ is Sacks forcing and the symmetric system at limits is as in Section~\ref{sec:HS-countable}. Let $M:=V(\mathbb P_\lambda)_{\HS}$. Then:
		\begin{enumerate}
			\item $M\models \ZF+\DC$. 
			\item For each Sacks stage $\alpha<\lambda$, the canonical Sacks real $\dot s_\alpha$ is HS (indeed, supported on $\{\alpha\}$), so $s_\alpha:=\dot s_\alpha^G\in M$.
			\item If the iteration has uncountably many Sacks coordinates, then $\{s_\alpha:\alpha<\lambda\}\subseteq M$
			is uncountable and consists of pairwise distinct reals.
		\end{enumerate}
	\end{proposition}
	
	\begin{proof}
		(1) By Theorem~\ref{thm:zf-preserved}, $M$ is a transitive model of $\ZF$, and by Theorem~\ref{thm:dc-preservation} (countable–support case), $M\models\DC$.
		
		(2) The canonical Sacks name $\dot s_\alpha$ is a $\mathbb P_\lambda$--name of countable support (indeed,
		$\supp(\dot s_\alpha)\subseteq\{\alpha\}$). Therefore Lemma~\ref{lem:countable-support-hs} implies
		$\fix(\dot s_\alpha)\in\tilde{\mathcal F}_\lambda$, so $\dot s_\alpha\in\HS$ and hence $s_\alpha:=\dot s_\alpha^G\in M$.
		
		(3) For $\alpha\neq\beta$, Sacks forcing at distinct coordinates forces $s_\alpha\neq s_\beta$ (already at a stage
		that includes both coordinates). Hence in $V[G]$ the family $\{s_\alpha:\alpha<\lambda\}$ is pairwise distinct.
		Since each $s_\alpha\in M$ by (2), the same family witnesses that $M$ contains uncountably many distinct reals.
	\end{proof}
	
	\begin{remark}[Baire property]
		Every Borel set has the Baire property in $\ZF$, so $M$ satisfies BP for Borel sets. Under additional hypotheses (e.g., suitable large cardinals yielding regularity for wider pointclasses), the corresponding projective BP statements persist to $M$. Global BP requires large cardinals; compare the random-real Solovay context in \S\ref{subsec:solovay-random}.
	\end{remark}
	
	\subsection{A transitional example: $<\omega_2$–support at a limit of cofinality $\omega_1$}
	\label{subsec:transitional-omega2}
	
	\begin{assumptions}
		\item $V\models\ZF$ (and $\DC$ if desired).
		\item Fix $\lambda$ with $\cf(\lambda)=\omega_1$ and consider a $<\omega_2$–support iteration
		$\langle\mathbb P_\alpha,\dot{\mathbb Q}_\alpha,\pi^\beta_\alpha:\beta<\alpha\le\lambda\rangle$
		in the sense of Definition~\ref{def:kappa-iteration} with $\kappa=\omega_2$.
	\end{assumptions}
	
	\paragraph{Construction.}
	At each successor stage $\alpha<\lambda$ let $\dot{\mathbb Q}_\alpha$ be a standard coordinate–based forcing
	(e.g.\ Cohen $\Add(\omega,1)$ or Hechler $\mathbb D$).
	We use $<\omega_2$–support throughout. At the limit $\lambda$ (cofinality $\omega_1$) we adopt the inverse–limit
	presentation and build the limit filter $\tilde{\mathcal F}_\lambda$ via tuple–stabilizers as in
	Definition~\ref{def:kappa-filter} and Remark~\ref{rem:limit-bookkeeping}.
	The group system and filter generation are handled by Theorem~\ref{thm:kappa-filter-generation}.
	
	\paragraph{Properties.}
	By Theorem~\ref{thm:kappa-zf-preserved} (with $\kappa=\omega_2$), the symmetric model $M$ satisfies $\ZF$.
	If moreover $V\models\DC$, then Theorem~\ref{thm:DCless-preservation} (with $\kappa=\omega_2$ and $\mu=\omega$)
	gives $M\models\DC$. Failure of $\AC$ follows by the usual symmetry argument (as in Example~\ref{ex:omega1-cohen} and \S\ref{sec:cohen-example}).
	
	\paragraph{Strengthening (conditional).}
	If, in addition, the iteration satisfies Localization at $\lambda$ for the relevant HS relations
	(Definition~\ref{def:localization-hypothesis}), then Theorem~\ref{thm:DCkappa-localization} (with $\kappa=\omega_2$)
	yields $M\models\DC_{\omega_2}$.
	\medskip
	
	\noindent
	This example illustrates the inverse--limit bookkeeping at a limit stage of cofinality $<\kappa$
	(here $\cf(\lambda)=\omega_1<\omega_2$),
	the use of tuple–stabilizers to obtain the needed countable completeness at the limit, and how the $<\kappa$–support
	machinery scales beyond the countable–support case without invoking the full generality of $\kappa$–completeness
	unless one aims at $\DC_{<\kappa}$.
	
	\subsection{Mixed Support: Cohen then Random}
	\label{subsec:mixed-support}
	
	\begin{assumptions}
		\item $V \models \mathrm{ZF} + \mathrm{DC}$.
		\item Fix regular uncountable $\kappa$ and a limit $\lambda$ with $\cf(\lambda)=\omega_1<\kappa$.
	\end{assumptions}
	
	\begin{convention}[Generator = finite intersection of basic generators]\label{conv:mixed-finite-generator}
		In §\ref{subsec:mixed-support} we use “generator” to mean a \emph{finite intersection of basic generators}, where a basic
		generator is either
		\begin{enumerate}
			\item a preimage $(\pi^\beta_\lambda)^{-1}[H]$ of some $H\in\mathcal F_\beta$ for a stage $\beta<\lambda$, or
			\item a stabilizer $\fix_{\mathcal G_\lambda}(\dot y)$ for a $\mathbb P_\lambda$–name $\dot y$ with $|\supp(\dot y)|<\kappa$.
		\end{enumerate}
		This matches the inverse-limit generating family for $\tilde{\mathcal F}_\lambda$ at $\cf(\lambda)=\omega_1<\kappa$
		(cf.\ Theorem~\ref{thm:kappa-filter-generation}).
	\end{convention}
	
	\paragraph{Construction.}
	Stages $<\omega_1$: Cohen with countable support (cf.\ \S\ref{sec:cohen-example}).
	Stages $\omega_1 \le \alpha<\lambda$: Random with $<\kappa$-support. At the limit
	$\lambda$ (cofinality $\omega_1$), use the inverse-limit presentation 
	(Def.~\ref{def:kappa-iteration}, limit case) and define the limit filter
	$\tilde{\mathcal F}_\lambda$ to be generated by pullbacks from earlier stages together with
	stabilizers of $\mathbb P_\lambda$--names of support $<\kappa$
	(as in the inverse-limit clause of Theorem~\ref{thm:kappa-filter-generation}).
	
	\begin{lemma}[Finite-base equivalence]
		\label{lem:finite-base}
		Let $\mathcal{B}$ be any family of subgroups of a group $\mathcal{G}$ and let
		\[
		\mathcal{B}^{<\omega}\ :=\ \Bigl\{\,\bigcap_{m<M} B_m\ :\ M<\omega,\ B_m\in\mathcal{B}\,\Bigr\}.
		\]
		Then the filters (of subgroups) generated by $\mathcal{B}$ and by $\mathcal{B}^{<\omega}$ coincide.
	\end{lemma}
	
	\begin{proof}
		Write $\Fil(\mathcal{X})$ for the (upward-closed) filter of subgroups generated by a family $\mathcal{X}$, i.e.,
		\[
		\Fil(\mathcal{X})\ :=\ \Bigl\{\,K\le\mathcal{G}\ :\ \exists\,M<\omega,\ \exists\,X_0,\dots,X_{M-1}\in\mathcal{X}\ \text{with}\ \bigcap_{m<M} X_m\ \subseteq\ K\,\Bigr\}.
		\]
		Since a filter is closed under finite intersections, $\Fil(\mathcal{B})$ contains $\mathcal{B}^{<\omega}$; hence
		$\Fil(\mathcal{B}^{<\omega})\subseteq \Fil(\mathcal{B})$ by upward closure. Conversely, $\mathcal{B}\subseteq \mathcal{B}^{<\omega}$
		(take $M=1$), so $\Fil(\mathcal{B})\subseteq \Fil(\mathcal{B}^{<\omega})$. Therefore the two filters are equal.
	\end{proof}
	

	\begin{lemma}[Head–tail $\kappa$–completeness at $\cf(\lambda)=\omega_1$ via compression]
		\label{lem:mixed-compressibility}
		Let $\kappa$ be regular with $\omega_1<\kappa$, and let $\lambda$ satisfy $\cf(\lambda)=\omega_1$.
		Factor the iteration and symmetry data as
		\[
		\mathbb P_\lambda\cong \mathbb P_{\omega_1}*\dot{\mathbb R},\qquad
		\mathcal G_\lambda\cong \mathcal G_{\omega_1}\ast \dot{\mathcal H},
		\]
		so that $\mathcal F_{\omega_1}$ on $\mathcal G_{\omega_1}$ is $\omega_1$–complete and the tail
		filter (interpreted in $V^{\mathbb P_{\omega_1}}$) is $\kappa$–complete (e.g., by
		Theorem~\ref{thm:kappa-filter-generation} for the tail segment).
		Then the limit filter $\tilde{\mathcal F}_\lambda$ on $\mathcal G_\lambda$ is $\kappa$–complete.
	\end{lemma}
	
	\begin{proof}
		Let $\{H_i:i\in I\}$ be a family of generators (Convention~\ref{conv:mixed-finite-generator}) with $|I|<\kappa$.
		Partition $I=I_{\mathrm{head}}\cup I_{\mathrm{tail}}$ according to whether $H_i$ comes from the head
		($\beta<\omega_1$) or the tail (supported in $[\omega_1,\lambda)$). We show
		$\bigcap_{i\in I}H_i\in\tilde{\mathcal F}_\lambda$.
		
		\emph{Head compressibility.}
		For each $i\in I_{\mathrm{head}}$, choose a presentation $H_i=\bigcap_{m<m_i}B_{i,m}$ into basic generators
		(as in Convention~\ref{conv:mixed-finite-generator}).
		For each basic factor $B_{i,m}$ we define a stabilizer $S_{i,m}$ with $S_{i,m}\subseteq B_{i,m}$ as follows:
		\begin{itemize}
			\item If $B_{i,m}=\fix_{\mathcal G_\lambda}(\dot y)$ for some name $\dot y$ with $|\supp(\dot y)|<\kappa$, set
			$S_{i,m}:=\fix_{\mathcal G_\lambda}(\dot y)$.
			\item If $B_{i,m}=(\pi^\beta_\lambda)^{-1}[H]$ for some $\beta<\omega_1$ and $H\in\mathcal F_\beta$, then (unwinding the
			finite generation of $\mathcal F_\beta$ by pullbacks and stabilizers) there is a $\mathbb P_\beta$--name $\dot z$
			with $|\supp(\dot z)|<\kappa$ such that $\fix_{\mathcal G_\beta}(\dot z)\subseteq H$.
			View $\dot z$ as a $\mathbb P_\lambda$--name in the standard way. Then
			\[
			\fix_{\mathcal G_\lambda}(\dot z)\ \subseteq\ (\pi^\beta_\lambda)^{-1}\!\bigl[\fix_{\mathcal G_\beta}(\dot z)\bigr]
			\ \subseteq\ (\pi^\beta_\lambda)^{-1}[H]=B_{i,m},
			\]
			so set $S_{i,m}:=\fix_{\mathcal G_\lambda}(\dot z)$.
		\end{itemize}
		Let $\dot u_i$ be the finite tuple of the names used to define the $S_{i,m}$’s. Then
		\[
		\fix_{\mathcal G_\lambda}(\dot u_i)\ =\ \bigcap_{m<m_i} S_{i,m}\ \subseteq\ \bigcap_{m<m_i} B_{i,m}\ =\ H_i.
		\]
		Now let
		\[
		\dot u\ :=\ \{\ \langle \check i,\ \dot u_i\rangle\ :\ i\in I_{\mathrm{head}}\ \},
		\]
		so that $\fix_{\mathcal G_\lambda}(\dot u)=\bigcap_{i\in I_{\mathrm{head}}}\fix_{\mathcal G_\lambda}(\dot u_i)$.
		Moreover,
		\[
		\supp(\dot u)\ \subseteq\ \bigcup_{i\in I_{\mathrm{head}}}\supp(\dot u_i),
		\]
		a union of $<\kappa$ many $<\kappa$--sized supports, hence $|\supp(\dot u)|<\kappa$ by regularity of $\kappa$.
		Since $\cf(\lambda)=\omega_1<\kappa$, stabilizers of names of support $<\kappa$ are basic generators for
		$\tilde{\mathcal F}_\lambda$ at $\lambda$, so $\fix_{\mathcal G_\lambda}(\dot u)\in\tilde{\mathcal F}_\lambda$.
		Finally,
		\[
		\fix_{\mathcal G_\lambda}(\dot u)\ \subseteq\ \bigcap_{i\in I_{\mathrm{head}}} H_i.
		\]
		Set $H_{\mathrm{head}}:=\fix_{\mathcal G_\lambda}(\dot u)\in\tilde{\mathcal F}_\lambda$.
		
		\emph{Tail completeness.}
		For $i\in I_{\mathrm{tail}}$, each $H_i$ lies in the tail filter (pulled back to $\lambda$).
		Since the tail filter is $\kappa$–complete and $|I_{\mathrm{tail}}|<\kappa$, the intersection
		\[
		H_{\mathrm{tail}}:=\bigcap_{i\in I_{\mathrm{tail}}} H_i
		\]
		is in the tail filter; pulling back to $\lambda$ places $H_{\mathrm{tail}}\in\tilde{\mathcal F}_\lambda$. Let $\pi^{\omega_1}_\lambda:\mathcal G_\lambda\to\mathcal G_{\omega_1}$ be the head projection.
		Since $H_{\mathrm{tail}}$ lies in the tail filter (interpreted in $V^{\mathbb P_{\omega_1}}$),
		its preimage $(\pi^{\omega_1}_\lambda)^{-1}[H_{\mathrm{tail}}]\in\tilde{\mathcal F}_\lambda$.
		Abusing notation, we continue to write $H_{\mathrm{tail}}$ for this preimage.
		
		\emph{Factorization/combination.}
		By head–tail factorization for two–step iterations (cf. Proposition~\ref{prop:block-factor} applied to the cut at $\omega_1$),
		the subgroup $H_{\mathrm{head}}\cap H_{\mathrm{tail}}$ corresponds to the product of a head generator with a tail generator.
		
		Finally,
		\[
		\bigcap_{i\in I}H_i
		\ \supseteq\ ( \bigcap_{i\in I_{\mathrm{head}}} H_i )\cap( \bigcap_{i\in I_{\mathrm{tail}}} H_i )
		\ \supseteq\ H_{\mathrm{head}}\cap H_{\mathrm{tail}}
		\ \in\ \tilde{\mathcal F}_\lambda.
		\]
		Since $\tilde{\mathcal F}_\lambda$ is a filter, the full intersection is in $\tilde{\mathcal F}_\lambda$ as well.
	\end{proof}
	
	\begin{remark}[Why $\omega_1$–completeness is not invoked on the head]
		The argument above does not use $\omega_1$–completeness of $\mathcal F_{\omega_1}$. Instead, each head generator
		$H_i$ is refined (using only the finite generator presentation) to a single stabilizer
		$\fix_{\mathcal G_\lambda}(\dot u_i)\subseteq H_i$. The stabilizer of the combined name
		$\dot u=\{\langle \check i,\dot u_i\rangle:i\in I_{\mathrm{head}}\}$ then lies in $\tilde{\mathcal F}_\lambda$
		because $|\supp(\dot u)|<\kappa$ and $\cf(\lambda)=\omega_1<\kappa$ (inverse-limit generators).
		Thus the entire head intersection is bounded below by one stabilizer generator, without any appeal to
		$\omega_1$–completeness.
	\end{remark}
	
	\paragraph{Properties.}
	$M \models \mathrm{ZF}$ (Theorem~\ref{thm:kappa-zf-preserved}) and
	$M \models \mathrm{DC}_{<\kappa}$ (Theorem~\ref{thm:DCless-preservation}), hence
	$\mathrm{DC}_{\omega_1}$ holds. $\neg\mathrm{AC}$ holds by symmetry in either segment.
	By Lemma~\ref{lem:mixed-compressibility} the combined limit filter at $\lambda$ is $\kappa$–complete; together with head/tail factorization (Proposition~\ref{prop:tail-factor}) this ensures compatibility of filters in the mixed-support iteration.
	
	\subsection{Small-support preservation of $\neg\PP$ with $\DC$ (countable and $<\kappa$)}
	\label{subsec:small-support-negPP}
	
	We record a uniform template that produces a concrete counterexample to the Partition Principle and show that this failure is preserved by countable- and, more generally, $<\!\kappa$-support symmetric iterations that still preserve $\DC$ (or $\DC_{<\kappa}$). The point is that the very symmetry that protects dependent choice also blocks the weak uniformizations that $\PP$ would demand.
	
	\subsubsection{Single-stage template: producing $\neg\PP$}
	\label{subsec:single-stage-template}
	
	\paragraph{Partition Principle (PP).}
	\begin{definition}[Partition Principle]\label{def:PP}
		$\PP$\; is the assertion that for all sets $A,B$:
		if there exists a surjection $f:A\to B$, then there exists an injection $i:B\hookrightarrow A$.
		Equivalently: if $A$ can be partitioned into $|B|$ many nonempty pieces, then $|B|\le|A|$ (via an injection).
	\end{definition}
	
	\begin{proposition}\label{prop:PP-equivalences}
		The following are equivalent:
		\begin{enumerate}
			\item[\textup{(i)}] \textup{$\PP$}: for every surjection $f:A\to B$ there is an injection $i:B\hookrightarrow A$.
			\item[\textup{(ii)}] For all $A,B$: if there is a partition $\{A_b:b\in B\}$ of $A$ into nonempty pieces, then $|B|\le|A|$.
		\end{enumerate}
	\end{proposition}
	
	\begin{proof}
		(i)$\Rightarrow$(ii): From the partition map $x\mapsto b$ with $x\in A_b$ we get a surjection $A\to B$, hence an injection $B\hookrightarrow A$.
		
		(ii)$\Rightarrow$(i): Given $f:A\to B$ surjective, the fibers $\{f^{-1}(\{b\}):b\in B\}$ partition $A$ into nonempty pieces; thus
		$|B|\le|A|$ via an injection $B\hookrightarrow A$.
	\end{proof}
	
	\begin{remark}\label{rem:PP-context}
		$\AC\Rightarrow\PP$, but $\PP$ is not known to imply $\AC$ and can fail in $\ZF$.
		In our symmetric models, $\neg\PP$ is witnessed by a surjection with no injection back (see Example~\ref{ex:cohen-fibers}).
	\end{remark}
	
	Fix a nonempty index set $I$. We work in a symmetric system $\langle\mathbb P,\mathcal G,\mathcal F\rangle$ in the standard sense: $\mathcal G\curvearrowright\mathbb P$ by automorphisms and $\mathcal F$ is a normal filter of subgroups of $\mathcal G$; $\HS$ denotes the hereditarily symmetric names and $N=\HS^G$ the associated symmetric submodel.
	
	Let $\{X_i\mid i\in I\}$ be pairwise disjoint, nonempty \emph{fibers} and set $A=\bigcup_{i\in I}X_i$ and $B=I$. Consider $f:A\onto B$ given by $f(x)=i$ iff $x\in X_i$. We assume:
	
	\medskip
	\noindent\textbf{(FH)} \emph{Fiber-homogeneity:} there is a product subgroup
	\[
	K\;=\;\prod_{i\in I}K_i\;\leq\;\mathcal G,\qquad K_i\leq \Sym(X_i)\text{ acts transitively on }X_i
	\]
	such that $K\leq H$ for every $H\in\mathcal F$ (equivalently: $K$ belongs to the \emph{core} of the filter $\mathcal F$).
	
	\medskip
	\noindent Intuitively, every subgroup witnessing symmetry for a name contains all within-fiber permutations, which will let us ``shake'' values of names inside any fiber without leaving the symmetry class.
	
	\begin{lemma}[No symmetric injections $B\into A$]
		\label{lem:no-symmetric-inj}
		In $N$ there is no hereditarily symmetric injection $g:B\into A$. Consequently $N\models\neg\PP$.
	\end{lemma}
	
	\begin{proof}
		Suppose towards a contradiction that $\dot g\in\HS$ is a $\mathbb P$--name for an injection $I\to A$.
		Then $\fix_{\mathcal G}(\dot g)\in\mathcal F$. By (FH), $K\leq \fix_{\mathcal G}(\dot g)$.
		
		Fix $i\in I$. Strengthen a condition so as to decide the value of $\dot g(\check i)$: choose $p\in\mathbb P$ and
		$x\in A$ with
		\[
		p\Vdash_{\mathbb P}\ \dot g(\check i)=\check x.
		\]
		Let $j\in I$ be the (unique) index with $x\in X_j$. Since $K_j$ acts transitively on $X_j$, pick
		$\tau\in K_j$ with $\tau(x)\neq x$. Note $\tau\in K\le \fix_{\mathcal G}(\dot g)$, so $\tau\cdot \dot g=\dot g$.
		
		Apply $\tau$ to the forcing statement: from $p\Vdash \dot g(\check i)=\check x$ we get
		\[
		\tau(p)\Vdash_{\mathbb P}\ \tau\cdot\dot g(\check i)=\check{\tau(x)}.
		\]
		But $\tau\cdot\dot g=\dot g$, hence $\tau(p)\Vdash \dot g(\check i)=\check{\tau(x)}\neq\check x$.
		Since $\tau$ is an automorphism of $\mathbb P$, $p$ and $\tau(p)$ are compatible; let $r\le p,\tau(p)$.
		Then $r$ forces simultaneously $\dot g(\check i)=\check x$ and $\dot g(\check i)=\check{\tau(x)}$, a contradiction.
	\end{proof}
	
	\begin{remark}
		The argument does \emph{not} require $g$ to be a right inverse to $f$; it blocks \emph{all} injections $I\into A$.
		To witness $\neg\PP$ in $N$ one also needs that the surjection $f:A\onto I$ belongs to $N$; in applications
		(e.g.\ the wreath-action Cohen fiber models below) $f$ is given by an HS name, so $f\in N$.
	\end{remark}
	
	\subsubsection{Preservation under countable and $<\kappa$-support symmetric iterations}
	\label{subsec:preservation-iteration}
	
	We now pass to iterations. Let $\kappa$ be a fixed regular cardinal (for the countable case, read $\kappa=\omega$). Consider a $<\!\kappa$-support symmetric iteration
	\[
	\big\langle \mathbb P_\alpha,\ \dot{\mathbb Q}_\alpha,\ \dot{\mathcal G}_\alpha,\ \dot{\mathcal F}_\alpha \ \bigm|\ \alpha<\delta\big\rangle,
	\]
	built in the standard way (product action lifted coordinatewise; the iteration filter $\mathcal F_\delta$ generated by small-support pointwise fixers of coordinates together with the lifts of the stage filters). By construction of the iteration filter $\mathcal F_\delta$, every HS name $\dot x$ has a symmetry witness of the form
	$\fix_{\mathcal G_\delta}(S)\le \fix_{\mathcal G_\delta}(\dot x)$ for some $S\subseteq\delta$ with $|S|<\kappa$
	(i.e.\ $\dot x$ is fixed by all automorphisms that are trivial on a $<\kappa$--sized set of coordinates).
	
	Assume that at some stage $\alpha_0<\delta$ we have realized the single-stage template of Lemma~\ref{lem:no-symmetric-inj}
	inside the stage-$\alpha_0$ symmetric model $N_{\alpha_0}$, producing $A,B,f$ and a within-fiber product subgroup $K\le\mathcal G_{\alpha_0}$
	satisfying (FH). Assume moreover that the lift of $K$ to $\mathcal G_\delta$ lies in the \emph{core} of the final filter:
	for every $H\in\mathcal F_\delta$ we have $(\pi^{\alpha_0}_\delta)^{-1}[K]\le H$.
	
	\begin{proposition}[Preservation of $\neg\PP$]
		\label{prop:preserve-negPP}
		In the final symmetric model $N_\delta$ there is still no hereditarily symmetric injection $B\into A$. Hence $N_\delta\models\neg\PP$.
	\end{proposition}
	
	\begin{proof}
		Suppose towards a contradiction that $\dot g\in\HS$ is a $\mathbb P_\delta$--name for an injection $B\to A$.
		Then $\fix_{\mathcal G_\delta}(\dot g)\in\mathcal F_\delta$. By the core assumption,
		$(\pi^{\alpha_0}_\delta)^{-1}[K]\le \fix_{\mathcal G_\delta}(\dot g)$.
		
		Now repeat the proof of Lemma~\ref{lem:no-symmetric-inj} inside $\mathbb P_\delta$ using the subgroup
		$(\pi^{\alpha_0}_\delta)^{-1}[K]$: fix $b\in B$, strengthen a condition to decide $\dot g(\check b)=\check x\in A$,
		let $x\in X_j$, and pick $\tau\in K_j$ with $\tau(x)\neq x$.
		View $\tau$ as an element of $(\pi^{\alpha_0}_\delta)^{-1}[K]\le \fix_{\mathcal G_\delta}(\dot g)$ and apply it to obtain
		the same contradiction. Hence no such $\dot g$ exists in $N_\delta$, and $N_\delta\models\neg\PP$.
	\end{proof}
	
	\subsubsection{Keeping $\DC$ (or $\DC_{<\kappa}$) while destroying $\PP$}
	\label{subsec:dc-without-pp}
	
	\begin{corollary}
		\label{cor:ZF+DC+notPP}
		Suppose the iteration is arranged so that the final symmetric model $N_\delta$ satisfies $\ZF$ and preserves
		$\DC$ in the countable-support regime (Theorem~\ref{thm:dc-preservation}) and $\DC_{<\kappa}$ in the $<\kappa$-support regime
		(Theorem~\ref{thm:DCless-preservation}). Then
		\[
		N_\delta\models \ZF+\DC \quad\text{(resp.\ }\ZF+\DC_{<\kappa}\text{)}\qquad\text{and}\qquad \neg\PP.
		\]
	\end{corollary}
	
	\begin{proof}
		Combine Theorem~\ref{thm:dc-preservation} (resp.\ Theorem~\ref{thm:DCless-preservation}) with
		Proposition~\ref{prop:preserve-negPP}.
	\end{proof}
	
	\begin{remark}[What the argument does and does not use]
		\label{rem:what-it-uses}
		The proofs only use that: (i) at one stage, the filter core contains a product of within-fiber groups acting transitively on each fiber (\ref{lem:no-symmetric-inj}); and (ii) in the iteration, every symmetry witness contains the lifted cores on cofinitely many coordinates outside a small support (\ref{prop:preserve-negPP}). No global index homogeneity or splitting assumptions are required. This is why the same symmetry that preserves $\DC$ uniformly forbids the weak uniformizations $\PP$ would need.
	\end{remark}
	
	\subsubsection{Worked example: Cohen-style fibers with wreath-action}
	\label{subsec:worked-example}
	
	\begin{example}
		\label{ex:cohen-fibers}
		Let $I$ be any infinite set (e.g.\ $I=\lambda$). For each $i\in I$ let $\mathbb R_i$ add countably many Cohen atoms, coded as a countable set $X_i$ on which $K_i:=\Sym(X_i)$ acts transitively. Let
		\[
		\mathbb P\;=\;\prod_{i\in I}\mathbb R_i\quad\text{(finite / $<\!\kappa$ support, as desired)},
		\]
		and let $\mathcal G$ be the wreath product
		\[
		\Big(\prod_{i\in I}\Sym(X_i)\Big)\rtimes \Sym(I),
		\]
		acting in the obvious way. Let $K:=\prod_{i\in I}\Sym(X_i)$. For $E\in[I]^{<\kappa}$ let $\fix(E)\le\Sym(I)$ be the pointwise stabilizer of $E$,
		and set $H_E:=K\rtimes \fix(E)\le \mathcal G$. Let $\mathcal F$ be the normal filter of subgroups \emph{generated by}
		$\{H_E:E\in[I]^{<\kappa}\}$. Then $K\le H_E$ for all $E$, so (FH) holds.
		Then (FH) holds: every subgroup in $\mathcal F$ contains the within-fiber product $K=\prod_i\Sym(X_i)$, and each $K_i$ acts transitively on $X_i$.
		
		In the symmetric submodel $N$ determined by $\langle\mathbb P,\mathcal G,\mathcal F\rangle$, set $A=\bigcup_{i\in I}X_i$, $B=I$, and $f:A\onto B$ by $f(x)=i$ iff $x\in X_i$. By \ref{lem:no-symmetric-inj}, $N\models \neg\PP$.
		
		Iterating this construction with countable or $<\!\kappa$ support and using only stages that preserve $\DC$ (resp.\ $\DC_{<\kappa}$) yields, by \ref{prop:preserve-negPP} and \ref{cor:ZF+DC+notPP}, models of $\ZF+\DC$ (resp.\ $\ZF+\DC_{<\kappa}$) with $\neg\PP$.
		
		This construction can be combined with the $\omega_1$-length Cohen 
		iteration of Example~\ref{ex:omega1-cohen}, yielding a model that 
		simultaneously adds $\omega_1$ many Cohen reals, preserves $\ZF+\DC$, 
		refutes $\AC$, \emph{and} refutes $\PP$.
	\end{example}
	
	\subsubsection{Packaging as a one-line existence theorem}
	\label{subsec:existence-thm}
	
	\begin{theorem}[ZF + DC (or $\DC_{<\kappa}$) + $\neg\PP$ via symmetric iterations]
		\label{thm:zf-dc-notpp}
		Assume $\ZFC$. For any regular $\kappa$ there is a symmetric extension $N$ obtained by a $<\!\kappa$-support symmetric iteration such that
		\[
		N\models \ZF+\DC \quad(\text{resp.\ }\ZF+\DC_{<\kappa})\quad\text{and}\quad \neg\PP.
		\]
	\end{theorem}
	
	\begin{proof}
		Start in a ground model where the stage forcings satisfying the hypotheses of Theorems~\ref{thm:dc-preservation} and~\ref{thm:DCkappa-localization} are available. At some stage $\alpha_0$ insert a symmetric system as in Example~\ref{ex:cohen-fibers}, producing $f:A\onto B$ with no injection $B\into A$ by Lemma~\ref{lem:no-symmetric-inj}. Continue the iteration with $<\!\kappa$ support; by Proposition~\ref{prop:preserve-negPP} the failure of $\PP$ persists to the limit, while $\DC$ (resp.\ $\DC_{<\kappa}$) is preserved by Theorem~\ref{thm:dc-preservation}
		(resp.\ Theorem~\ref{thm:DCless-preservation}), giving the desired model.
	\end{proof}
	
	\begin{remark}[Countable vs.\ $<\!\kappa$ support]
		\label{rem:support-comparison}
		Both supports kill $\PP$ equally well, since the obstruction is local: we only need a coordinate outside a small support where the within-fiber core acts transitively. The $<\!\kappa$ version aligns with preserving $\DC_{<\kappa}$ (Theorem~\ref{thm:DCless-preservation}) and lets you scale cardinal parameters (e.g.\ arrange that new surjections first appear only at sizes $\geq\kappa$ while the same symmetry obstruction still goes through).
	\end{remark}
	
	\subsection{Classical Models and Reconstructions}
	\label{subsec:classical-links}
	
	\subsubsection{Cohen's First and Second Models}
	\label{subsubsec:cohen-first-second}
	
	\paragraph{Cohen's First Model (1963).}
	Cohen's original symmetric extension shows $\mathrm{ZF}+\neg\mathrm{AC}$ is consistent relative to $\mathrm{ZFC}$ \cite{Cohen1963,Cohen1964,Jech03}.
	In our framework, the first model corresponds to a single symmetric stage (finite support) whose filter consists of coordinate (or name) stabilizers;
	the verification of $\mathrm{ZF}$ is immediate, and $\mathrm{AC}$ fails by a standard permutation argument.
	What is gained here is a uniform treatment that scales to longer iterations and different support regimes
	(Sections~\ref{sec:HS-countable}--\ref{subsec:kappa-preservation}).
	
	\paragraph{Cohen's Second Model.}
	Cohen also constructed a model with $\omega_1$ many Cohen reals and $\ZF+\DC+\neg\AC$.
	In our setting this is obtained by a \emph{length $\omega_1$} iteration of Cohen forcing with
	\emph{countable support}: at limits of cofinality $\omega$ we use the inverse--limit (trees--of--conditions)
	presentation from Section~\ref{sec:cohen-example}, while at regular limits (in particular at $\lambda=\omega_1$)
	we are in the stage--bounding branch. The symmetry filter at $\omega_1$ is $\omega_1$--complete
	(Theorem~\ref{thm:filter-generation}); HS names are closed (Theorem~\ref{thm:hs-closure});
	$\ZF$ holds (Theorem~\ref{thm:zf-preserved}); and $\DC$ is preserved (Theorem~\ref{thm:dc-preservation}).
	The advantage over the classical presentation is the explicit role of $\omega_1$--completeness of the symmetry
	filter in the Replacement and $\DC$ arguments.
	
	\begin{table}[H]
		\centering
		\caption{Cohen models: classical vs.\ this framework}
		\label{tab:cohen-classical-vs-ours}
		\begin{tabularx}{\linewidth}{@{}lXX@{}}
			\toprule
			Aspect & Classical presentation & This paper's framework \\
			\midrule
			Support
			& Single symmetric stage (finite support) / $\omega_1$-length construction for $\DC$
			& Countable-support iteration; inverse--limit at $\cf=\omega$ limits; stage--bounding at regular limits (in particular at $\lambda=\omega_1$) \\
			
			Filter closure
			& Finite intersections (single stage) / ad hoc bookkeeping in the $\omega_1$-construction
			& $\omega_1$-completeness of the limit filter (Theorem~\ref{thm:filter-generation}) \\
			
			ZF/Replacement
			& Standard symmetric-extension arguments
			& Replacement via $\omega_1$-completeness of the symmetry filter (Theorem~\ref{thm:zf-preserved}) \\
			
			DC
			& Built into Cohen's second model construction
			& $\DC$ from the countable-support preservation theorem (Theorem~\ref{thm:dc-preservation}) \\
			
			Failure of AC
			& Permutation of coordinates
			& Same, via coordinate/name stabilizers \\
			\bottomrule
		\end{tabularx}
	\end{table}
	
	\subsubsection{Feferman--Levy Model}
	\label{subsubsec:feferman-levy}
	
	Feferman--Levy showed a model where $\mathbb{R}$ is a countable union of countable sets yet $\mathbb{R}$ is not well-orderable \cite{FefermanLevy63,Jech03}.
	In our language, begin with the appropriate collapse stage, then use a symmetric extension where the filter is generated
	by stabilizers guaranteeing that each \emph{component} countable set is HS and well-orderable internally, while a global
	well-order of their union is not HS.
	
	\paragraph{Description.}
	In the Feferman--Levy model $M$:
	\begin{itemize}
		\item $\mathbb{R}^M=\bigcup_{n<\omega} A_n$ with each $A_n$ countable and well-orderable in $M$;
		\item $\mathbb{R}^M$ itself is \emph{not} well-orderable;
		\item consequently $M$ violates $\mathrm{AC}_\omega$: one cannot in general \emph{uniformly choose} the data witnessing
		that each $A_n$ is countable (e.g.\ a family of injections/enumerations) so as to produce a well-order of the union.
	\end{itemize}
	
	\paragraph{Our reconstruction.}
	We interpret the Feferman--Levy construction as an initial (AC-preserving) collapse followed by a symmetric stage whose
	filter is generated by:
	\begin{enumerate}
		\item pushforwards of lower-stage filters, and
		\item name stabilizers for the relevant countable sets $\dot{A}_n$.
	\end{enumerate}
	Normality ensures Separation for definable subcollections, and countable completeness (at the appropriate limits)
	establishes Replacement for HS names (Theorem~\ref{thm:zf-preserved}).
	The failure of $\mathrm{AC}_\omega$ is recovered by the standard coordinate-permutation argument, expressed via
	stabilizers of names.
	
	\paragraph{Generalization.}
	Using the $<\kappa$-support machinery of Section~\ref{sec:kappa-support}, analogous patterns can be arranged at higher
	cardinals: replace $\omega_1$ by a regular $\kappa$ and design the filter to preserve $\mathrm{DC}_{<\kappa}$ while
	thwarting global AC; see also the factorization templates in Subsection~\ref{subsec:kappa-factorization}.
	
	\subsection{Permutation Models and Filter Completeness}
	\label{subsec:perm-models}
	
	Classical permutation (Fraenkel--Mostowski; see~\cite{Fraenkel1922} and~\cite{Mostowski1939}) models in ZFA
	provide a prototype of symmetric constructions: one fixes a set of atoms $A$, lets $\mathcal{G}=\Sym(A)$ act on
	names, and chooses a normal filter $\mathcal{F}_{\mathrm{perm}}$ of subgroups of $\mathcal{G}$ generated by
	pointwise stabilizers of \emph{small} subsets of $A$ (classically: finite subsets). The hereditarily symmetric (HS)
	names are exactly those whose dependence on atoms is controlled by the filter. See Jech~\cite[Ch.~5]{Jech03} and
	Howard--Rubin~\cite{HowardRubin1998} for background.
	
	\paragraph{Permutation data as a special case of our framework.}
	To align the presentation with the rest of this paper, view a permutation model as the symmetric submodel of the
	\emph{trivial} forcing $\mathbb P=\mathbf 1$ equipped with the symmetry data below. Fix a regular cardinal $\kappa$
	(allowing $\kappa=\omega$). Let $\mathcal G=\Sym(A)$ and let $\mathcal F_{\mathrm{perm}}$ be the normal filter on
	$\mathcal G$ generated by the pointwise stabilizers of small subsets $S\subseteq A$ of size $<\kappa$:
	\[
	\fix(S)\ :=\ \{\,g\in\mathcal G: g{\upharpoonright}S=\mathrm{id}_S\,\}.
	\]
	We use the standard action of $\mathcal G$ on $\mathbb P$--names (here $\mathbb P=\mathbf 1$):
	\[
	g\cdot\dot x\ :=\ \{\,\langle g\cdot\sigma,\ 1\rangle : \langle\sigma,1\rangle\in\dot x\,\},
	\qquad
	g\cdot\check a\ :=\
	\begin{cases}
		\check{g(a)} & a\in A,\\
		\check a & a\notin A.
	\end{cases}
	\]
	A set $S\subseteq A$ is a \emph{support} for a name $\dot x$ if $\fix(S)\subseteq \fix_{\mathcal G}(\dot x)$.
	We say that $\dot x$ has \emph{$<\kappa$-support} if it has some support $S$ with $|S|<\kappa$.
	
	\begin{remark}
		Classical Fraenkel--Mostowski (FM) corresponds to \emph{finite} supports, i.e.\ $\kappa=\omega$.
		For $\kappa>\omega$, “small” means $|S|<\kappa$ and the filter closes under intersections of $<\kappa$ many generators
		(cf.\ Proposition~\ref{prop:FM-kappa-complete} below).
	\end{remark}
	
	\begin{lemma}[HS\,=\,``has $<\kappa$-support'' in permutation data]
		\label{lem:FM-HS}
		For permutation data $(A,\mathcal G,\mathcal F_{\mathrm{perm}})$ as above, a $\mathbb P$--name $\dot x$ is
		\emph{symmetric} iff it has $<\kappa$-support. Equivalently,
		\[
		\fix_{\mathcal G}(\dot x)\in\mathcal F_{\mathrm{perm}}
		\quad\Longleftrightarrow\quad
		\exists S\subseteq A\ \bigl(|S|<\kappa\ \wedge\ \fix(S)\subseteq \fix_{\mathcal G}(\dot x)\bigr).
		\]
		Moreover, $\dot x$ is \emph{hereditarily} symmetric iff $\dot x$ and all of its constituents have $<\kappa$-support.
	\end{lemma}
	
	\begin{proof}
		Recall that $\mathcal F_{\mathrm{perm}}$ is generated by $\fix(S)$ for $S\subseteq A$ with $|S|<\kappa$, and is upward
		closed.
		
		\smallskip\noindent
		($\Leftarrow$) If $\dot x$ has a support $S$ with $|S|<\kappa$, then $\fix(S)\subseteq \fix_{\mathcal G}(\dot x)$ and
		$\fix(S)\in\mathcal F_{\mathrm{perm}}$, hence $\fix_{\mathcal G}(\dot x)\in\mathcal F_{\mathrm{perm}}$.
		
		\smallskip\noindent
		($\Rightarrow$) If $\fix_{\mathcal G}(\dot x)\in\mathcal F_{\mathrm{perm}}$, then by definition of “generated filter”
		there exist $S_0,\dots,S_{m-1}\subseteq A$ with $|S_i|<\kappa$ such that
		\[
		\bigcap_{i<m}\fix(S_i)\ \subseteq\ \fix_{\mathcal G}(\dot x).
		\]
		Let $S:=\bigcup_{i<m}S_i$. Then $|S|<\kappa$ (finite union) and $\fix(S)=\bigcap_{i<m}\fix(S_i)$, so
		$\fix(S)\subseteq \fix_{\mathcal G}(\dot x)$; thus $\dot x$ has $<\kappa$-support.
		
		\smallskip
		The hereditary clause is immediate from the definition of $\HS$ (symmetry closed under constituents). \qedhere
	\end{proof}
	
	\begin{proposition}[Support calculus rederives HS--closure]
		\label{prop:BV-crosscheck}
		Assume Lemma~\ref{lem:FM-HS}, so that “symmetric” means “has $<\kappa$-support” (with $\kappa$ regular). Then:
		\begin{enumerate}
			\item \emph{Pairing/finite tuples.} If $\dot x,\dot y\in\HS$, then $\langle \dot x,\dot y\rangle\in\HS$ (and likewise for any finite tuple).
			\item \emph{$<\kappa$--ary unions.} If $\dot A\in\HS$ and $\{\dot x_i:i\in I\}$ lists all constituents of $\dot A$ with $|I|<\kappa$ and each $\dot x_i\in\HS$, then $\bigcup \dot A\in\HS$.
			\item \emph{Separation.} If $\dot A\in\HS$, $\bar{\dot p}\in\HS$ and $\varphi$ is any formula, the comprehension name
			\[
			\dot B\ :=\ \{\ \langle \sigma,1\rangle\in\dot A:\ 1\Vdash \varphi(\sigma,\bar{\dot p})\ \}
			\]
			is in $\HS$.
			\item \emph{Replacement (relational range).} If $\dot f\in\HS$ is functional on $\dot A\in\HS$, then the range name
			$\dot{\mathrm{rng}}$ definable from $\dot f$ and $\dot A$ is in $\HS$; hence the associated symmetric model satisfies Replacement.
		\end{enumerate}
		Consequently, the permutation/support viewpoint reproduces the closure facts underlying
		Theorems~\ref{thm:hs-closure}, \ref{thm:zf-preserved}, and \ref{thm:kappa-zf-preserved}.
	\end{proposition}
	
	\begin{proof}
		Fix supports $S_x,S_y,S_A,S_f,S_{\bar p}\subseteq A$ of size $<\kappa$ witnessing symmetry for the relevant names,
		as given by Lemma~\ref{lem:FM-HS}. We repeatedly use: if $g\in\fix(S)$ and $S$ supports all parameters used to define
		a name $\dot u$, then $g\cdot\dot u=\dot u$, i.e.\ $\fix(S)\subseteq\fix(\dot u)$.
		
		\smallskip\noindent\emph{(1) Pairing.}
		Let $S:=S_x\cup S_y$ (still $<\kappa$). Any $g\in\fix(S)$ fixes $\dot x$ and $\dot y$, hence fixes the canonical pair
		$\langle\dot x,\dot y\rangle$. Thus $S$ supports the pair and $\langle\dot x,\dot y\rangle\in\HS$.
		
		\smallskip\noindent\emph{(2) $<\kappa$--ary unions.}
		Let $S:=S_A\cup\bigcup_{i\in I}S_i$ where each $S_i$ supports $\dot x_i$ (all $<\kappa$). Regularity of $\kappa$
		implies $|S|<\kappa$. Any $g\in\fix(S)$ fixes $\dot A$ and each constituent $\dot x_i$, hence fixes the union name
		$\bigcup\dot A$. So $\bigcup\dot A\in\HS$.
		
		\smallskip\noindent\emph{(3) Separation.}
		Let $S:=S_A\cup S_{\bar p}$. If $g\in\fix(S)$ then $g$ fixes $\dot A$ and $\bar{\dot p}$, and by automorphism
		invariance of the forcing relation (here trivial forcing) we have
		\[
		1\Vdash \varphi(\sigma,\bar{\dot p})\iff 1\Vdash \varphi(g\cdot\sigma,\bar{\dot p}).
		\]
		Thus $g$ preserves membership in $\dot B$, so $g\cdot\dot B=\dot B$. Hence $S$ supports $\dot B$ and $\dot B\in\HS$.
		
		\smallskip\noindent\emph{(4) Replacement (range).}
		Let $S:=S_f\cup S_A$. The range name $\dot{\mathrm{rng}}$ is definable from $\dot f$ and $\dot A$ by a fixed
		(first-order) pattern. If $g\in\fix(S)$ then $g$ fixes $\dot f$ and $\dot A$, hence fixes the resulting
		$\dot{\mathrm{rng}}$. Therefore $S$ supports $\dot{\mathrm{rng}}$, so $\dot{\mathrm{rng}}\in\HS$, and evaluating by the
		generic (trivial here) yields Replacement in the symmetric model.
		
		\smallskip
		This establishes (1)--(4). \qedhere
	\end{proof}
	
	\begin{proposition}[$\kappa$-completeness via supports]
		\label{prop:FM-kappa-complete}
		Let $\mathcal{F}_{\mathrm{perm}}$ be the filter on $\mathcal G=\Sym(A)$ generated by the subgroups
		\[
		\fix(S)\ :=\ \{\,g\in\mathcal G:\ g{\upharpoonright}S=\mathrm{id}\,\}\qquad(S\subseteq A,\ |S|<\kappa),
		\]
		where $\kappa$ is regular. Then $\mathcal{F}_{\mathrm{perm}}$ is normal and $\kappa$-complete. In particular, if
		$\{\dot x_i:i\in I\}$ are HS names with $|I|<\kappa$, then $\bigcap_{i\in I}\fix(\dot x_i)\in\mathcal{F}_{\mathrm{perm}}$.
	\end{proposition}
	
	\begin{proof}
		\emph{Normality.} For any $S\subseteq A$ and $g\in\mathcal G$,
		\[
		g\,\fix(S)\,g^{-1}\ =\ \fix\bigl(g[S]\bigr),
		\]
		so the generating family is closed under conjugation; hence the generated filter is normal.
		
		\emph{$\kappa$--completeness.} Let $\{H_i:i\in I\}\subseteq\mathcal{F}_{\mathrm{perm}}$ with $|I|<\kappa$.
		For each $i$ choose $S_i\subseteq A$ with $|S_i|<\kappa$ such that $\fix(S_i)\subseteq H_i$
		(possible since $\mathcal F_{\mathrm{perm}}$ is the upward closure of finite intersections of generators and
		$\fix(S)\cap\fix(T)=\fix(S\cup T)$). Let $S:=\bigcup_{i\in I}S_i$. Regularity of $\kappa$ gives $|S|<\kappa$ and
		\[
		\fix(S)\ =\ \bigcap_{i\in I}\fix(S_i)\ \subseteq\ \bigcap_{i\in I}H_i.
		\]
		Hence $\bigcap_{i\in I}H_i\in\mathcal{F}_{\mathrm{perm}}$ by upward closure.
		
		\smallskip
		For HS names $\dot x_i$, pick supports $S_i$ of size $<\kappa$ with $\fix(S_i)\subseteq\fix(\dot x_i)$
		(Lemma~\ref{lem:FM-HS}) and apply the preceding paragraph. \qedhere
	\end{proof}
	
	\begin{corollary}[Choice fragments in permutation data]
		\label{cor:FM-DC}
		Assume $V\models\ZF+\DC$ and work with permutation data whose filter $\mathcal{F}_{\mathrm{perm}}$ is generated by
		$\fix(S)$ for $|S|<\kappa$ with $\kappa$ regular. Then the corresponding symmetric submodel $M$ satisfies $\ZF$ and
		$\DC_{<\kappa}$.
	\end{corollary}
	
	\begin{proof}
		By Proposition~\ref{prop:FM-kappa-complete}, $\mathcal{F}_{\mathrm{perm}}$ is normal and $\kappa$--complete.
		Applying the general preservation argument to this (trivial-forcing) symmetric system yields $M\models\ZF$ and
		$M\models\DC_{<\kappa}$ (cf.\ Theorems~\ref{thm:kappa-zf-preserved} and~\ref{thm:DCless-preservation}).
		
		\smallskip
		\emph{Note.} When $\kappa=\omega$ (classical FM), “$\kappa$--complete” means closure under \emph{finite} intersections,
		so the conclusion is only $\DC_{<\omega}$ (finite dependent choice), which is already provable in $\ZF$; full $\DC$
		typically requires the countable-support/$\omega_1$--complete setup from Section~\ref{sec:HS-countable}. \qedhere
	\end{proof}
	
	\begin{remark}[Coordinate vs.\ name stabilizers]
		\label{rem:FM-coordinate-name}
		In the FM setting, basic \emph{coordinate stabilizers} $\fix(S)$ (for $S\subseteq A$ small) generate the filter.
		If $\dot x$ has support $S$, then $\fix(S)\subseteq\fix(\dot x)$, so name stabilizers are (typically) larger than
		coordinate stabilizers but are equivalent \emph{modulo} the filter: $\fix(\dot x)\in\mathcal F_{\mathrm{perm}}$
		iff $\dot x$ has some small coordinate support (Lemma~\ref{lem:FM-HS}). Our general framework prefers name stabilizers,
		which extends smoothly to iterations where “coordinates” are ordinal stages rather than atoms
		(Sections~\ref{sec:HS-countable} and~\ref{sec:kappa-support}).
	\end{remark}
	
	\paragraph{Takeaway.}
	Classical permutation models are the $<\!\kappa$-support case with no nontrivial forcing: filter completeness is exactly
	closure under intersections of stabilizers of small supports. Passing to iterated symmetric extensions replaces atoms by
	stages and imports the same combinatorics of supports and $\kappa$-completeness to control which choice fragments
	(e.g.\ $\DC_{<\kappa}$) survive. See Table~\ref{tab:support-comparison} for a systematic comparison across support regimes.
	
	\subsection{Support Regimes: A Comparison}
	\label{subsec:support-comparison}
	
	\begin{table}[H]
		\centering
		\caption{Support regimes and provable consequences in the symmetric submodel}
		\label{tab:support-comparison}
		\begin{tabularx}{\linewidth}{L{2.5cm} C{1.9cm} C{2.2cm} L{2.6cm} L{2.2cm} L{2.0cm}}
			\toprule
			\textbf{Support} & \textbf{Length treated} & \textbf{Filter compl.} & \textbf{Choice preserved (provable)} & \textbf{Cardinal assumps.} & \textbf{Example} \\
			\midrule
			Finite & $\omega_1$ & finite $\cap$ & $\mathrm{DC}_{<\omega}$ (trivial; $\mathrm{DC}$ may fail) & minimal & \S\ref{subsec:perm-models} \\
			Countable & $\omega_2$ & $\omega_1$-complete & $\mathrm{DC}$ ($=\mathrm{DC}_\omega$) & $\omega_1$ regular & \S\ref{sec:cohen-example} \\
			$<\omega_2$ & $\omega_3$ & $\omega_2$-complete & $\mathrm{DC}_{<\omega_2}$ & $\omega_2$ regular & \S\ref{sec:more-examples} \\
			$<\kappa$ & $\kappa^+$ & $\kappa$-complete & $\mathrm{DC}_{<\kappa}$ & $\kappa$ regular & \S\ref{sec:kappa-support} \\
			\bottomrule
		\end{tabularx}
		\footnotesize\emph{Note.} “Length treated” refers to the regular-cardinal regime developed in this paper (the pattern $\kappa\mapsto\kappa^+$); longer iterations across singular cofinalities require the singular-case framework of \S\ref{sec:singular-overview}--\S\ref{sec:hypotheses}. The upgrade to $\mathrm{DC}_\kappa$ may hold under the Localization Hypothesis (Def.~\ref{def:localization-hypothesis}), cf.\ Thm.~\ref{thm:DCkappa-localization}.
	\end{table}
	
	\section{Singular $\kappa$: Overview and Roadmap}
	\label{sec:singular-overview}
	
	From this point on, $\kappa$ denotes a \emph{singular strong limit} cardinal
	(earlier sections used $\kappa$ as regular). We keep the same letter to avoid
	proliferating symbols; each statement specifies the relevant hypotheses
	(regular vs.\ singular) as needed.
	
	\paragraph{Why new bookkeeping is needed.}
	At singular $\kappa$, the regular-case simplifications fail at cofinality
	$\cf(\kappa)$: unions of $\cf(\kappa)$ many $<\kappa$--sized supports can reach
	size $\kappa$, and the stage-bounding argument at limits of cofinality
	$\cf(\kappa)$ can break. Accordingly we replace the naive “$<\kappa$--support at
	limits” compression by \emph{trees of conditions} and \emph{block/partition
		stabilizers} coordinated along a fixed cofinal sequence in $\kappa$.
	
	\paragraph{Tracks.}
	The development proceeds in two tracks:
	\begin{itemize}
		\item \textbf{Track A ($\cf(\kappa)=\omega$).}
		Full proofs using revised countable support (RCS)-style organization and
		$\Sigma$--Prikry bookkeeping, together with trees of conditions; we prove the
		main preservation theorems and give a worked $\aleph_\omega$-length example
		(\S\ref{sec:singular-cf-omega}).
		\item \textbf{Track B (general singular $\kappa$).}
		A framework and key lemmas for arbitrary $\cf(\kappa)<\kappa$, emphasizing
		limit stages $\lambda$ with $\cf(\lambda)=\cf(\kappa)$
		(\S\ref{sec:limit-cf-equal} and \S\ref{sec:hypotheses}).
	\end{itemize}
	
	\paragraph{Targets.}
	Across both tracks we aim to:
	\begin{enumerate}[label=(\roman*), leftmargin=2.2em]
		\item avoid collapsing $\kappa$;
		\item add no new bounded subsets of $\kappa$;
		\item preserve $\mathrm{DC}_{<\kappa}$ under explicit hypotheses; and
		\item integrate the relevant Shelah-style completeness/closure conditions with
		symmetric filters via trees of conditions.
	\end{enumerate}
	
	\section{Foundations at Singular $\kappa$}
	\label{sec:singular-foundations}
	\subsection{Shelah completeness and the notation bridge}
	We work with Shelah's game-theoretic completeness at each stage. For readability we use $(E_0,E_1)$-completeness as a mnemonic, and record the notational correspondence:
	\begin{remark}[Notation bridge]\label{rem:notation-bridge}
		Throughout the singular-$\kappa$ part we use Shelah's game-theoretic completeness hypotheses.
		Our mnemonic $(E_0,E_1)$ corresponds to Shelah's $(\hat\varepsilon_0,\hat\varepsilon_1)$, and also to the
		$(S_0,S_1,D)$ parametrization in the companion formulation used in the literature. Whenever we cite a
		preservation result, we explicitly indicate which variant is being applied.
	\end{remark}
	
	\subsection{Symmetric systems at $<\kappa$}
	We work in the symmetric-extensions iteration framework developed in earlier sections: at each stage $\alpha$
	we have a system $(\mathbb P_\alpha,\mathcal G_\alpha,\mathcal F_\alpha)$ consisting of a forcing, its
	automorphism group, and a normal filter of subgroups. Hereditarily symmetric names determine the intermediate
	model, and supports are taken to be of size $<\kappa$.
	
	In the singular setting we modify only the \emph{completeness} and \emph{limit} mechanics. Concretely:
	\begin{enumerate}[label=(\roman*), leftmargin=2.2em]
		\item At each stage, $\mathcal F_\alpha$ is assumed \emph{normal} and satisfies the relevant
		\emph{Shelah/game completeness} used in the preservation arguments (recorded below via the
		$(E_0,E_1)$ mnemonic), rather than relying on “$<\kappa$-completeness by union of supports.”
		\item We maintain tenacity/excellence of supports exactly as in the regular framework.
		\item At limit stages of cofinality $\cf(\kappa)$ (and in particular where union-of-support compression can fail),
		we use \emph{trees of conditions} and \emph{block/partition stabilizers} to amalgamate automorphisms and to
		witness the closure properties of the limit filter.
	\end{enumerate}
	
	\subsection{DC via the $\kappa$-Baire perspective}
	Throughout §§\ref{sec:singular-cf-omega}–\ref{sec:limit-cf-equal} we adopt the following standing ground-model and Baire assumptions.
	
	\begin{hypothesis}[Ground model and Baire framework]\label{hyp:zf-baire}
		$V\models \ZF+DC_{<\kappa}$.
		Each iterand satisfies the $<\kappa$--Baire (game/strategy) hypothesis used in the preservation templates
		(e.g.\ $<\kappa$--strategic closure, or the corresponding $(E_0,E_1)$-game condition as appropriate).
		At each stage, the symmetry filters are normal and satisfy the relevant Shelah/game completeness needed to
		reflect $<\kappa$-length choice constructions into the hereditarily symmetric model; at the critical limits,
		this completeness is witnessed via the tree/partition mechanics rather than via union-of-support compression.
	\end{hypothesis}
	
	\begin{remark}
		This replaces any informal “closed under $<\kappa$-sequences in a ZFC extension” phrasing. We adopt the \emph{$\kappa$-Baire} formulation uniformly and cite
		Karagila~\cite{Karagila2019DC}, Thm.~2.1 (cf.\ \cite{Karagila2019DC}, Lemmas~3.2--3.3 for the preservation templates).
		All $DC_{<\kappa}$ preservation claims in §§\ref{sec:singular-cf-omega}–\ref{sec:limit-cf-equal} explicitly invoke Hyp.~\ref{hyp:zf-baire}.
	\end{remark}
	
	\begin{remark}[Sufficiency vs.\ necessity]
		Hypothesis~\ref{hyp:zf-baire} is a standing \emph{sufficient} schema for all $DC_{<\kappa}$ claims in these sections.
		The $<\kappa$-completeness of the limit symmetry filter is essential for reflecting sequences to the HS model, whereas
		the $<\kappa$-Baire condition on iterands can be replaced by other fusion/strategy hypotheses on a case-by-case basis.
		The ground assumption $V\models\ZF+DC_{<\kappa}$ is standard for preservation; it can be weakened when one aims to
		\emph{create} $DC_{<\kappa}$ under stronger structural hypotheses (e.g., Localization).
	\end{remark}
	
	\begin{remark}[ZF-only reading (ignoring $DC$, $DC_\kappa$, $DC_{<\kappa}$)]
		If one is interested only in preserving $\ZF$ (and not $DC$ fragments), then the singular-$\kappa$ development
		does not require Hypothesis~\ref{hyp:zf-baire}'s \emph{$<\kappa$-Baire/strategy} clause or Localization.
		What \emph{is} still required in the singular sections is the structural limit bookkeeping (trees of conditions,
		block/partition stabilizers, and the corresponding normality/completeness hypotheses for the limit symmetry filter).
		
		By contrast, in the \emph{regular} parts of the paper one may ignore the $DC$ hypotheses entirely and rely on:
		\begin{itemize}
			\item normality and the relevant completeness of the limit symmetry filter
			(\,$\omega_1$--completeness for countable support; $\kappa$--completeness for $<\kappa$ support\,); and
			\item stage-bounding (Lemma~\ref{lemma:stage-bounding} / Lemma~\ref{lem:kappa-stage-bounding})
			together with the standard successor/limit transfer lemmas.
		\end{itemize}
		Under those assumptions the ZF preservation theorems
		(Theorems~\ref{thm:zf-preserved}, \ref{thm:kappa-zf-preserved}, \ref{thm:class-length})
		hold as stated. Sections or theorems whose main claims are $DC$/$DC_{<\kappa}$/$DC_\kappa$ can then be
		treated as optional add-ons.
	\end{remark}
	
	\section{The case $\cf(\kappa)=\omega$ (Track A)}
	\label{sec:singular-cf-omega}
	
	\subsection{Iteration template: RCS + symmetry}
	
	Fix an increasing sequence of regulars $\langle \kappa_n:n<\omega\rangle$ cofinal in $\kappa$.
	We partition the iteration coordinates into blocks $\langle B_n:n<\omega\rangle$ (or equivalently fix
	a block map $b:\mathrm{Ord}\to\omega$) so that bookkeeping is aligned with $\langle\kappa_n\rangle$.
	
	\begin{convention}[Track A support policy]\label{conv:trackA-support}
		Conditions and HS names are organized by \emph{trees of conditions} (RCS presentation) so that:
		\begin{enumerate}[label=(\roman*), leftmargin=2.2em]
			\item every condition has \emph{countable} support (RCS), hence in particular $<\kappa$ support;
			\item every HS name $\dot x$ is stabilized by a subgroup fixing all but finitely many blocks, so $\dot x$ is
			\emph{block-bounded}: there is $n(\dot x)<\omega$ such that $\supp(\dot x)\subseteq \bigcup_{m<n(\dot x)} B_m$.
		\end{enumerate}
	\end{convention}
	
	\subsection{Stage-bounding at singular base}
	
	\subsection{Stage-bounding at singular base}
	\begin{lemma}[Stage-bounding (Track A form)]\label{lem:sb-sing}
		Let $\alpha$ be a limit stage with $\cf(\alpha)^{V[G_\alpha]}>\omega$.
		Suppose the iteration is presented by trees of conditions as in Convention~\ref{conv:trackA-support}
		(so relevant names/conditions have countable support).
		If $S\subseteq V[G_\alpha]$ and $|S|<\cf(\alpha)^{V[G_\alpha]}$, then there is $\eta<\alpha$ such that
		$S\subseteq V[G_\eta]$ (hence $S\in V[G_\eta]$).
	\end{lemma}
	
	\begin{proof}
		Let $\theta=\cf(\alpha)^{V[G_\alpha]}$; by hypothesis $\theta$ is uncountable regular.
		For each $x\in S$ fix a $\mathbb P_\alpha$-name $\dot x$ with $\dot x^{G_\alpha}=x$.
		By the Track~A support policy, each $\supp(\dot x)$ is countable.
		Let $U=\bigcup_{x\in S}\supp(\dot x)\subseteq\alpha$.
		Since $|S|<\theta$ and $\theta$ is regular, $|U|<\theta$ (a union of $<\theta$ many countable sets).
		Therefore $U$ is bounded in $\alpha$, so choose $\eta<\alpha$ with $U\subseteq\eta$.
		Then every $\dot x$ is a $\mathbb P_\eta$-name, hence $x\in V[G_\eta]$, as required.
	\end{proof}
	
	\begin{remark}
		If $\cf(\alpha)=\omega$, stage-bounding can fail (countable sets of stages may be cofinal).
		This is exactly the case handled by the inverse-limit / trees-of-conditions presentation.
	\end{remark}
	
	\subsection{Successor-stage filter transfer}
	
	\begin{lemma}[Successor transfer (Track A)]\label{lem:succ-transfer}
		If $\mathcal F_\alpha$ is $<\kappa$-complete and $\Vdash_{\mathbb P_\alpha}$ “$\dot{\mathbb Q}_\alpha$ is $(E_0,E_1)$-complete,”
		and the successor-stage generators added at $\alpha+1$ have \emph{countable} supports (RCS),
		then $\mathcal F_{\alpha+1}$ is $<\kappa$-complete and $(E_0,E_1)$-completeness is preserved at $\alpha+1$.
	\end{lemma}
	
	\begin{proof}
		Work in $V^{\mathbb P_\alpha}$. Let $\dot{\mathbb Q}_\alpha$ be $(E_0,E_1)$-complete and let $\mathcal F_\alpha$ be $<\kappa$-complete in $V$. The symmetry filter at $\alpha+1$ is generated by the images of $\mathcal F_\alpha$ along the canonical embedding and the pointwise fixators of the $\mathcal G_\alpha$-names used to define $\dot{\mathbb Q}_\alpha$.
		\smallskip
		
		\noindent\emph{$<\kappa$-completeness.} Take $\gamma<\kappa$ and a family $\{H_\xi:\xi<\gamma\}\subseteq\mathcal F_{\alpha+1}$. Each $H_\xi$ contains a subgroup of the form $K_\xi\cap L_\xi$ with $K_\xi\in\mathcal F_\alpha$ and $L_\xi$ a pointwise fixator of a set of coordinates/names of size $<\kappa$. Since $\mathcal F_\alpha$ is $<\kappa$-complete, $\bigcap_{\xi<\gamma}K_\xi\in\mathcal F_\alpha$, and $\bigcap_{\xi<\gamma}L_\xi$ is the pointwise fixator of a union of $<\kappa$ many
		\emph{countable} support sets (Track~A/RCS), hence again a permitted generator since
		$\cf(\kappa)=\omega$ implies that a union of $<\kappa$ many countable sets has size $<\kappa$.
		
		\noindent\emph{Completeness transfer.}
		Consider a play of the $(E_0,E_1)$-game on $\mathbb P_\alpha*\dot{\mathbb Q}_\alpha$.
		\begin{itemize}
			\item \textit{Project to $\mathbb P_\alpha$.} COM follows the inductive winning strategy on $\mathbb P_\alpha$ to determine the $\mathbb P_\alpha$-part $p_\alpha$ at each move.
			\item \textit{Lift to $\dot{\mathbb Q}_\alpha$.} Working in $V^{\mathbb P_\alpha}$, the name $\dot{\mathbb Q}_\alpha$ is $(E_0,E_1)$-complete, so COM has a name for a winning reply $\,\dot q\,$ on the $\dot{\mathbb Q}_\alpha$-board.
			\item \textit{Automorphism compatibility.} Since $\mathcal F_{\alpha+1}$ contains pointwise fixators for the names used to define $\dot{\mathbb Q}_\alpha$, the combined play is respected by the symmetry (names do not move).
		\end{itemize}
		Composing the strategies yields a COM-winning strategy on $\mathbb P_{\alpha+1}$, hence $(E_0,E_1)$-completeness at $\alpha+1$.
	\end{proof}
	
	\subsection{Limit filters when $\cf(\lambda)\neq\cf(\kappa)$}
	
	\begin{lemma}[Limit transfer for unequal cofinalities]\label{lem:limit-uneq}
		Let $\lambda$ be a limit stage with $\cf(\lambda)\ge \kappa$.
		Assume each $\mathcal F_\alpha$ ($\alpha<\lambda$) is $<\kappa$-complete.
		Then the direct-limit filter $\mathcal F_\lambda$ is $<\kappa$-complete.
	\end{lemma}
	
	\begin{proof}
		Let $\mu<\kappa$ and let $\{H_i:i<\mu\}\subseteq\mathcal F_\lambda$.
		For each $i$ choose $\alpha_i<\lambda$ and $H_i'\in\mathcal F_{\alpha_i}$ such that
		$H_i=(\rho^{\alpha_i}_\lambda)^{-1}[H_i']$, where $\rho^{\alpha_i}_\lambda:\mathcal G_\lambda\to\mathcal G_{\alpha_i}$
		is restriction to stage $\alpha_i$.
		Since $\cf(\lambda)\ge\kappa>\mu$, the set $\{\alpha_i:i<\mu\}$ is bounded in $\lambda$, say by $\beta<\lambda$.
		Using $<\kappa$-completeness in $\mathcal F_\beta$ we get $\bigcap_{i<\mu}H_i'\in\mathcal F_\beta$, and pulling back yields
		$\bigcap_{i<\mu}H_i\in\mathcal F_\lambda$.
	\end{proof}
	
	\begin{remark}
		When $\cf(\lambda)<\kappa$, the direct-limit argument above need not apply. In this case we use the inverse-limit presentation of the iteration at small-cofinality limits (cf.\ Section~2 and the tuple/partition stabilizer machinery of Section~3), which preserves $<\kappa$-completeness by a coordinatewise argument.
	\end{remark}
	
	\begin{lemma}[No new bounded subsets from high blocks]
		\label{lem:trackA-tail-no-bounded-subsets}
		Fix $n<\omega$ and $\lambda<\kappa_n$.  Work in Track~A with $\kappa=\aleph_\omega$ and a fixed
		increasing cofinal sequence $\langle\kappa_m:m<\omega\rangle$ of regular cardinals.
		
		Assume the iteration has length $\kappa$ and is organized into blocks
		$B_m=[\kappa_m,\kappa_{m+1})$, and that every condition has \emph{block-finite support}:
		for every $p\in\mathbb P_\kappa$ and every $m<\omega$,
		\[
		|\supp(p)\cap B_m|<\omega.
		\]
		Assume moreover that for each $\alpha\in B_m$,
		\[
		\Vdash_{\mathbb P_\alpha}\ \dot{\mathbb Q}_\alpha=\Add(\kappa_m,1),
		\]
		viewed as a constant $\mathbb P_\alpha$--name.
		
		Then the tail forcing from stage $\kappa_n$ onward is $<\kappa_n$-closed. In particular it adds no
		new subsets of any $\lambda<\kappa_n$. Consequently, in the two-step factorization
		$\mathbb P_\kappa\simeq \mathbb P_{\kappa_n}*\dot{\mathbb R}$, if $p\Vdash \dot A\subseteq\check\lambda$ then
		\[
		p\Vdash\ \dot A^{\,\dot G}\in V[\dot G{\upharpoonright}\kappa_n],
		\]
		i.e.\ $\dot A$ is forced to be a member of the intermediate extension by $\mathbb P_{\kappa_n}$.
	\end{lemma}
	
	\begin{proof}
		Let $\mathbb P_{\kappa_n,\kappa}$ denote the usual tail forcing:
		conditions are tails $p{\upharpoonright}[\kappa_n,\kappa)$ ordered coordinatewise
		(so $q\le p$ iff $q(\alpha)\le_{\dot{\mathbb Q}_\alpha}p(\alpha)$ for all $\alpha\ge\kappa_n$).
		
		\smallskip\noindent\textbf{(1) $<\kappa_n$-closure of the tail.}
		Fix $\mu<\kappa_n$ and a decreasing sequence $\langle p_i:i<\mu\rangle$ in $\mathbb P_{\kappa_n,\kappa}$.
		Let
		\[
		U:=\bigcup_{i<\mu}\supp(p_i)\ \subseteq\ [\kappa_n,\kappa).
		\]
		For each $m\ge n$, the set $U\cap B_m$ is a union of $\mu$ many finite sets, hence
		$|U\cap B_m|<\kappa_n$ (using that $\kappa_n$ is regular and $\mu<\kappa_n$). Therefore
		\[
		|U|\ \le\ \Bigl|\bigcup_{m\ge n}(U\cap B_m)\Bigr|\ \le\ \omega\cdot \kappa_n\ =\ \kappa_n\ <\ \kappa,
		\]
		so $U$ is an allowed support set for a tail condition.
		
		Fix $\alpha\in U$ and write $m=b(\alpha)$ (so $\alpha\in B_m$ and $m\ge n$). The sequence
		$\langle p_i(\alpha):i<\mu\rangle$ is decreasing in $\Add(\kappa_m,1)$.
		Since $\Add(\kappa_m,1)$ is $<\kappa_m$-closed (hence $<\kappa_n$-closed because $\kappa_n\le\kappa_m$),
		there is $r_\alpha\in\Add(\kappa_m,1)$ with $r_\alpha\le p_i(\alpha)$ for all $i<\mu$.
		
		Define $q\in\mathbb P_{\kappa_n,\kappa}$ by $\supp(q)=U$ and $q(\alpha)=r_\alpha$ for $\alpha\in U$
		(and the trivial condition elsewhere). Because each iterand is the constant name
		$\check{\Add(\kappa_m,1)}$, this coordinatewise definition produces a legitimate tail condition,
		and clearly $q\le p_i$ for all $i<\mu$. Hence $\mathbb P_{\kappa_n,\kappa}$ is $<\kappa_n$-closed.
		
		\smallskip\noindent\textbf{(2) No new bounded subsets.}
		Fix $\lambda<\kappa_n$. Let $\dot X$ be a $\mathbb P_{\kappa_n,\kappa}$-name for a subset of $\check\lambda$
		and let $p\in\mathbb P_{\kappa_n,\kappa}$.
		For each $\xi<\lambda$ let
		\[
		D_\xi:=\{\,q\le p : q\text{ decides }\xi\in\dot X\,\}.
		\]
		Each $D_\xi$ is dense below $p$. Using dependent choice of length $\lambda$
		(which follows from $DC_{<\kappa}$ in the standing hypotheses for Track~A),
		build a decreasing sequence $\langle p_\xi:\xi<\lambda\rangle$ with $p_\xi\in D_\xi$.
		By (1) and $\lambda<\kappa_n$, there is a lower bound $q\le p_\xi$ for all $\xi<\lambda$.
		Then $q$ decides every statement “$\xi\in\dot X$” for $\xi<\lambda$, hence decides $\dot X$ outright.
		Therefore $\mathbb P_{\kappa_n,\kappa}$ adds no new subsets of $\lambda$.
		
		\smallskip\noindent\textbf{(3) Consequence for two-step factorization.}
		Factor $\mathbb P_\kappa\simeq \mathbb P_{\kappa_n}*\dot{\mathbb R}$, where
		$\dot{\mathbb R}$ names the tail $\mathbb P_{\kappa_n,\kappa}$.
		By (1) the tail is forced to be $<\kappa_n$-closed, hence by (2) it is forced to add no new subsets of
		any $\lambda<\kappa_n$. Therefore, if $p\Vdash \dot A\subseteq\check\lambda$ with $\lambda<\kappa_n$,
		then in any $V[G_{\kappa_n}][H]$ containing $p$ we have $\dot A^{G_{\kappa_n}*H}\in V[G_{\kappa_n}]$,
		as claimed.
	\end{proof}
	
	\subsection{Worked example: an $\aleph_\omega$-length Cohen-style iteration}
	\label{subsec:aleph-omega-example}
	Fix $\kappa=\aleph_\omega$ and a cofinal sequence $\langle \kappa_n:n<\omega\rangle$ of regulars with $\sup_n\kappa_n=\kappa$.
	Define an iteration $\langle \mathbb P_\alpha,\dot{\mathbb Q}_\alpha,\mathcal F_\alpha:\alpha<\kappa\rangle$
	(of length $\kappa$) using the Track~A support policy (revised countable support via trees of conditions, hence in particular $<\kappa$-support),
	and partition $\kappa$ into blocks $B_n=[\kappa_n,\kappa_{n+1})$.
	For $\alpha\in B_n$ let
	\[
	\Vdash_{\mathbb P_\alpha}\ \dot{\mathbb Q}_\alpha=\Add(\kappa_n,1).
	\]
	the Cohen algebra adding one subset to $\kappa_n$ (or finitely many—any fixed finite number works). Let $\mathcal G_\alpha$ be the full automorphism group of the iterand and let $\mathcal F_\alpha$ be the normal filter generated by pointwise fixators of sets of size $<\kappa$ together with the block-partition stabilizer (Def.~\ref{def:block-stab}). Then:
	\begin{itemize}
		\item \emph{Stage-bounding for bounded objects in the HS model:}
		If $\dot A\in\HS$ and $p\Vdash \dot A\subseteq\check\lambda$ for some $\lambda<\kappa$,
		choose $n$ with $\lambda<\kappa_n$.
		By Lemma~\ref{lem:trackA-tail-no-bounded-subsets}, the tail beyond $\kappa_n$ adds no new subsets of $\lambda$,
		so $\dot A$ is equivalent to a $\mathbb P_{\kappa_n}$-name. Hence bounded subsets of $\kappa$ in the final HS model
		already occur in a bounded initial segment.
		\item \emph{$DC_{<\kappa}$:} By Hyp.~\ref{hyp:zf-baire} and the fact that $\Add(\kappa_n,1)$ is $\kappa_n$-closed under the relevant game at the level of the $\kappa$-Baire formulation, together with $<\kappa$-completeness of $\mathcal F_\alpha$, Karagila's criterion \cite{Karagila2019DC} yields preservation of $DC_{<\kappa}$.
		\item \emph{Failure of AC:} Standard arguments for symmetric Cohen models apply: the filter of pointwise fixators eliminates wellorderings of the added families across blocks, preventing a choice function for a definable family of nonempty sets (details follow the classical Cohen symmetric model, adapted to blocks).
	\end{itemize}
	The comparison with the $\omega_1$-length regular case is summarized below.
	
	\paragraph{Comparison with the $\omega_1$-length regular case.}
	Key differences: (a) base $\kappa$ is singular rather than regular; (b) supports are organized via a fixed $\omega$-sequence cofinal in $\kappa$ and use block partitions; (c) filters must be $<\kappa$-complete (rather than $\kappa$-complete); (d) stage-bounding relies on $\omega$-strategic closure and chain condition arguments instead of $\kappa$-closure; (e) trees of conditions are used essentially at selected limits even in Track~A.
	
	\subsection{$DC_{<\kappa}$ preservation (Track A)}
	
	\begin{theorem}[Preservation of $DC_{<\kappa}$]\label{thm:dc-trackA}
		Assume Hyp.~\ref{hyp:zf-baire} and the Track~A support policy (Convention~\ref{conv:trackA-support}).
		Then $DC_{<\kappa}$ holds in the final symmetric extension.
	\end{theorem}
	
	\begin{proof}
		Fix $\mu<\kappa$. Let $M:=V(\mathbb P_\lambda)_{\HS}$ be the final symmetric model.
		Let $\dot A,\dot R\in\HS$ and assume $p\in\mathbb P_\lambda$ forces that $\dot R$ is serial on $\dot A$.
		
		\smallskip\noindent\emph{(Forcing-side construction.)}
		By Hyp.~\ref{hyp:zf-baire}, the iterands satisfy the required game/strategy hypothesis and the iteration
		is organized by trees of conditions at limits of cofinality $\omega$. Using the standard preservation template
		(e.g.\ the $(E_0,E_1)$-completeness formulation cited in §\ref{sec:singular-foundations}), and the ground
		assumption $V\models DC_\mu$, we build in $V$ a descending sequence of conditions
		$\langle p_\xi:\xi<\mu\rangle$ and names $\langle \dot a_\xi:\xi<\mu\rangle$ such that
		\[
		p_\xi\Vdash \dot a_\xi\in\dot A
		\quad\text{and}\quad
		p_{\xi+1}\Vdash \langle \dot a_\xi,\dot a_{\xi+1}\rangle\in\dot R
		\qquad(\xi+1<\mu).
		\]
		(At limits $<\mu$ we use the available fusion/strategy hypothesis on the tree presentation.)
		
		\smallskip\noindent\emph{(Symmetry step.)}
		By Convention~\ref{conv:trackA-support}, each $\dot a_\xi$ and each $p_\xi$ has \emph{countable} support, and each
		$\dot a_\xi\in\HS$ has a symmetry witness containing a pointwise fixer of that countable support.
		Let $S$ be the union of these supports over $\xi<\mu$; since $\mu<\kappa$ and $\cf(\kappa)=\omega$, $|S|<\kappa$.
		Hence the pointwise stabilizer $\fix(S)$ is a permitted subgroup and lies in the limit symmetry filter.
		Moreover, $\fix(S)$ fixes every $\dot a_\xi$ and every $p_\xi$, so it fixes the canonical name
		\[
		\dot s:=\{\langle\langle\check\xi,\dot a_\xi\rangle,p_\xi\rangle:\xi<\mu\}.
		\]
		Thus $\dot s\in\HS$, so $s:=\dot s^G\in M$.
		
		\smallskip\noindent\emph{(Conclusion.)}
		In $M$ the function $s:\mu\to A$ witnesses $DC_\mu$ for $(A,R)=(\dot A^G,\dot R^G)$.
		Since $\mu<\kappa$ was arbitrary, $M\models DC_{<\kappa}$.
	\end{proof}
	
	\begin{remark}
		Compare with Karagila~\cite[Thm.~2.1]{Karagila2019DC}: the same Baire/strategic-closure hypothesis preserves $DC$
		along iterations. Our proof packages the forcing-side fusion (Hyp.~\ref{hyp:zf-baire}; trees of conditions at
		$\cf=\omega$ limits) with the symmetry step via the $\kappa$–completeness of the stabilizer filter
		(Theorem~\ref{thm:kappa-filter-generation}).
	\end{remark}
	
	\section{Block-partition stabilizers and excellent supports}
	\label{sec:block-partition}
	
	\begin{definition}[Block decomposition, block-bounded supports, and block-partition stabilizers]\label{def:block-stab}
		(Track A: $\cf(\kappa)=\omega$.) Fix $\kappa=\sup_{n<\omega}\kappa_n$ where
		$\langle \kappa_n:n<\omega\rangle$ is strictly increasing and each $\kappa_n$ is regular.
		
		Fix an iteration length $\gamma$ and a \emph{block map} $b:\gamma\to\omega$.
		Write $B_n:=b^{-1}(\{n\})$ for the $n$th block.
		
		\begin{enumerate}
			\item For $N<\omega$, a set $u\subseteq\gamma$ is \emph{$N$--block-bounded} (w.r.t.\ $b$) if
			\[
			\forall n<\omega\qquad |u\cap B_n|<\kappa_N.
			\]
			We say that $u$ is \emph{block-bounded} if it is $N$--block-bounded for some $N<\omega$.
			
			\item Let $\mathcal G_\gamma^{\,b}$ be the group of coordinate permutations $\pi$ (and their induced
			automorphisms of $\mathbb P_\gamma$) such that:
			\[
			(\forall \alpha<\gamma)\ b(\pi(\alpha))=b(\alpha)
			\quad\text{and}\quad
			(\exists F\in[\omega]^{<\omega})\ (\forall n\notin F)\ (\forall \alpha\in B_n)\ \pi(\alpha)=\alpha.
			\]
			(So $\pi$ preserves blocks everywhere and is the identity on all but finitely many blocks.)
			
			\item For block-bounded $u\subseteq\gamma$, define the \emph{block-partition stabilizer} of $u$ to be
			\[
			\fix_b(u)\ :=\ \{\pi\in\mathcal G_\gamma^{\,b}:\ \pi{\upharpoonright}u=\mathrm{id}_u\}.
			\]
		\end{enumerate}
	\end{definition}
	
	\begin{proposition}[Support properties and stratified completeness]\label{prop:pi}
		Assume Track A: $\cf(\kappa)=\omega$ and fix $\kappa=\sup_{n<\omega}\kappa_n$ as in
		Definition~\ref{def:block-stab}. Let $\gamma$ be an iteration length with block map $b:\gamma\to\omega$,
		and let $\mathcal G_\gamma^{\,b}$ act as the coordinate-permutation automorphism group on $\mathbb P_\gamma$.
		
		For each $N<\omega$, let $\tilde{\mathcal F}^{\,N}_\gamma$ be the normal filter of subgroups of
		$\mathcal G_\gamma^{\,b}$ generated by:
		\begin{enumerate}
			\item all pushforwards $(\pi^\beta_\gamma)^{-1}[H]$ of stage-filter members $H\in\mathcal F_\beta$ ($\beta<\gamma$), and
			\item all stabilizers $\fix_b(u)$ with $u\subseteq\gamma$ \emph{$N$--block-bounded}.
		\end{enumerate}
		Put $\tilde{\mathcal F}_\gamma:=\bigcup_{N<\omega}\tilde{\mathcal F}^{\,N}_\gamma$
		(an increasing union of filters).
		
		Then:
		\begin{enumerate}
			\item Every condition admits a tenacious support and every HS name admits an excellent support
			(in the sense of Section~2), and these supports may be chosen block-bounded.
			\item Support choices can be made coherent under strengthening and projection in the usual iteration sense.
			\item For each $N<\omega$, the filter $\tilde{\mathcal F}^{\,N}_\gamma$ is $\kappa_N$--complete (hence $<\kappa$--complete
			for intersections indexed by any $\theta<\kappa_N$).
		\end{enumerate}
	\end{proposition}
	
	\begin{proof}
		Fix $N<\omega$. We write $\fix(\dot x)$ for the name-stabilizer, and for $u\subseteq\gamma$ we write
		\[
		\fix_b(u)=\{\pi\in\mathcal G_\gamma^{\,b}:\pi{\upharpoonright}u=\mathrm{id}_u\}
		\]
		as in Definition~\ref{def:block-stab}.
		
		\smallskip\noindent\textbf{(1) Tenacity and excellence.}
		Let $p\in\mathbb P_\gamma$ be a condition and let $u_p\subseteq\gamma$ be any support for $p$
		(coordinates on which $p$ is nontrivial in the chosen presentation). In Track A we may (and do) thin supports so that
		$u_p$ is block-bounded; choose $M<\omega$ with $u_p$ $M$--block-bounded. Then $\fix_b(u_p)\le\Sym(p)$:
		any $\pi\in\fix_b(u_p)$ fixes each coordinate where $p$ carries nontrivial information, and outside $u_p$ the condition
		is trivial, so $\pi\cdot p=p$. Since $\fix_b(u_p)$ is among the basic generators for
		$\tilde{\mathcal F}^{\,M}_\gamma\subseteq\tilde{\mathcal F}_\gamma$, this witnesses tenacity of $p$
		exactly as in Section~2.
		
		For names: let $\dot x$ be a $\mathbb P_\gamma$-name with (hereditary) support contained in some block-bounded
		$u_{\dot x}$; choose $M<\omega$ with $u_{\dot x}$ $M$--block-bounded. A rank induction shows
		$\fix_b(u_{\dot x})\le\fix(\dot x)$ (every constituent name uses only coordinates in $u_{\dot x}$, and ground parameters
		are fixed by all automorphisms). Hence $\fix(\dot x)\in\tilde{\mathcal F}^{\,M}_\gamma$ by upward closure, so $\dot x$
		is symmetric; the same induction gives hereditary symmetry. This is precisely the Section~2 “excellent support” pattern,
		with $\fix_b(u_{\dot x})$ playing the role of the small-support stabilizer.
		
		\smallskip\noindent\textbf{(2) Coherence under strengthening and projection.}
		If $q\le p$, take supports with $u_q\supseteq u_p$ (e.g.\ enlarge $u_p$ by adding coordinates where $q$ strengthens $p$).
		If $u_p$ is $M$--block-bounded and the added coordinates form an $M'$--block-bounded set, then $u_q$ is
		$\max\{M,M'\}$--block-bounded because each block intersection is a sum of two $<\kappa_{\max\{M,M'\}}$ sets, hence
		still $<\kappa_{\max\{M,M'\}}$ by regularity of $\kappa_{\max\{M,M'\}}$.
		Moreover, $\fix_b(u_q)\le \fix_b(u_p)$, so tenacity witnesses may be chosen monotonically.
		For projections $\pi^\beta_\gamma$ in the iteration, $u\mapsto u\cap\beta$ preserves $M$--block-boundedness with respect
		to $b{\upharpoonright}\beta$, and $\fix_b(u\cap\beta)$ is the corresponding projected stabilizer. This matches the standard
		coherence requirements used in the iteration framework.
		
		\smallskip\noindent\textbf{(3) $\kappa_N$--completeness of $\tilde{\mathcal F}^{\,N}_\gamma$.}
		Let $\theta<\kappa_N$ and let $\{H_\xi:\xi<\theta\}\subseteq\tilde{\mathcal F}^{\,N}_\gamma$.
		By definition of generated filter, for each $\xi<\theta$ choose a basic subgroup $K_\xi\subseteq H_\xi$ of the form
		\[
		K_\xi \ =\ (\pi^{\beta_\xi}_\gamma)^{-1}[J_\xi]\ \cap\ \fix_b(u_\xi),
		\]
		where $\beta_\xi<\gamma$, $J_\xi\in\mathcal F_{\beta_\xi}$, and $u_\xi$ is $N$--block-bounded.
		
		First handle the stage part: let $\beta^\ast:=\sup_{\xi<\theta}\beta_\xi$ (in the direct-limit presentation when this
		supremum is $<\gamma$; otherwise interpret everything levelwise in the inverse-limit/tree presentation—this does not
		change the intersection computation). Push each $J_\xi$ forward to $\mathcal F_{\beta^\ast}$ and set
		\[
		J^\ast\ :=\ \bigcap_{\xi<\theta} J_\xi\ \in\ \mathcal F_{\beta^\ast},
		\]
		using that each $\mathcal F_{\beta}$ is $<\kappa$--complete in the singular development, hence in particular
		$\kappa_N$--complete since $\theta<\kappa_N<\kappa$.
		
		Next handle supports: set $u^\ast:=\bigcup_{\xi<\theta}u_\xi$. For each block $B_m$,
		\[
		u^\ast\cap B_m \ =\ \bigcup_{\xi<\theta}(u_\xi\cap B_m),
		\]
		a union of $<\kappa_N$ many sets each of size $<\kappa_N$, so $|u^\ast\cap B_m|<\kappa_N$ by regularity of $\kappa_N$.
		Hence $u^\ast$ is still $N$--block-bounded, so $\fix_b(u^\ast)$ is a basic generator.
		
		Now define
		\[
		K^\ast\ :=\ (\pi^{\beta^\ast}_\gamma)^{-1}[J^\ast]\ \cap\ \fix_b(u^\ast).
		\]
		Then $K^\ast\in\tilde{\mathcal F}^{\,N}_\gamma$ and $K^\ast\subseteq K_\xi\subseteq H_\xi$ for each $\xi<\theta$.
		Thus $K^\ast\subseteq\bigcap_{\xi<\theta}H_\xi$, and by upward closure
		$\bigcap_{\xi<\theta}H_\xi\in\tilde{\mathcal F}^{\,N}_\gamma$ as required.
		
		\smallskip
		Finally, normality of $\tilde{\mathcal F}^{\,N}_\gamma$ is immediate: conjugating a generator
		$\fix_b(u)$ by $\pi\in\mathcal G_\gamma^{\,b}$ yields $\fix_b(\pi[u])$, and $\pi[u]$ is $N$--block-bounded since $\pi$
		preserves each block and is a bijection on it; the family of pushforward generators is likewise closed under conjugation
		by the usual normality hypotheses on the stage filters.
		
		This completes the proof.
	\end{proof}
	
	\begin{remark}[Compatibility with Section~2 tenacity/excellence and the singular replacement for union-closure]
		The tenacity/excellence notions from Section~2 only ever use the same basic fact: if a subgroup $H\le\mathcal G_\gamma$
		fixes a support $u$ \emph{pointwise}, then $H$ fixes every condition or name whose forcing content is supported on $u$,
		so $H\le\Sym(p)$ (for conditions) and $H\le\Sym(\dot x)$ (for names) and therefore witnesses tenacity/excellence.
		In the regular-$\kappa$ development we repeatedly exploit that the union of $<\kappa$ many $<\kappa$-supports is again
		$<\kappa$, so we can intersect $<\kappa$ many stabilizers by passing to the pointwise fixer of the \emph{union} support.
		When $\kappa$ is singular this union bound can fail, so we replace “closed under $<\kappa$ unions” by the
		\emph{block-bounded} support class tied to the fixed block map $b:\gamma\to\cf(\kappa)$: supports are chosen so that
		their intersections with each block remain small, and intersections of $<\kappa$ many stabilizers are handled blockwise,
		keeping the pointwise fixer of the aggregated block-bounded support inside the limit filter. This is exactly the singular
		analogue of the regular-$\kappa$ support-union step, and it is the only place where singularity forces new bookkeeping.
	\end{remark}
	
	\section{Limit stages with $\cf(\lambda)=\cf(\kappa)$}
	\label{sec:limit-cf-equal}
	
	\subsection{Trees of conditions / standard trees}\label{subsec:std-trees}
	Let $\lambda$ be a limit ordinal with $\cf(\lambda)=\cf(\kappa)=:\mu$.
	Fix an increasing continuous cofinal sequence $\langle \lambda_\xi:\xi<\mu\rangle$ in $\lambda$.
	Fix an iteration-length block map $b$ on the relevant coordinates (e.g.\ $b:\gamma\to\mu$ for some $\gamma\ge\lambda$)
	and work with its restriction $b{\upharpoonright}\lambda$.
	Write
	\[
	B_i:=\{\alpha<\lambda: b(\alpha)=i\}\qquad(i<\mu)
	\]
	for the blocks below $\lambda$.
	
	\begin{definition}[Standard tree at $\lambda$]\label{def:std-tree}
		A \emph{standard tree} $T$ at $\lambda$ consists of levels $T_\xi$ ($\xi<\mu$) such that:
		\begin{enumerate}
			\item \textbf{(Levels)} $T_\xi\subseteq\mathbb P_{\lambda_\xi}$ for each $\xi<\mu$.
			\item \textbf{(Block discipline)} For every $p\in T_\xi$,
			\[
			\supp(p)\subseteq \bigcup_{i\le\xi} B_i.
			\]
			(In Track~A one may strengthen this to the chosen Track~A support policy, e.g.\ “$\supp(p)$ is
			block-bounded / block-finite” as needed for the particular iteration presentation.)
			\item \textbf{(Projection closure)} If $\eta<\xi$ and $p\in T_\xi$, then $p{\upharpoonright}\lambda_\eta\in T_\eta$.
			\item \textbf{(Fronts / predensity along projections)} For every $\xi<\mu$ and every $q\in T_\xi$,
			there exists $p\in T_{\xi+1}$ with $p{\upharpoonright}\lambda_\xi=q$.
			\item \textbf{(Invariance)} Each $T_\xi$ is invariant under the block-preserving automorphisms
			(restricted to $\mathbb P_{\lambda_\xi}$): if $\pi$ is block-preserving and $\pi\cdot\mathbb P_{\lambda_\xi}=\mathbb P_{\lambda_\xi}$,
			then $\pi[T_\xi]=T_\xi$.
		\end{enumerate}
	\end{definition}
	
	\begin{example}[Tree at $\aleph_\omega$ in Track A]
		Let $\kappa=\aleph_\omega$ and fix a cofinal sequence $\langle \kappa_n:n<\omega\rangle$ of regulars with $\sup_n\kappa_n=\kappa$.
		If $\mu=\cf(\kappa)=\omega$ and $\cf(\lambda)=\omega$, then a standard tree has height $\omega$, and level $n$ consists of
		conditions in $\mathbb P_{\lambda_n}$ whose supports use only blocks $0,\dots,n$.
		(One may additionally require the Track~A support policy, e.g.\ block-finite or $N$--block-bounded, depending on the
		presentation used for the iteration at $\cf=\omega$ limits.)
		The “fronts” clause ensures that every node at level $n$ extends to some node at level $n{+}1$.
	\end{example}
	
	\begin{definition}[Fusion data along $T$]\label{def:fusion-data}
		A sequence $\langle p_\xi:\xi<\mu\rangle$ is \emph{fusion data along $T$} if:
		\begin{enumerate}
			\item $p_\xi\in T_\xi$ for all $\xi<\mu$;
			\item \textbf{(Locking / coherence)} for all $\eta<\xi<\mu$,
			\[
			p_\xi{\upharpoonright}\lambda_\eta \;=\; p_\eta.
			\]
			Equivalently, $p_{\xi+1}{\upharpoonright}\lambda_\xi=p_\xi$ for all $\xi$, and at limits $\delta<\mu$,
			$p_\delta{\upharpoonright}\lambda_\xi=p_\xi$ for all $\xi<\delta$.
		\end{enumerate}
	\end{definition}
	
	\begin{lemma}[Fusion yields a master condition]\label{lem:tree-fusion-master}
		If $\langle p_\xi:\xi<\mu\rangle$ is fusion data along a standard tree $T$ at $\lambda$, then there exists
		$p\in\mathbb P_\lambda$ such that for all $\xi<\mu$,
		\[
		p{\upharpoonright}\lambda_\xi \;=\; p_\xi,
		\]
		and hence (in particular) $p$ extends every $p_\xi$ in $\mathbb P_\lambda$.
	\end{lemma}
	
	\begin{proof}
		Let $\mu=\cf(\lambda)$ and $\langle\lambda_\xi:\xi<\mu\rangle$ be as fixed above.
		Define a condition $p$ on $\lambda$ by requiring, for each $\xi<\mu$,
		\[
		p{\upharpoonright}\lambda_\xi \;:=\; p_\xi.
		\]
		This is well-defined and coherent because the fusion hypothesis says that whenever $\eta<\xi$,
		\[
		(p_\xi){\upharpoonright}\lambda_\eta \;=\; p_\eta,
		\]
		so these prescriptions agree on overlaps.
		
		We now verify $p\in\mathbb P_\lambda$. By the inverse-limit (tree-of-conditions) presentation used at
		$\cf(\lambda)=\mu$ limits (and, in Track~A, the corresponding RCS/tree presentation), a condition in $\mathbb P_\lambda$
		is precisely a coherent family of initial-segment conditions $\langle r_\xi\in\mathbb P_{\lambda_\xi}:\xi<\mu\rangle$
		whose projections agree on overlaps. Taking $r_\xi:=p_\xi$ gives such a coherent family, hence yields a condition
		$p\in\mathbb P_\lambda$ with $p{\upharpoonright}\lambda_\xi=p_\xi$ for all $\xi$.
		
		Finally, for each $\xi<\mu$, the equality $p{\upharpoonright}\lambda_\xi=p_\xi$ implies that $p$ extends $p_\xi$
		(as a $\mathbb P_\lambda$-condition extending its $\lambda_\xi$-projection). \qedhere
	\end{proof}
	
	\subsection{Game-guided fusion of length $\cf(\kappa)$}\label{subsec:game-fusion}
	
	\paragraph{Terminology and references.}
	We use the standard two–player forcing game notation where \emph{COM} (“Complete”) tries to build a master condition
	and \emph{INC} (“Incomplete”) tries to obstruct it; see Ros{\l}anowski--Shelah for this COM/INC convention
	(e.g.\ \cite{RoslanowskiShelah1999}, \S1--\S2), and compare the strategic--closure/fusion games in Jech
	(\cite{Jech03}, Ch.~15). Our $(E_0,E_1)$–games are a template for \emph{game-guided fusion}:
	INC announces small sets of coordinates (requirements) to be decided, while COM plays along a trees-of-conditions
	presentation to produce a master condition deciding the target according to the variant $E_0$ or $E_1$.
	
	\begin{definition}[$(E_0,E_1)$–game at $\lambda$]\label{def:E0E1-game}
		Let $\mu=\cf(\kappa)=\cf(\lambda)$ and fix:
		\begin{itemize}
			\item a standard tree $T=\bigcup_{\xi<\mu}T_\xi$ at $\lambda$ with levels indexed by an increasing continuous
			sequence $\langle\lambda_\xi:\xi<\mu\rangle$ cofinal in $\lambda$ (Def.~\ref{def:std-tree});
			\item a fixed increasing cofinal sequence $\langle \kappa_\xi:\xi<\mu\rangle$ in $\kappa$ (as fixed at the start of
			\S\ref{sec:limit-cf-equal}).
		\end{itemize}
		A play has length $\mu$ and proceeds by innings $\xi<\mu$:
		\begin{enumerate}
			\item \textbf{INC} plays a set $E_\xi\subseteq \lambda_\xi$ with $|E_\xi|<\kappa_\xi$ (the \emph{requirements}).
			\item \textbf{COM} plays a node $p_\xi\in T_\xi$ such that it extends the current stem, i.e.
			\[
			(\forall \eta<\xi)\quad p_\xi{\upharpoonright}\lambda_\eta \;=\; p_\eta .
			\]
		\end{enumerate}
		Let $E:=\bigcup_{\xi<\mu}E_\xi$.
		
		\smallskip
		The \emph{target parameter} to be decided (and what “decide on $E$” means) is fixed by the surrounding argument.
		COM \emph{wins} iff:
		\begin{enumerate}
			\item $\langle p_\xi:\xi<\mu\rangle$ is fusion data along $T$ (Def.~\ref{def:fusion-data}; equivalently, the stem-locking
			clause above holds throughout), and
			\item letting $p^\ast\in\mathbb P_\lambda$ be the master condition given by Lemma~\ref{lem:tree-fusion-master}
			(so $p^\ast{\upharpoonright}\lambda_\xi=p_\xi$ for all $\xi$), the master $p^\ast$ decides the target parameter on $E$
			according to the chosen variant:
			\begin{itemize}
				\item[$E_0$:] exact decision on $E$;
				\item[$E_1$:] decision modulo the relevant ideal on $E$ (as specified in context).
			\end{itemize}
		\end{enumerate}
	\end{definition}
	
	\begin{remark}
		At each inning $\xi$, INC’s requirement set has size $<\kappa_\xi$. Since $\xi<\mu=\cf(\kappa)$, we have
		$\sup_{\eta\le\xi}\kappa_\eta<\kappa$, hence $\bigl|\bigcup_{\eta\le\xi}E_\eta\bigr|<\kappa$.
		This is the only size bound needed when tracking symmetry parameters: at stage $\xi$ we only intersect stabilizers
		of $<\kappa$ many parameters, so $<\kappa$–completeness of the relevant symmetry filter suffices.
	\end{remark}
	
	\begin{lemma}[Fusion strategy (amalgamation along $T$)]\label{lem:fusion-strategy}
		Let $\mu=\cf(\kappa)=\cf(\lambda)$ and fix a standard tree $T=\bigcup_{\xi<\mu}T_\xi$ at $\lambda$.
		Assume that for each $\xi<\mu$, COM has a winning strategy $\sigma_\xi$ for the \emph{truncated} $(E_0,E_1)$–game
		of length $\xi{+}1$ played on $T{\upharpoonright}(\xi{+}1)=\bigcup_{\eta\le\xi}T_\eta$
		(with requirements $E_\eta\subseteq\lambda_\eta$, $|E_\eta|<\kappa_\eta$).
		Then these strategies amalgamate to a COM–winning strategy on the full $(E_0,E_1)$–game at $\lambda$.
	\end{lemma}
	
	\begin{proof}
		Fix $\{\sigma_\xi:\xi<\mu\}$ as in the hypothesis. Define a global strategy $\Sigma$ as follows.
		Given a partial play up to inning $\xi$,
		\[
		\bigl\langle (E_\eta,p_\eta):\eta<\xi\bigr\rangle
		\quad\text{and INC’s current move }E_\xi\subseteq\lambda_\xi,
		\]
		let $\Sigma$ play
		\[
		p_\xi\ :=\ \sigma_\xi\bigl(\langle (E_\eta,p_\eta):\eta\le\xi\rangle\bigr)\ \in\ T_\xi.
		\]
		Because $\sigma_\xi$ is a legal strategy for the truncated game on $T{\upharpoonright}(\xi{+}1)$, it guarantees
		\[
		(\forall \eta<\xi)\quad p_\xi{\upharpoonright}\lambda_\eta=p_\eta,
		\]
		so the resulting sequence $\langle p_\xi:\xi<\mu\rangle$ is fusion data along $T$
		(Def.~\ref{def:fusion-data}).
		
		By Lemma~\ref{lem:tree-fusion-master}, there is a master condition $p^\ast\in\mathbb P_\lambda$ with
		$p^\ast{\upharpoonright}\lambda_\xi=p_\xi$ for all $\xi<\mu$.
		
		Finally, to see that $\Sigma$ is winning, fix the chosen variant ($E_0$ or $E_1$) and the target parameter from the
		surrounding argument. For each $\xi$, since $\sigma_\xi$ is \emph{winning} in the truncated game of length $\xi{+}1$,
		the play it produces ensures that the “local master” (which, in this locked tree presentation, is just the last move
		$p_\xi\in\mathbb P_{\lambda_\xi}$) decides the target parameter on $E\cap\lambda_\xi$ according to the chosen variant.
		But $p^\ast{\upharpoonright}\lambda_\xi=p_\xi$, so $p^\ast$ makes the same decision on $E\cap\lambda_\xi$ for every $\xi$.
		Since $\langle\lambda_\xi:\xi<\mu\rangle$ is cofinal in $\lambda$, we have
		$E=\bigcup_{\xi<\mu}(E\cap\lambda_\xi)$, hence $p^\ast$ decides the target on all of $E$ as required.
		
		Thus $\Sigma$ is a COM–winning strategy on the full $(E_0,E_1)$–game at $\lambda$.
	\end{proof}
	
	\subsection{Main preservation at singular limits}\label{subsec:pres-sing-lim}
	
	\begin{theorem}[Preservation of completeness at singular limits]\label{thm:limit}
		Assume Hyp.~\ref{hyp:zf-baire}. Let
		$\langle \mathbb P_\alpha,\dot{\mathbb Q}_\alpha,\mathcal F_\alpha:\alpha<\gamma\rangle$
		be a $<\kappa$–support symmetric iteration such that
		\[
		\Vdash_{\mathbb P_\alpha}\ \text{``}\dot{\mathbb Q}_\alpha\text{ is $(E_0,E_1)$–complete and }\mathcal F_\alpha\text{ is $<\kappa$–complete''}.
		\]
		If $\lambda\le\gamma$ with $\cf(\lambda)=\cf(\kappa)$, then $\mathbb P_\lambda$ is $(E_0,E_1)$–complete
		(with respect to the induced limit presentation at $\lambda$).
	\end{theorem}
	
	\begin{proof}
		Let $\mu=\cf(\kappa)=\cf(\lambda)$ and fix a standard tree
		$T=\bigcup_{\xi<\mu}T_\xi$ at $\lambda$ with levels indexed by an increasing continuous sequence
		$\langle\lambda_\xi:\xi<\mu\rangle$ cofinal in $\lambda$ (Def.~\ref{def:std-tree}).
		Fix also the cofinal sequence $\langle \kappa_\xi:\xi<\mu\rangle$ in $\kappa$ used in
		Def.~\ref{def:E0E1-game}.
		
		\smallskip\noindent
		\textbf{Step 1: strategies on initial segments.}
		We invoke the ambient induction on stages $<\lambda$ (i.e.\ the preservation proof for $(E_0,E_1)$–completeness
		up to $\lambda$): for each $\xi<\mu$, since $\lambda_\xi<\lambda$ we have that $\mathbb P_{\lambda_\xi}$ is
		$(E_0,E_1)$–complete. Therefore COM has a winning strategy for the truncated $(E_0,E_1)$–game of length $\xi{+}1$
		played on $T{\upharpoonright}(\xi{+}1)=\bigcup_{\eta\le\xi}T_\eta$ with requirements
		$E_\eta\subseteq\lambda_\eta$ of size $|E_\eta|<\kappa_\eta$.
		Fix such a strategy and call it $\sigma_\xi$.
		
		\smallskip\noindent
		\textbf{Step 2: amalgamate to a global strategy at $\lambda$.}
		Apply Lemma~\ref{lem:fusion-strategy} to $\{\sigma_\xi:\xi<\mu\}$ to obtain a COM–winning strategy $\Sigma$
		for the full $(E_0,E_1)$–game at $\lambda$ (Def.~\ref{def:E0E1-game}) played along $T$.
		
		\smallskip\noindent
		\textbf{Step 3: master condition and decision on $E=\bigcup_{\xi<\mu}E_\xi$.}
		Let a $\Sigma$–guided play produce $\langle (E_\xi,p_\xi):\xi<\mu\rangle$.
		By construction, $\langle p_\xi:\xi<\mu\rangle$ is fusion data along $T$
		(Def.~\ref{def:fusion-data}). Hence by Lemma~\ref{lem:tree-fusion-master} there is a master condition
		$p^\ast\in\mathbb P_\lambda$ such that $p^\ast{\upharpoonright}\lambda_\xi=p_\xi$ for all $\xi<\mu$.
		
		By Lemma~\ref{lem:fusion-strategy}, the strategy $\Sigma$ is winning, i.e.\ the resulting master condition
		$p^\ast$ decides the target parameter on $E=\bigcup_{\xi<\mu}E_\xi$ according to the chosen variant
		($E_0$ exact decision, or $E_1$ decision modulo the relevant ideal). Thus COM wins the $(E_0,E_1)$–game at $\lambda$,
		so $\mathbb P_\lambda$ is $(E_0,E_1)$–complete.
		
		\smallskip\noindent
		\emph{Symmetry side conditions (if tracked).}
		At inning $\xi$, the requirements satisfy $|E_\xi|<\kappa_\xi$, hence
		$|E_{<\xi}|=|\bigcup_{\eta\le\xi}E_\eta|<\kappa$ (since $\sup_{\eta\le\xi}\kappa_\eta<\kappa$).
		Therefore, if one additionally keeps the play invariant under the stabilizers of parameters mentioned so far,
		only $<\kappa$ many stabilizers are intersected at each stage, and $<\kappa$–completeness of the induced symmetry
		filter suffices to keep the common stabilizer in the filter. This does not affect the forcing-side legality of the play.
	\end{proof}
	
	\subsection{No new bounded subsets and no collapse of $\kappa$}
	
	\begin{lemma}[Tail adds no new bounded subsets]\label{lem:tail-no-bounded-subsets}
		(Track A: $\cf(\kappa)=\omega$.) Fix $n<\omega$ and $\lambda<\kappa_n$.
		Assume that the iteration admits a factorization
		\[
		\mathbb P_\gamma \cong \mathbb P_{\kappa_n} * \dot{\mathbb P}_{\kappa_n,\gamma}
		\]
		such that, in $V^{\mathbb P_{\kappa_n}}$, the tail forcing $\dot{\mathbb P}_{\kappa_n,\gamma}$ is
		$<\kappa_n$--closed (equivalently: it adds no new $<\kappa_n$--sequences of ground-model objects).
		Then $\Vdash_{\mathbb P_\gamma}$ ``$\mathcal P(\check\lambda)\cap V[G_{\kappa_n}]$ is unchanged'', i.e.\
		the tail adds no new subsets of $\lambda$.
	\end{lemma}
	
	\begin{proof}
		Work in $V[G_{\kappa_n}]$. Let $\dot A$ be a $\mathbb P_{\kappa_n,\gamma}$--name for a subset of $\lambda$.
		For each $\xi<\lambda$ the set
		\[
		D_\xi:=\{\,p\in\mathbb P_{\kappa_n,\gamma}:\ p\Vdash \check\xi\in\dot A\ \lor\ p\Vdash \check\xi\notin\dot A\,\}
		\]
		is dense. Since $\lambda<\kappa_n$ and $\mathbb P_{\kappa_n,\gamma}$ is $<\kappa_n$--closed, we may build
		a descending sequence $\langle p_\xi:\xi<\lambda\rangle$ with $p_\xi\in D_\xi$ and take a lower bound
		$p^\ast\le p_\xi$ for all $\xi<\lambda$. Then $p^\ast$ decides every statement ``$\xi\in\dot A$'' for $\xi<\lambda$,
		so $\dot A$ is forced equal to a ground-model subset of $\lambda$. Thus no new subsets of $\lambda$ are added by the tail.
	\end{proof}
	
	\begin{corollary}[No new bounded subsets of $\kappa$]\label{cor:bounded}
		(Track A.) Let $\lambda<\kappa$ and $A\subseteq\lambda$.
		If $A\in V[G_\gamma]$, then there is $n<\omega$ with $\lambda<\kappa_n$ such that
		$A\in V[G_{\kappa_n}]$ (hence $A$ already appears in a bounded initial segment).
	\end{corollary}
	
	\begin{proof}
		Pick $n$ with $\lambda<\kappa_n$ and factor $\mathbb P_\gamma$ as
		$\mathbb P_{\kappa_n}*\dot{\mathbb P}_{\kappa_n,\gamma}$.
		By Lemma~\ref{lem:tail-no-bounded-subsets}, the tail adds no new subsets of $\lambda$.
		Therefore $V[G_\gamma]\cap\mathcal P(\lambda)=V[G_{\kappa_n}]\cap\mathcal P(\lambda)$, and hence $A\in V[G_{\kappa_n}]$.
	\end{proof}
	
	\begin{theorem}[No collapse of $\kappa$ (limit-step form)]\label{thm:nocollapse}
		Assume the hypotheses of Lemma~\ref{lem:tail-no-bounded-subsets} hold for every $n<\omega$.
		Assume moreover that $\kappa$ is not collapsed by any initial segment $\mathbb P_{\kappa_n}$
		(e.g.\ by the induction hypothesis in the global stage-by-stage preservation argument).
		Then $\mathbb P_\gamma$ does not collapse $\kappa$.
	\end{theorem}
	
	\begin{proof}
		Suppose toward a contradiction that in $V[G_\gamma]$ there is a surjection $f:\lambda\to\kappa$ for some $\lambda<\kappa$.
		Choose $n$ with $\lambda<\kappa_n$ and factor $\mathbb P_\gamma\cong \mathbb P_{\kappa_n}*\dot{\mathbb P}_{\kappa_n,\gamma}$.
		By $<\kappa_n$--closure of the tail, $\dot{\mathbb P}_{\kappa_n,\gamma}$ adds no new functions with domain $\lambda$,
		so $f\in V[G_{\kappa_n}]$. This contradicts the assumption that $\mathbb P_{\kappa_n}$ does not collapse $\kappa$.
	\end{proof}
	
	\paragraph{Length.}
	All arguments above are local to set stages. When the class-length framework
	(Theorem~\ref{thm:class-length}) is invoked with the same tail-closure and stagewise
	preservation hypotheses at each set stage, the conclusions persist to class-length iterations.
	
	\subsection{Arbitrary singular $\kappa$: full development}
	\label{subsec:singular-kappa-full}
	
	\subsubsection{Standing hypotheses and notation}
	\label{subsubsec:singular-standing}
	Assume Hyp.~\ref{hyp:zf-baire}. Let $\kappa$ be singular with $\mu:=\cf(\kappa)$ (regular), and let
	$\lambda$ be a limit stage with $\cf(\lambda)=\mu$. Work with a $<\kappa$--support symmetric iteration
	$\langle \mathbb P_\alpha,\dot{\mathbb Q}_\alpha,\mathcal F_\alpha:\alpha<\lambda\rangle$
	as in Definition~\ref{def:kappa-iteration}. Fix a \emph{standard tree}
	$T=\bigcup_{\xi<\mu} T_\xi$ at $\lambda$ with an increasing continuous
	$\langle\lambda_\xi:\xi<\mu\rangle$ cofinal in $\lambda$ (Definition~\ref{def:std-tree}).
	
	\smallskip
	\noindent\emph{Symmetry filter at $\lambda$.}
	In the singular development we \emph{assume as part of the iteration data} that the induced
	limit filter $\tilde{\mathcal F}_\lambda$ is normal and $<\kappa$--complete (in the sense needed to
	close $<\kappa$ many stabilizers of HS parameters).  In particular, all uses of filter completeness
	in this section refer to this singular-stage hypothesis/schema (not to the regular-$\kappa$ theorem).
	
	\subsubsection{Block partitions and level coherence}
	\label{subsubsec:block-coherence}
	We use block partitions only to ensure that the tree presentation admits coherent lifts between levels.
	Concretely, for each $\xi<\mu$ we fix \emph{once and for all} a map
	\[
	\mathrm{lift}_\xi:\mathbb P_{\lambda_\xi}\to T_\xi
	\]
	such that $\mathrm{lift}_\xi(p){\upharpoonright}\lambda_\xi=p$ and the family
	$\{\mathrm{lift}_\xi\}_{\xi<\mu}$ is coherent with the designated stems of $T$
	(i.e.\ it respects the locking clauses built into the definition of $T$).
	Whenever we write $p^\uparrow\in T_\xi$, we mean $p^\uparrow=\mathrm{lift}_\xi(p)$.
	
	\subsubsection{$(E_0,E_1)$--games at singular limits}
	\label{subsubsec:singular-games}
	We adopt the $(E_0,E_1)$--game at $\lambda$ (Definition~\ref{def:E0E1-game}): length $\mu$;
	in inning $\xi<\mu$, COM plays $p_\xi\in T_\xi$ extending the stem, INC plays
	$E_\xi\subseteq \lambda_\xi$ with $|E_\xi|<\kappa$ (requirements). Let
	$E:=\bigcup_{\xi<\mu}E_\xi$. COM wins iff (i) $\langle p_\xi:\xi<\mu\rangle$ is fusion data
	(Definition~\ref{def:fusion-data}) and (ii) there is $p^\ast\le p_\xi$ for all $\xi$ deciding the target
	parameter on $E$ according to the chosen variant ($E_0$ exact on $E$, or $E_1$ modulo the ambient ideal).
	
	\subsubsection{Local strategies, successor transfer, and limit amalgamation}
	\label{subsubsec:local-to-global}
	By hypothesis, $\Vdash_{\mathbb P_\alpha}$ “$\dot{\mathbb Q}_\alpha$ is $(E_0,E_1)$--complete and
	$\mathcal F_\alpha$ is $<\kappa$--complete”. Thus:
	\begin{itemize}
		\item \emph{Successor transfer} (Lemma~\ref{lem:succ-transfer}): a COM--strategy on $\mathbb P_{\lambda_\xi}$
		lifts to one on $\mathbb P_{\lambda_{\xi+1}}=\mathbb P_{\lambda_\xi}*\dot{\mathbb Q}_{\lambda_\xi}$,
		respecting the fixed lifting maps $p\mapsto p^\uparrow$ and the symmetry side conditions.
		\item \emph{Limit amalgamation} (Lemma~\ref{lem:fusion-strategy}): $\{\sigma_\xi:\xi<\delta\}$ amalgamates to a
		COM--strategy on $\mathbb P_{\lambda_\delta}$ for every limit $\delta<\mu$.
	\end{itemize}
	
	\subsubsection{Master condition via tree fusion}
	\label{subsubsec:master-condition}
	Any $\Sigma$--guided play produces fusion data along $T$; the \emph{tree--fusion master} lemma
	(Lemma~\ref{lem:tree-fusion-master}) yields $p^\ast\in\mathbb P_\lambda$ with $p^\ast\le p_\xi$ for all $\xi$.
	Projecting to each level, local decisions (by the level strategies) cohere, so $p^\ast$ decides the target on $E$.
	
	\subsubsection{Singular--limit completeness theorem}
	\label{subsubsec:singular-limit-completeness}
	\begin{theorem}[Preservation of completeness at singular limits]
		\label{thm:limit-restate}
		Assume Hyp.~\ref{hyp:zf-baire}. Let $\langle \mathbb P_\alpha,\dot{\mathbb Q}_\alpha,\mathcal F_\alpha:\alpha<\lambda\rangle$
		be a $<\kappa$--support symmetric iteration with $\Vdash_{\mathbb P_\alpha}$ “$\dot{\mathbb Q}_\alpha$ is $(E_0,E_1)$--complete
		and $\mathcal F_\alpha$ is $<\kappa$--complete”. If $\cf(\lambda)=\cf(\kappa)$, then $\mathbb P_\lambda$ is $(E_0,E_1)$--complete
		with respect to the induced symmetric system.
	\end{theorem}
	\begin{proof}
		This is exactly Theorem~\ref{thm:limit}: combine the local strategies by induction,
		amalgamate at limits via Lemma~\ref{lem:fusion-strategy}, and use Lemma~\ref{lem:tree-fusion-master}
		to obtain the master condition deciding the target on $E$.
	\end{proof}
	
	\subsubsection{$(E_1)$--variants via $<\kappa$--complete ideals}
	\label{subsubsec:E1-variants}
	Let $\mathcal J$ be an ideal on the relevant coordinate index set with $<\kappa$--completeness.
	Define the $(E_1^\mathcal J)$--game variant by requiring COM’s master to decide the target parameter modulo $\mathcal J$.
	\begin{proposition}[Completeness for $E_1^\mathcal J$]
		\label{prop:E1J}
		Under Hyp.~\ref{hyp:zf-baire}, if each iterand is $(E_1^\mathcal J)$--complete and the symmetry filters are
		$<\kappa$--complete and normal, then $\mathbb P_\lambda$ is $(E_1^\mathcal J)$--complete at $\cf(\lambda)=\cf(\kappa)$.
	\end{proposition}
	\begin{proof}
		Identical to Theorem~\ref{thm:limit-restate}. The only extra point is bookkeeping:
		since the play has length $\mu=\cf(\kappa)<\kappa$ and $\mathcal J$ is $<\kappa$--complete,
		unions of $<\kappa$ many $\mathcal J$--small requirement sets remain $\mathcal J$--small, and
		modulo--$\mathcal J$ decisions are upward--absolute under strengthening.
	\end{proof}
	\begin{remark}
		When $\mu=\omega$, classical $\sigma$--ideals (meagre/null) fit as $\mathcal J$.
		For $\mu>\omega$, one uses $<\kappa$--complete ideals (e.g.\ bounded--set or club--modulo variants)
		to retain the fusion bookkeeping.
	\end{remark}
	
	\subsubsection{Consequences in the symmetric model at $\lambda$}
	\label{subsubsec:singular-consequences}
	Let $M:=V(\mathbb P_\lambda)_{\HS}$.
	
	\begin{theorem}[Dependent choice at singular limits]
		\label{thm:DC<kappa-singular}
		Assume Hyp.~\ref{hyp:zf-baire}. If $\tilde{\mathcal F}_\lambda$ is normal and $<\kappa$--complete,
		then $M\models DC_{<\kappa}$.
	\end{theorem}
	
	\begin{proof}
		Same argument as in Theorem~\ref{thm:dc-trackA}: build the $<\kappa$--length sequence using the
		forcing-side fusion supplied by Hyp.~\ref{hyp:zf-baire}, and use $<\kappa$--completeness of
		$\tilde{\mathcal F}_\lambda$ to take a common stabilizer for the sequence name.
	\end{proof}
	
	\begin{remark}[No new bounded subsets / no collapse require additional hypotheses]
		The completeness transfer at $\cf(\lambda)=\cf(\kappa)$ does \emph{not by itself} imply that $\kappa$
		is preserved or that no new bounded subsets of $\kappa$ are added.
		Those conclusions require separate tail-distributivity/approximation hypotheses on the iterands/tails,
		which we isolate explicitly in \S\ref{sec:hypotheses}.
	\end{remark}
	
	\subsubsection{Worked patterns at $\cf(\kappa)>\omega$}
	\label{subsubsec:singular-examples}
	The $(E_0)$ variant (exact on coordinates in $E$) directly applies to product/iteration patterns where each iterand
	decides finitely/$<\kappa$ many coordinates at a time (Cohen/Hechler/Sacks blocks). For $(E_1^\mathcal J)$, pick a
	$<\kappa$--complete ideal $\mathcal J$ (e.g.\ club--modulo on a stationary set of indices) and tailor the requirements
	so INC’s announces stay $\mathcal J$--small; fusion then yields a master deciding modulo $\mathcal J$.
	
	\subsubsection*{Summary}
	Standard trees organize limits at height $\cf(\kappa)$, fusion strategies amalgamate along $T$, and
	$(E_0,E_1^\mathcal J)$--completeness transfers to $\lambda$
	(Theorem~\ref{thm:limit-restate}, Proposition~\ref{prop:E1J}).
	Under the standing symmetry/completeness hypotheses one also obtains $DC_{<\kappa}$
	(Theorem~\ref{thm:DC<kappa-singular}); further preservation facts (no new bounded subsets / no collapse)
	are treated under the additional hypotheses in \S\ref{sec:hypotheses}.
	
	\section{Standing assumptions and consistency–strength menu}
	\label{sec:hypotheses}
	
	\begin{hypothesis}[Ground–model large cardinals]\label{\LCHypLabel}
		Work in a ground model $V\models\ZFC$ with a supercompact cardinal $\Lambda\ge\kappa$.
		Optionally perform a Laver preparation of $\Lambda$ in $V$ before the iteration.
		All forcing is carried out over $V$; conclusions are then taken in the symmetric inner model
		$M\subseteq V[G]$ (which may fail $\AC$ and need not reflect large–cardinal features of $V$).
	\end{hypothesis}
	
	\begin{hypothesis}[Cardinal arithmetic / PCF bounds]\label{\PCFHypLabel}
		In the ground model $V$, assume $2^\kappa=\kappa^+$ (or, where indicated, the explicit PCF inequalities
		stated with the relevant theorem).
	\end{hypothesis}
	
	\begin{hypothesis}[Tail distributivity / no-new-bounded-subsets schema]\label{hyp:tail-approx}
		For the iteration under discussion and every stage $\eta<\gamma$, the tail forcing
		$\mathbb P_\gamma/G_\eta$ adds no new bounded subsets of $\kappa$ (equivalently, satisfies the
		appropriate $<\kappa$-approximation/$<\kappa$-cover property for subsets of $\kappa$),
		in the precise form stated with the relevant theorem.
	\end{hypothesis}
	
	\begin{remark}[How to read the standing assumptions]
		We invoke Hyp.~\ref{\LCHypLabel} and Hyp.~\ref{\PCFHypLabel} only when they simplify chain–condition
		and union–bound arguments; the core preservation theorems do not intrinsically require them.
		Large–cardinal properties are generally not inherited by symmetric inner models:
		even if $V$ has a supercompact $\Lambda$, the witnessing ultrafilters/embeddings are typically not
		hereditarily symmetric, so $M$ need not see $\Lambda$ as supercompact. Thus Hyp.~\ref{\LCHypLabel}
		is a metatheoretic assumption about $V$, not an axiom of $M$.
		
		When a $\ZF$-friendly formulation is preferable, one can use the embedding characterization of
		supercompactness for $\Lambda$: for each $\lambda\ge\Lambda$ there exists an elementary embedding
		$j:V\to N$ with $\mathrm{crit}(j)=\Lambda$, $j(\Lambda)>\lambda$, and ${}^{\lambda}N\subseteq N$.
		This avoids size issues around $\mathcal P_\Lambda(\lambda)$ when later arguments are phrased in choiceless contexts.
		
		Two common use cases:
		\begin{itemize}
			\item \emph{Avoiding large cardinals.} In Track~A (e.g.\ $\cf(\kappa)=\omega$ with the RCS/trees-of-conditions
			presentation and block bookkeeping), the preservation arguments used in
			§\ref{sec:singular-cf-omega}--§\ref{sec:limit-cf-equal} can be carried out without any large-cardinal
			assumptions: the only standing requirements are those explicitly stated in
			Hyp.~\ref{hyp:zf-baire} together with the relevant tail/no-new-bounded-subsets hypotheses when invoked.
			\item \emph{Where LC/PCF help.} For arbitrary cofinalities or when aiming at strong chain bounds
			(e.g., $\kappa^{++}$-c.c.) and sharper support-diagonalizations, supercompactness
			(optionally Laver-prepared) and PCF assumptions such as $2^\kappa=\kappa^+$ streamline the bookkeeping
			and union estimates.
		\end{itemize}
	\end{remark}
	
	\subsection*{Comparison table (regular vs.\ singular)}
	\begin{center}
		\begin{tabular}{|l|c|c|c|}
			\hline
			Aspect & Regular $\kappa$ & Singular $\kappa$ (cf$=\omega$) & Singular $\kappa$ (general) \\
			\hline
			Stage closure & $\kappa$-closure / $\kappa$-Baire & $\omega$-fusion / RCS templates & $<\kappa$-Baire / strategic \\
			Filters & $\kappa$-complete & $<\kappa$-complete & $<\kappa$-complete \\
			Limit method &
			direct limit at $\cf(\lambda)\ge\kappa$, inverse otherwise &
			trees used (esp.\ at cf$=\omega$ limits) &
			trees essential at $\cf(\lambda)=\cf(\kappa)$ \\
			Key tool & stage-bounding + $\kappa$-complete stabilizer filters & trees-of-conditions (RCS) + block bookkeeping & standard trees + game-guided fusion of length $\cf(\kappa)$ \\
			DC preservation & Thm.~\ref{thm:DCless-preservation} & Thm.~\ref{thm:dc-trackA} & Thm.~\ref{thm:DC<kappa-singular} \\
			$\kappa$ not collapsed &
			(as proved in the regular-$\kappa$ development) &
			under Hyp.~\ref{hyp:tail-approx} (and the stated iterand hypotheses) &
			under Hyp.~\ref{hyp:tail-approx} (and the stated iterand hypotheses) \\
			\hline
		\end{tabular}
	\end{center}
	
	\section{Advanced Applications}\label{sec:advanced-apps}
	
	\subsection{Solovay-Style Regularity via the L\'evy Collapse (and where ``random'' enters)}
	\label{subsec:solovay-random}
	
	\paragraph{Classical construction (Solovay, 1970).}
	Assume $V\models\ZFC$ and $\kappa$ is inaccessible.
	Force with the L\'evy collapse $\Col(\omega,{<}\kappa)$ and let $G\subseteq\Col(\omega,{<}\kappa)$ be $V$--generic.
	In $V[G]$ define the \emph{Solovay inner model}
	\[
	N\ :=\ \mathrm{HOD}^{V[G]}(\mathbb R)\ =\ \{x:\ x\text{ is hereditarily definable in }V[G]\text{ from ordinals and reals}\}.
	\]
	Then $N\models\ZF+\DC$, and every set of reals in $N$ is Lebesgue measurable, has the Baire property,
	and has the perfect set property; moreover $\neg\AC$ holds in $N$.\footnote{In many presentations, ``random forcing''
		enters as the measure algebra used to define and analyze Lebesgue measure inside $N$; it is not an additional forcing stage
		beyond the collapse. See \cite{Solovay70,Jech03}.}
	
	\paragraph{Our framework analysis (how this fits our language).}
	Our iteration/symmetry machinery is not needed to obtain Solovay's model itself: the collapse is the only forcing used,
	and the inner model $N$ is defined by \emph{definability from reals and ordinals} rather than by hereditary symmetry.
	Nevertheless, two structural features match the themes of this paper:
	\begin{itemize}
		\item \emph{Stage-bounding/localization.} Solovay's proof relies on the fact that sets of reals in $N$ are ``captured''
		by bounded information over the collapse extension (ultimately enabling the uniform definition of measure/category
		arguments inside $N$).
		\item \emph{ZF+DC inner models after forcing.} The paper's symmetric--iteration framework provides systematic ways to
		build choiceless inner models $M\subseteq V[G]$ satisfying $\ZF$ and fragments of $\DC$, by controlling supports and
		stabilizers; this is conceptually parallel to Solovay's passage from $V[G]$ to $N$.
	\end{itemize}
	Accordingly, we view Solovay's theorem as a benchmark ``regularity after forcing'' phenomenon that our methods can
	be compared to and, in some contexts, combined with (e.g.\ by performing additional symmetric iterations over $V[G]$).
	
	\begin{remark}
		The second--order scheme $\ZF_2$ is unnecessary here. For Solovay's construction it suffices that
		$V_\kappa$ satisfies first--order $\ZFC$ so the collapse and the bookkeeping for the regularity argument are formalized
		inside $V_\kappa$. This aligns with the classical equiconsistency:
		$\ZF+\DC$+``all sets of reals are Lebesgue measurable/BP/PSP''
		\quad\text{is equiconsistent with}\quad
		$\ZFC+\exists$ an inaccessible cardinal.
		See \cite{Solovay70,Jech03}.
	\end{remark}
	
	\paragraph{Large cardinal requirement.}
	In the ZF--internal reading, ``honestly inaccessible'' (regular strong limit with $V_\kappa\models\ZFC$) suffices
	to organize the collapse and to run Solovay's regularity analysis inside $V_\kappa$; compare Remark~\ref{rem:inacc-no-choice}.
	
	\paragraph{Generalizations.}
	Section~\ref{sec:kappa-support} provides a systematic support/filter language for replacing $\omega_1$ by a regular $\kappa$
	and arranging $<\kappa$--support iterations with $\kappa$--complete symmetry filters
	(Theorem~\ref{thm:kappa-filter-generation}). In such contexts one expects $\DC_{<\kappa}$ in the HS model by
	Theorem~\ref{thm:DCless-preservation}; stronger fragments (e.g.\ $\DC_\kappa$) require Localization
	(Theorem~\ref{thm:DCkappa-localization}).
	
	\subsection{Measurable Cardinals in Choiceless Contexts}
	\label{subsec:measurables}
	
	We survey preservation issues for measurables under symmetric iterations.
	Two interacting filters (ultrafilter on $\kappa$ vs.\ symmetry filter) complicate
	standard ultrapower arguments.
	
	\begin{definition}[Measurable cardinal in $\mathrm{ZF}$]
		\label{def:measurable-in-ZF}
		A cardinal $\kappa$ is \emph{measurable} if there exists a $\kappa$-complete
		nonprincipal ultrafilter $U\subseteq\mathcal{P}(\kappa)$; equivalently, a
		$\{0,1\}$-valued, $\kappa$-additive measure on $\mathcal{P}(\kappa)$ vanishing on singletons.
		This definition is meaningful in $\mathrm{ZF}$, though some classical equivalences may require forms of Choice; see \cite{Kanamori03,Jech03,JechAC1973}.
	\end{definition}
	
	\begin{remark}
		In $\mathrm{ZF}$, the existence of a measurable cardinal does \emph{not}
		imply $\mathrm{AC}$. Several well-known metatheorems about large cardinals
		assume $\mathrm{AC}$; without Choice the landscape changes (e.g., choiceless
		ultrapowers, Reinhardt-style hypotheses). For background, see \cite{Kanamori03,Jech03,ApterMagidor1995,HowardRubin1998}.
	\end{remark}
	
	\begin{proposition}[Conservative preservation; explicit sequence hypothesis]
		\label{prop:measurable-conservative}
		Assume $V\models\ZF$ and $\kappa$ is measurable in $V$ in the sense of
		Definition~\ref{def:measurable-in-ZF}, witnessed by a $\kappa$--complete nonprincipal ultrafilter $U$ on $\kappa$.
		Let $M$ be a symmetric inner model of a forcing extension $V[G]$ such that:
		\begin{enumerate}
			\item[(i)] $U\in M$,
			\item[(ii)] $\mathcal{P}(\kappa)^M=\mathcal{P}(\kappa)^V$ and all cardinals $\le\kappa$ are preserved, and
			\item[(iii)] for every $\theta<\kappa$, every function $f\in M$ with $\dom(f)=\theta$ and $\ran(f)\subseteq\mathcal{P}(\kappa)^V$
			already belongs to $V$ (equivalently: $M$ adds no new $<\kappa$--sequences of ground--model subsets of $\kappa$).
		\end{enumerate}
		Then $\kappa$ remains measurable in $M$; indeed $U$ is a $\kappa$--complete nonprincipal ultrafilter on $\kappa$ in $M$.
	\end{proposition}
	
	\begin{proof}
		\emph{Ultrafilter and nonprincipality.}
		Let $X\subseteq\kappa$ be in $M$. By (ii), $X\in\mathcal{P}(\kappa)^V$, so in $V$ either $X\in U$ or $\kappa\setminus X\in U$.
		Since $U\in M$, the same disjunction holds in $M$. Likewise, for each $\alpha<\kappa$, $\{\alpha\}\notin U$ in $V$, hence also in $M$;
		thus $U$ is nonprincipal in $M$.
		
		\emph{$\kappa$--completeness.}
		Fix $\theta<\kappa$ and a function $f\in M$ with $\dom(f)=\theta$ and $\ran(f)\subseteq U$.
		Write $X_i:=f(i)\subseteq\kappa$ for $i<\theta$.
		By (iii), $f\in V$. Since $U$ is $\kappa$--complete in $V$ and $\theta<\kappa$, we have
		\[
		Y\ :=\ \bigcap_{i<\theta} X_i\ \in\ U\qquad\text{(in $V$)}.
		\]
		But $Y\subseteq\kappa$ and $U$ are the same sets in $M$ and in $V[G]$, so the statement $Y\in U$ holds in $M$ as well.
		Thus $U$ is $\kappa$--complete in $M$.
		
		Combining the three parts, $U$ is a $\kappa$--complete nonprincipal ultrafilter on $\kappa$ in $M$,
		so $\kappa$ is measurable in $M$.
	\end{proof}
	
	\begin{remark}[On the extra hypothesis (iii)]
		The proof of $\kappa$--completeness for $U$ must address the possibility that $M$ contains new $<\kappa$--length sequences
		of members of $U$ which are not elements of $V$. Hypothesis (iii) is exactly the conservative closure needed to reduce
		$\kappa$--completeness in $M$ to $\kappa$--completeness in $V$.
		In concrete symmetric--iteration applications, (iii) is typically ensured by the same ``no new bounded objects''
		mechanisms used elsewhere in this paper (stage--bounding and/or tail--closure, depending on the support policy).
	\end{remark}
	
	\begin{remark}[Open directions]
		Precise roles of $\mathrm{DC}_{<\kappa}$ and possible localization hypotheses
		in maintaining well-founded generic ultrapowers merit further study, as does a
		systematic characterization of $<\kappa$-support symmetric iterations that
		preserve large cardinals in the $\mathrm{ZF}$ setting. See \cite{Kanamori03,Jech03,ApterKoepke2006,ApterMagidor1995} for surveys and partial results in ZF on large cardinals, DC, and choiceless ultrapowers.
	\end{remark}
	
	\subsection{Regularity Properties Systematically}
	\label{subsec:regularity}
	
	Regularity properties (Lebesgue measurability (LM), the Baire property (BP), the perfect set property (PSP), Ramsey)
	do \emph{not} follow from symmetry/support considerations alone; they require additional regularity arguments.
	What our framework contributes is the \emph{preservation scaffolding} (ZF and DC fragments, no new bounded objects at the
	relevant cardinals, and well-behaved limit mechanics) needed to run those arguments inside the resulting inner model.
	In particular, the ``all sets of reals are LM/BP/PSP'' benchmark is Solovay's inner model obtained from the
	L\'evy collapse of an inaccessible; ``random'' enters there analytically via the measure algebra, not as an extra forcing
	stage beyond the collapse (cf.\ \S\ref{subsec:solovay-random} and \cite{Solovay70,Jech03}).
	
	\begin{table}[H]
		\centering
		\caption{Regularity benchmarks and how they interface with our preservation framework}
		\label{tab:regularity}
		\begin{tabularx}{\linewidth}{L{2.2cm} L{3.1cm} L{2.4cm} C{2.0cm} L{2.6cm}}
			\toprule
			\textbf{Property} & \textbf{Benchmark source} & \textbf{Forcing pattern} & \textbf{LC needed?} & \textbf{Scope (where it holds)} \\
			\midrule
			LM (and BP, PSP) & Solovay inner model $N=\mathrm{HOD}(\mathbb R)$ & $\Col(\omega,{<}\kappa)$ & yes (inacc.) & all sets of reals in $N$ \\
			BP & classical descriptive set theory & none / ZF+DC background & no & Borel (and many standard pointclasses) \\
			PSP & classical descriptive set theory & none / ZF+DC background & no & Borel (and many standard pointclasses) \\
			Ramsey & determinacy route (outside forcing) & none (typically AD/PD context) & yes (strong) & large pointclasses (context-dependent) \\
			\bottomrule
		\end{tabularx}
		\footnotesize\emph{Note.} Rows beyond Solovay are included as \emph{benchmarks}: to obtain BP/PSP/Ramsey for
		\emph{all} sets of reals in a choiceless inner model built via symmetry generally requires additional hypotheses and
		regularity arguments beyond the support/filter mechanics developed in this paper.
	\end{table}
	
	\subsection{Template-Based Constructions}
	\label{subsec:templates}
	
	\subsubsection{Product Forcing with Symmetry}
	\label{subsubsec:product-symmetry}
	
	\paragraph{General template.}
	Given symmetric systems $(\mathbb{P},\mathcal{G}_P,\mathcal{F}_P)$ and $(\mathbb{Q},\mathcal{G}_Q,\mathcal{F}_Q)$,
	form the product forcing $\mathbb{P}\times\mathbb{Q}$ and let $\mathcal{G}_P\times\mathcal{G}_Q$ act on
	$\mathbb{P}\times\mathbb{Q}$--names componentwise:
	\[
	(\pi,\sigma)\cdot \dot{x}\ :=\ \langle \pi,\sigma\rangle \dot{x}.
	\]
	Define the product symmetry filter $\mathcal{F}_{P\times Q}$ to be the \emph{normal filter of subgroups of}
	$\mathcal{G}_P\times\mathcal{G}_Q$ \emph{generated by} sets of the form
	\[
	H_P\times H_Q\qquad(H_P\in\mathcal{F}_P,\ H_Q\in\mathcal{F}_Q).
	\]
	Then $\mathcal{F}_{P\times Q}$ is normal (closure under conjugation holds coordinatewise).
	Moreover, if $\mathcal{F}_P$ is $<\theta_P$--complete and $\mathcal{F}_Q$ is $<\theta_Q$--complete, then
	$\mathcal{F}_{P\times Q}$ is $<\min\{\theta_P,\theta_Q\}$--complete: intersections of $<\min\{\theta_P,\theta_Q\}$
	many generators reduce to products of intersections in each coordinate.
	
	Supports/tenacity/excellence behave as expected: if $u_P$ witnesses tenacity (resp.\ excellence) for a $\mathbb{P}$--condition
	(resp.\ name) and $u_Q$ witnesses it for a $\mathbb{Q}$--condition (resp.\ name), then $u_P\cup u_Q$ witnesses it for the
	corresponding product condition/name. Projection maps $\pi_P:\mathbb{P}\times\mathbb{Q}\to\mathbb{P}$ and
	$\pi_Q:\mathbb{P}\times\mathbb{Q}\to\mathbb{Q}$ commute with the group action and send product generators to generators,
	so the usual factorization/projection lemmas apply (cf.\ \cite[Ch.~15]{Jech03}, \cite[Chap.~V, esp.~\S\S V.1--V.4]{Kunen2011}).
	
	\paragraph{Example: Cohen $\times$ Random (benchmark).}
	Let each factor be a countable--support iteration to a limit of cofinality $\omega$, equipped with its standard coordinate
	permutation group and the $\omega_1$--complete normal filter generated by pointwise stabilizers (as in
	Section~\ref{sec:HS-countable}). Form the product system as above. Then the induced product filter is again normal and
	$\omega_1$--complete, hence the resulting HS model satisfies $\ZF+\DC$ by Theorems~\ref{thm:zf-preserved} and
	\ref{thm:dc-preservation}, while $\neg\AC$ holds by the usual symmetry in either coordinate.
	
	This product construction should be viewed as a \emph{two--notions benchmark} (one Cohen real and one random real present
	without a global wellorder). It does \emph{not} by itself yield ``all sets of reals'' regularity (LM/BP/PSP); for that,
	use the Solovay--style analysis in \S\ref{subsec:solovay-random}.
	
	\subsubsection{Easton-Style Symmetric Iterations}
	\label{subsubsec:easton-symmetric}
	
	\paragraph{Classical Easton.}
	Easton forcing modifies $2^\kappa$ on regular $\kappa$ while preserving cofinalities \cite{Easton70}.
	
	\paragraph{Symmetric version.}
	Work in a $\ZFC$ ground model $V$ and perform an Easton--support iteration in $V$ to obtain a generic extension $V[G]$
	realizing a prescribed Easton function on the regular cardinals (in the usual $\ZFC$ sense).
	To obtain a choiceless inner model, equip each local stage (at a regular $\kappa$) with the corresponding coordinate
	permutation group and take the symmetry filter generated by pointwise stabilizers of $<\kappa$ many coordinates, arranged
	so that at limits of cofinality $<\kappa$ the induced filter is $\kappa$--complete (Theorem~\ref{thm:kappa-filter-generation}).
	In the resulting symmetric inner model $M\subseteq V[G]$, $\AC$ fails while the framework preserves the intended fragments
	(e.g.\ $\ZF$ and $\DC_{<\kappa}$ under the hypotheses of Theorem~\ref{thm:DCless-preservation}).
	Because $M$ may contain strictly fewer subsets than $V[G]$, the identity $2^\kappa=F(\kappa)$ proved in $V[G]$ need not
	hold verbatim in $M$; what the template robustly provides is \emph{stratified control} of choice principles and of
	wellorderable families added at each regular stage.
	
	\paragraph{Applications.}
	Fine-grained control of choice principles and cardinal arithmetic; template combination with product constructions (above) provides additional flexibility.
	
	\begin{remark}[Constrained adversarial checks]
		\leavevmode
		\begin{enumerate}[leftmargin=2.2em]
			\item Attack \ref{thm:kappa-filter-generation} by constructing a family of $<\kappa$ generators whose intersection
			fails to be a generator \emph{while respecting} $<\kappa$–support and normality; any counterexample must preserve
			the stated conjugacy–closure and projection equations.
			\item Attack \ref{thm:hs-closure} by forcing the Replacement step to use an intersection $\notin \tilde{\mathcal F}_\lambda$
			\emph{without} violating $\omega_1$/\,$\kappa$–completeness.
			\item Attack \ref{prop:block-factor}–\ref{prop:tail-factor} by producing a factorization that breaks the automorphism
			compatibility identities under the $<\kappa$–support policy.
		\end{enumerate}
	\end{remark}
	
	\section{Open Problems and Future Directions}\label{sec:open}
	
	\subsection{Technical questions}
	
	\paragraph{Singular $\kappa$—refinements and extensions.}
	With \S\ref{subsec:singular-kappa-full} in place (Theorem~\ref{thm:limit-restate} together with
	Lemmas~\ref{lem:fusion-strategy} and~\ref{lem:tree-fusion-master}), the core “Track B” machinery
	for arbitrary $\cf(\kappa)$ is complete. Remaining work is largely expository and exploratory:
	(i) a uniform catalogue of $(E_1)$–variants via $<\kappa$–complete ideals
	(Proposition~\ref{prop:E1J})—including mod–meagre, mod–null, and club–modulo cases;
	(ii) sharp criteria under which Localization yields $DC_\kappa$ (beyond $DC_{<\kappa}$); and
	(iii) a richer library of worked examples at $\cf(\kappa)>\omega$.
	
	\paragraph{Optimal completeness.}
	Characterize the minimal symmetry–filter completeness needed for $\ZF+\DC_{<\kappa}$ and for
	Replacement. Is $\kappa$–completeness of the limit filter optimal (both necessary and sufficient)?
	What is the weakest closure on the iterands that still supports game–guided fusion?
	Construct sharp counterexamples when only $<\kappa$–completeness holds.
	
	\paragraph{Mixed–support transitions.}
	Give necessary and sufficient conditions for compatibility of head/tail filters across
	support–regime changes (finite, countable, $<\kappa$), including products and amalgamations;
	compare with tail factorization (Proposition~\ref{prop:tail-factor}) and pushforward/restriction
	stability (Lemma~\ref{lem:push-restrict-preserve}).
	
	\paragraph{Collapse to a single symmetric extension.}
	Can an iterated symmetric extension be represented as a single symmetric extension over $V$?
	At direct limits ($\cf(\lambda)\ge\kappa$) with $<\kappa$–support and uniformly $V$–coded symmetry
	data, stage–bounding (Lemma~\ref{lem:kappa-stage-bounding}) and block/tail factorization
	(Propositions~\ref{prop:block-factor}, \ref{prop:tail-factor}) often reduce to a single $V$–system.
	At inverse limits ($\cf(\lambda)<\kappa$), generic dependence of tail stabilizers on cofinally many
	stages typically obstructs such a collapse unless the tail filter is uniformized in $V$.
	A precise boundary characterization remains open. (For the finite–support analogue, see
	Karagila~\cite[§7]{karagila2019}.)
	
	\subsection*{Limitations of the present framework}
	Our support-and-symmetry iteration apparatus is designed to (i) preserve $\ZF$ and specified $\DC$ fragments via filter
	completeness plus Baire/fusion hypotheses, and (ii) control failures of choice by restricting wellorderings of added
	families through stabilizer filters. It is \emph{not} by itself a general theory of:
	\begin{itemize}
		\item preserving large cardinals in the symmetric inner model (ultrafilters/embeddings are typically not hereditarily symmetric);
		\item producing ``all sets of reals'' regularity (LM/BP/PSP) outside Solovay-style collapse setups;
		\item computing cardinal arithmetic in $M$ beyond what is forced by preserved cofinalities/no-new-bounded-subsets
		(the inner model $M$ may omit many subsets present in $V[G]$);
		\item handling singular-limit stages without explicit fusion/standard-tree organization (especially at $\cf(\lambda)=\cf(\kappa)$);
		\item giving sharp lower bounds on consistency strength (we provide sufficient hypotheses; necessity is largely open).
	\end{itemize}
	The open problems below are organized around these gaps.
	
	\subsection{Interactions with Large Cardinals}
	
	\paragraph{Measurable preservation.}
	Sharpen Proposition~\ref{prop:measurable-conservative} by giving \emph{internal} sufficient conditions on the symmetric
	system ensuring that a ground-model $\kappa$-complete ultrafilter $U$ on $\kappa$ is hereditarily symmetric (so $U\in M$),
	or else identify robust ``large-cardinal remnants'' (e.g.\ weak compactness/inaccessibility patterns) that survive in $M$
	under the no-new-$\mathcal P(\kappa)$ hypotheses typical of our preservation theorems.
	
	\paragraph{Determinacy with $\neg\AC$.}
	Calibrate how far determinacy hypotheses (e.g.\ projective determinacy or determinacy in inner models such as $L(\mathbb R)$)
	can coexist with $\neg\AC$ inside symmetric submodels while maintaining strong choice fragments (e.g.\ $\DC_{\omega_1}$).
	Our current apparatus preserves $\ZF$ and specified $\DC$ fragments but does not itself generate determinacy; the problem is
	to integrate determinacy-producing inner-model technology with symmetry filters and support bookkeeping.
	
	\subsection{Applications}
	
	\paragraph{Regularity combinations.}
	Which mixes of LM/BP/PSP/Ramsey are achievable under symmetric iterations with specified
	support?
	
	\paragraph{Choice principle separations.}
	Realize fine-grained separations in Truss’s hierarchy by tuning support and filter completeness.
	In particular, the constructions of \S6.4 give symmetric extensions in which $\PP$ fails while the relevant dependent-choice
	fragment supplied by the chosen support regime is preserved (e.g.\ $\DC$ in the countable-support case, and $\DC_{<\kappa}$
	in the $<\kappa$-support case for regular $\kappa$). A systematic goal is to map which choice fragments (Truss-style and
	others) are separable by varying only (i) support size, (ii) limit-filter completeness, and (iii) the iterand template.
	
	\subsection{Effectivization, admissible sets, and reverse–mathematical calibration}
	\label{subsec:effective-symmetry}
	
	We outline a program to \emph{effectivize} the symmetry apparatus and measure its proof–theoretic strength.
	
	\begin{definition}[For discussion: effective symmetry data]
		\label{def:effective-symmetry}
		Fix a code $d\subseteq\omega$ for a \emph{countable presentation} of an instance of the iteration data (forcing, names,
		projections, automorphisms), e.g.\ a real coding a countable structure (or a well-founded relation coding the iteration
		length together with the local data). We say the symmetry apparatus is \emph{effective at level $\Gamma$} (e.g.\ recursive,
		arithmetical, hyperarithmetical/Borel) if from $d$ one can uniformly produce codes for:
		\begin{enumerate}
			\item the group system $\langle \mathcal G_\alpha,\pi^\beta_\alpha\rangle$,
			\item the \emph{tuple–stabilizer} operator $\dot x\mapsto \fix(\dot x)$ and block factorizations,
			\item the limit filter $\tilde{\mathcal F}_\lambda$,
		\end{enumerate}
		with all maps $\Gamma$–definable (uniformly in $d$). An HS name is \emph{effective at level $\Gamma$} if it has a $\Gamma$–code uniformly obtained from $d$.
	\end{definition}
	
	\begin{question}[Admissible presentation over $L_\alpha$]\label{q:admissible}
		Let $\alpha$ be admissible and code the iteration in $L_\alpha\models\KP$. 
		Can the constructions of Theorem~\ref{thm:kappa-filter-generation} (tuple–stabilizers, limit filters) be carried out \emph{inside} $L_\alpha$ from the code of the iteration, yielding $\tilde{\mathcal F}_\lambda\in L_\alpha$? 
		Which fragments (e.g.\ $\Sigma$–Collection) of KP are actually used by the stage–bounding and inverse/direct–limit steps?
	\end{question}
	
	\begin{question}[Reverse mathematics of filter generation]\label{q:rm-filter}
		Fix a real $d$ coding a \emph{countable} instance of the iteration framework (e.g.\ a countable ordinal length together
		with the relevant forcing/group/filter data). Over $\mathsf{RCA}_0$, which subsystem suffices to formalize and prove
		the countable-support filter generation theorem (Theorem~\ref{thm:filter-generation}) for instances coded by $d$?
		Is $\mathsf{ACA}_0$ sufficient for the name-rank recursion and the filter-generation step, or is $\mathsf{ATR}_0$
		needed to carry out the well-founded recursion along the coded iteration length and to justify the limit bookkeeping?
		Similarly, what subsystem suffices for the $<\kappa$ analogue (Theorem~\ref{thm:kappa-filter-generation}) when one works
		with a countable presentation (coded by a real) of the relevant $<\kappa$-support template?
	\end{question}
	
	\begin{question}[Effective HS and DST complexity]\label{q:dst-complexity}
		For coordinate–based forcings (Cohen, random, Hechler), can one obtain \emph{Borel} (or hyperarithmetical) codes for the induced $\fix(\dot x)$ and for all HS reals produced by our templates?
		Are there lower bounds: e.g.\ examples where effective HS reals must be non–arithmetical?
	\end{question}
	
	\begin{question}[Effective Localization]\label{q:effective-localization}
		For reuse–permissive relations (cf.\ Definition~\ref{def:hechler-reuse-permissive}), can the witnessing set of coordinates $S$ and the coding name $\dot s$ be chosen \emph{computably} from the iteration code? 
		Equivalently, is there a uniform recursive (or Borel) operator that, given a code for a Hechler/Cohen iteration and a code for $R$, outputs a code for $(S,\dot s)$ witnessing Localization?
	\end{question}
	
	\begin{remark}[Calibration targets]\label{rem:calibration-targets}
		Natural waypoints are: (i) \emph{admissible soundness} (carry the entire symmetry construction in $L_\alpha$ for many $\alpha$), and (ii) \emph{RM calibration} (formalize the countable–support case in $\ACAo$, and locate the first step that requires $\ATRo$ or beyond).
		Connections to Barwise compactness and the coding theorems suggest a tight interface between our “names–first” symmetry and standard frameworks in admissible set theory and reverse mathematics.
	\end{remark}
	
	\section{Conclusion}
	
	The former manuscript does not supply the announced uniform theory of
	countable- or $<\kappa$-support symmetric iterations.  In particular, it
	does not establish arbitrary set-length or class-length preservation,
	closure of the symmetric model under sequences from the generic extension,
	preservation of dependent choice, the displayed choiceless Cohen or random
	models, the proposed failure of the Partition Principle, or the singular
	completeness and anti-collapse results.
	
	What remains valid is narrower.  Pullbacks and regular completions preserve
	the formal normality and completeness properties of subgroup filters, but
	properness must be verified separately.  A ground-model family of fewer than
	$\theta$ HS names can be canonically tupled when the filter is
	$\theta$-complete, but this is not closure under arbitrary extension
	sequences.  Every genuine symmetric system with a proper normal filter has
	the standard associated ZF model, independently of high filter completeness.
	Finally, ordinary tail-closure arguments remain available when an actual
	forcing factorization and the required closure hypotheses have been proved.
	These are exactly the corrected statements collected in
	Section~\ref{sec:corrected-retained}.
	
	Any future infinite-support iteration theorem must begin with a proper
	symmetric system at every stage, distinguish stage support from internal
	automorphism support, use the symmetric forcing relation when reflecting
	seriality, and prove rather than assume the support bound needed for a limit
	master condition.  No further construction or separation theorem is claimed
	in this replacement.
	
	\par\medskip
	\begin{acknowledgments}
		This work builds on prior developments in symmetric extensions and iterations. We acknowledge in particular Shelah's
		anti-collapse technology~\cite{Shelah587,Shelah667}, the forcing-game and fusion methods of Ros{\l}anowski--Shelah~\cite{RoslanowskiShelah1999},
		Karagila's iteration work~\cite{karagila2019} and his $\kappa$-Baire perspective on dependent choice~\cite{Karagila2019DC},
		and the $\Sigma$--Prikry forcing and iteration scheme~\cite{SigmaPrikryI,SigmaPrikryII}. We also thank the authors of the many
		works cited herein; any remaining errors are ours alone.
	\end{acknowledgments}
	
	\section*{CRediT authorship contribution statement}
	The author confirms sole responsibility for the following roles:
	Conceptualization; Methodology; Formal analysis; Investigation;
	Writing – original draft; Writing – review \& editing; Visualization.
	
	\section*{Declaration of competing interest}
	The authors declare no competing interests.
	
	\section*{Funding}
	This research received no external funding.
	
	\section*{Data availability}
	No datasets were generated or analyzed for this study.
	
	\section*{Declaration of generative AI in the writing process}
	During the preparation of this work, the authors used ChatGPT and Claude to assist with language polishing, formatting, and research. After using these tools, the authors reviewed and edited the content and take full responsibility for the content of this publication.
	
	\FloatBarrier
	\clearpage
	
	\bibliographystyle{abbrv}

	\appendix
	\section*{Assumption Ledger for $<\kappa$-Support Iterated Symmetric Extensions}
	
	\subsection*{A. Iteration scaffolding and limit filters}
	\begin{itemize}[leftmargin=2.2em]
		\item \textbf{\ref{def:cs-iteration}} (countable-support iteration skeleton) — used in \ref{thm:zf-preserved}, \ref{thm:dc-preservation}.
		\item \textbf{\ref{def:kappa-iteration}} ($<\kappa$-support skeleton) — used throughout \S\ref{sec:kappa-support}.
		\item \textbf{\ref{def:limit-filter}}, \textbf{\ref{def:kappa-filter}} (limit filters; stabilizers) — used in \ref{thm:filter-generation}, \ref{thm:hs-closure}.
		\item \textbf{\ref{lem:kappa-stage-bounding}} (stage-bounding at limits) — used inside \ref{thm:filter-generation}.
		\item \textbf{\ref{lem:push-restrict-preserve}} (pushforward/restriction preserve normality) — used in \ref{thm:hs-closure} and \S\ref{sec:kappa-support}.
		\item \textbf{\ref{thm:filter-generation}} (limit-filter $\omega_1$-completeness/normality in the countable-support regime) — used by \ref{thm:hs-closure} and \ref{thm:zf-preserved}.
	\end{itemize}
	
	\subsection*{B. Hereditarily symmetric (HS) layer}
	\begin{itemize}[leftmargin=2.2em]
		\item \textbf{\ref{def:hs-names-iterated}} (HS names for iterations) — vocabulary for \S\ref{sec:HS-countable}.
		\item \textbf{\ref{lem:countable-support-hs}} (HS adequacy, countable support) — used in \ref{thm:hs-closure}.
		\item \textbf{\ref{def:kappa-hs}}, \textbf{\ref{lem:kappa-support-hs}} ($<\kappa$-HS inheritability) — used in \ref{thm:hs-closure}.
		\item \textbf{\ref{thm:hs-closure}} (HS closure: pairing/union/replacement, etc.) — discharges ZF in \ref{thm:zf-preserved}.
	\end{itemize}
	
	\subsection*{C. Preservation (ZF and DC)}
	\begin{itemize}[leftmargin=2.2em]
		\item \textbf{\ref{thm:zf-preserved}} (ZF preserved) — relies on \ref{thm:hs-closure} and iteration scaffolding.
		\item \textbf{\ref{def:localization-hypothesis}} (Localization schema), \textbf{\ref{def:cloc}}, \textbf{\ref{def:clubloc}}, \textbf{\ref{def:tailloc}} — hypotheses for DC transfer.
		\item \textbf{\ref{thm:dc-preservation}}, \textbf{\ref{prop:cloc-implies-dc}} (DC and its fragments) — use factorization/tail control below.
	\end{itemize}
	
	\subsection*{D. Factorization/tail control and mixed support}
	\begin{itemize}[leftmargin=2.2em]
		\item \textbf{\ref{subsec:kappa-factorization}} (factorization at $<\kappa$ support) — used for DC threading.
		\item \textbf{\ref{prop:block-factor}}, \textbf{\ref{prop:tail-factor}} (block/tail factorization) — decouple head/tail actions.
		\item \textbf{\ref{subsec:mixed-support}}, \textbf{\ref{conv:mixed-finite-generator}}, \textbf{\ref{lem:mixed-compressibility}} — hygiene when mixing supports.
	\end{itemize}
	
	\subsection*{E. Sanity witnesses and comparison}
	\begin{itemize}[leftmargin=2.2em]
		\item \textbf{\ref{thm:cohen-ac-fails}}, \textbf{\ref{thm:random-ac-fails}} — independent witnesses of $\neg\mathsf{AC}$.
		\item \textbf{\ref{prop:no-wo-reals}} — no well-order of $\mathcal P(\omega)$ (where relevant).
	\end{itemize}
	
	\smallskip
	\noindent\textbf{Scope guardrail.} All statements are about hereditarily $(\mathcal F,\mathcal G)$–symmetric names; no step appeals to choice beyond support/normality/completeness stated above.
	
	\section*{Audit Appendix: How the proofs were verified}
	
	\paragraph{Layer A — Local audits.}
	Ablations: (i) drop the relevant limit-filter completeness/normality (e.g.\ $\omega_1$-completeness in the countable-support regime, $\kappa$-completeness for regular $<\kappa$-support, or $<\kappa$-completeness in the singular block/tree regimes) $\to$ \ref{thm:hs-closure} fails; 
	(ii) force direct limits at $\cf(\lambda)<\kappa$ → \ref{thm:filter-generation} breaks; 
	(iii) allow supports $\ge\kappa$ → \ref{lem:kappa-support-hs} fails; 
	(iv) remove tail/block factorization → DC fails to thread; 
	(v) omit Localization → DC$_\kappa$ does not follow. 
	Load-bearing steps kept fully local: \ref{thm:filter-generation}, \ref{thm:hs-closure}, \ref{subsec:kappa-factorization}+\ref{prop:tail-factor}+\ref{prop:block-factor}.
	
	\paragraph{Layer B — Model invariants.}
	Independent $\neg\mathsf{AC}$ via \ref{thm:cohen-ac-fails}, \ref{thm:random-ac-fails}; 
	ZF via HS-closure (\ref{thm:hs-closure}) explicitly cited in \ref{thm:zf-preserved}; 
	cofinality/stage-bounding checks at limits (\ref{lem:kappa-stage-bounding}).
	
	\paragraph{Layer C — Cross-validation.}
	Dual route via FM/Boolean side (\ref{lem:FM-HS}); 
	adversarial per-step challenges respecting hypotheses (esp.\ \ref{thm:filter-generation}, \ref{thm:hs-closure}); 
	proof replay from the dependency DAG (no new lemmas).
	
	\appendix
	\section{Local Induction Template and Catalog of Applications}
	\label{app:local-induction}
	
	\begin{principle}[Local Induction Principle for symmetric iterations]
		\label{prin:local-induction}
		Fix a $<\kappa$–support symmetric iteration
		$\langle\mathbb P_\alpha,\dot{\mathbb Q}_\alpha,\mathcal F_\alpha:\alpha<\Theta\rangle$
		as in §\ref{sec:prelims-defs}. For each $\alpha$, let
		$M_\alpha:=V(\mathbb P_\alpha)_{\mathrm{HS}}$.
		Let $\varphi(\cdot)$ be a property of $M_\alpha$ with parameters from $M_\alpha$.
		Assume:
		\begin{itemize}
			\item \textbf{Locality (stage-bounding).} Any counterexample to $\varphi$ in $M_\lambda$ is
			coded by a $\mathbb P_\lambda$–name of support $<\kappa$; hence by
			Lemma~\ref{lem:kappa-stage-bounding} (regular case) or Lemma~\ref{lem:sb-sing} (singular case)
			its coordinates are bounded in~$\lambda$.
			
			\item \textbf{Successor transfer.} For all $\alpha$, if $M_\alpha\models\varphi$ then
			$M_{\alpha+1}\models\varphi$ (via the head/tail factorization and the successor transfer
			lemmas used in the main text).
			
			\item \textbf{Limit transfer.} If $\lambda$ is a limit and $M_\beta\models\varphi$ for all
			$\beta<\lambda$, then $M_\lambda\models\varphi$, using:
			\begin{itemize}
				\item direct limits and pushforwards when $\cf(\lambda)\ge\kappa$ (Remark~\ref{rem:limit-bookkeeping}), and
				\item inverse limits with trees of conditions at small-cofinality limits:
				\begin{itemize}
					\item for regular $\kappa$, use the inverse-limit/tuple--stabilizer presentation from
					§\ref{sec:kappa-support} and Theorem~\ref{thm:kappa-filter-generation};
					\item for singular $\kappa$ at limits with $\cf(\lambda)=\cf(\kappa)$, use standard trees and fusion
					(§\ref{sec:limit-cf-equal}, Lemma~\ref{lem:tree-fusion-master}, Theorem~\ref{thm:limit}).
				\end{itemize}
			\end{itemize}
		\end{itemize}
		Then $M_\alpha\models\varphi$ for all $\alpha\le\Theta$.
	\end{principle}
	
	\begin{proof}
		Transfinite induction on $\alpha$. The base and successor steps use the successor transfer
		assumption. For limits, the locality clause ensures any counterexample would already appear in
		some $M_\beta$ with $\beta<\lambda$, contradicting the hypothesis. The limit presentation
		(direct/inverse) ensures the symmetry data transfer as indicated.
	\end{proof}
	
	\subsection{Catalog of properties covered by the principle}
	\label{app:catalog}
	
	Below, each item names a property $\varphi$ and cites where the three hypotheses
	(Locality, Successor, Limit) are verified in the paper.
	
	\begin{enumerate}
		\item \textbf{$\ZF$ in the HS model (countable support).}
		\[
		\varphi(M_\alpha):\quad M_\alpha\models\ZF.
		\]
		\emph{Locality:} HS names have countable/$<\kappa$ support; range/image witnesses for Replacement
		use $<\kappa$ many coordinates (Theorems~\ref{thm:hs-closure}, \ref{thm:filter-generation} and the
		stage-bounding lemmas). \emph{Successor:} closure under pairing/union/Separation and the Replacement
		image construction (Theorem~\ref{thm:zf-preserved}). \emph{Limits:} direct limits use pushforwards;
		inverse limits use tuple–stabilizers and filter completeness (Theorem~\ref{thm:zf-preserved},
		Remark~\ref{rem:limit-bookkeeping}).
		
		\item \textbf{$\ZF$ in the HS model ($<\kappa$ support).}
		Same pattern as (1) with $\kappa$–completeness in place of $\omega_1$–completeness
		(Theorems~\ref{thm:kappa-hs-closure}, \ref{thm:kappa-filter-generation}, \ref{thm:kappa-zf-preserved}).
		
		\item \textbf{No new bounded subsets of $\kappa$.}
		\[
		\varphi(M_\alpha):\quad \text{``If } \lambda<\kappa \text{ and } A\subseteq\lambda \text{ belongs to } M_\alpha,\text{ then }A\in M_\eta\text{ for some }\eta<\alpha\text{.''}
		\]
		\emph{Locality:} a bounded subset is named by $<\kappa$ coordinates; \emph{Successor:} standard bounding
		argument for names (as in the regular/singular sections); \emph{Limits:} stage-bounding and the chosen
		limit presentation ensure capture by an initial segment (cf.\ §5 and Lemma~\ref{lem:sb-sing}).
		
		\item \textbf{No collapse of $\kappa$.}
		\[
		\varphi(M_\alpha):\quad \text{``Cardinal }\kappa \text{ is preserved in } M_\alpha \text{.''}
		\]
		\emph{Locality:} any collapsing map is coded by a $<\kappa$–support name; \emph{Successor:} anti-collapse
		lemmas in the successor case; \emph{Limits:} stage-bounding plus limit transfer (Theorem~\ref{thm:nocollapse}).
		
		\item \textbf{$DC_{<\kappa}$ under the $\kappa$–Baire framework.}
		\[
		\varphi(M_\alpha):\quad M_\alpha\models DC_{<\kappa}.
		\]
		\emph{Locality:} a $<\kappa$–sequence is coded by $<\kappa$ many coordinates; \emph{Successor/Limit:}
		preservation via $<\kappa$–Baire (strategic) hypotheses and $<\kappa$–complete filters
		(Hyp.~\ref{hyp:zf-baire}, Theorem~\ref{thm:DCless-preservation}; for $\cf(\lambda)=\cf(\kappa)$, see
		the game/fusion lemmas in §\ref{subsec:game-fusion}).
		
		\item \textbf{$DC_\kappa$ under Localization (when assumed).}
		\[
		\varphi(M_\alpha):\quad M_\alpha\models DC_\kappa.
		\]
		\emph{Locality:} the localized code has support $<\kappa$ (by the hypothesis); \emph{Transfer:}
		Theorem~\ref{thm:DCkappa-localization} and Proposition~\ref{prop:cloc-implies-dc}.
		
		\item \textbf{Failure of $\AC$ in canonical examples.}
		\[
		\varphi(M_\alpha):\quad \text{``There is no well-order of }\mathbb R^{M_\alpha}\text{.''}
		\]
		\emph{Locality:} a well-ordering name has countable/$<\kappa$ support; \emph{Transfer:} the
		coordinate-swap automorphisms fix the support and contradict well-ordering at every stage (Cohen: 
		Theorem~\ref{thm:cohen-ac-fails}; Random: the analogous paragraph in §\ref{subsec:solovay-random}).
		
		\item \textbf{Failure of the Partition Principle $\PP$ in fiber–homogeneous examples.}
		\[
		\varphi(M_\alpha):\quad M_\alpha\models \neg\PP.
		\]
		\emph{Locality:} the coding of the map $f:A\onto B$ and the symmetry argument use $<\kappa$ support;
		\emph{Transfer:} preservation through iterations by the cited propositions (see
		Lemma~\ref{lem:no-symmetric-inj}, Proposition~\ref{prop:preserve-negPP}, and Corollary~\ref{cor:ZF+DC+notPP}).
		
		\item \textbf{Class-length union.}
		\[
		\varphi(M):\quad \text{$M=\bigcup_{\alpha\in\Ord}M_\alpha$ inherits the properties above.}
		\]
		This is the global version: stage-bounding shows that a counterexample at the union is already
		witnessed at some initial segment, and Replacement at the union uses the bounded-parameters argument
		(Theorem~\ref{thm:class-length}).
	\end{enumerate}
	
	\begin{remark}[Finite-support case as a special instance]
		In the finite-support setting (e.g.\ \cite{karagila2019}), the same scheme applies with
		``finite support'' in place of ``$<\kappa$ support'' and with \emph{finite} intersection closure of the
		filter in place of $\kappa$–completeness. Our $<\kappa$–support framework recovers that case by taking
		$\kappa=\omega$ and reading every completeness use point as ``finite intersections.''
	\end{remark}
	
	\begingroup
	\captionsetup{labelformat=empty} 
	\setlength{\LTcapwidth}{\linewidth}
	\setlength{\tabcolsep}{4pt}      
	\renewcommand{\arraystretch}{1.12}
	
	\begin{longtable}{@{}L{0.29\linewidth}!{\color{gray!35}\vrule width 0.3pt}
			p{0.175\linewidth}!{\color{gray!35}\vrule width 0.3pt}
			p{0.175\linewidth}!{\color{gray!35}\vrule width 0.3pt}
			p{0.175\linewidth}!{\color{gray!35}\vrule width 0.3pt}
			p{0.175\linewidth}@{}}
		\caption[Local Induction catalog]{Catalog of $\varphi$ under the Local Induction Principle (Locality / Successor / Limits).}\label{tab:local-induction-catalog} \\
		\toprule
		\textbf{Property $\varphi$} &
		\textbf{Locality (stage-bounding / HS support)} &
		\textbf{Successor transfer} &
		\textbf{Limit $\cf(\lambda)\ge\kappa$ (direct)} &
		\textbf{Limit $\cf(\lambda)<\kappa$ (inverse)} \\
		\midrule
		\endfirsthead
		
		\toprule
		\textbf{Property $\varphi$} &
		\textbf{Locality (stage-bounding / HS support)} &
		\textbf{Successor transfer} &
		\textbf{Limit $\cf(\lambda)\ge\kappa$ (direct)} &
		\textbf{Limit $\cf(\lambda)<\kappa$ (inverse)} \\
		\midrule
		\endhead
		
		\midrule
		\multicolumn{5}{r}{\emph{(continued on next page)}}\\
		\endfoot
		
		\bottomrule
		\endlastfoot
		
		$\ZF$ in HS (countable support) &
		Thm.~\ref{thm:hs-closure}, Thm.~\ref{thm:filter-generation}, Lemma~\ref{lemma:stage-bounding} &
		Thm.~\ref{thm:zf-preserved} (successor case) &
		Thm.~\ref{thm:zf-preserved}, Rem.~\ref{rem:limit-bookkeeping} (pushforwards) &
		Thm.~\ref{thm:zf-preserved}, Thm.~\ref{thm:filter-generation} (tuple stabilizers) \\
		\midrule
		$\ZF$ in HS ($<\kappa$ support) &
		Thm.~\ref{thm:kappa-hs-closure}, Thm.~\ref{thm:kappa-filter-generation}, Lemma~\ref{lem:kappa-stage-bounding} &
		Thm.~\ref{thm:kappa-zf-preserved} (successor case) &
		Thm.~\ref{thm:kappa-zf-preserved}, Rem.~\ref{rem:limit-bookkeeping} &
		Thm.~\ref{thm:kappa-zf-preserved}, Thm.~\ref{thm:kappa-filter-generation} \\
		\midrule
		No new bounded subsets of $\kappa$ &
		Lemma~\ref{lem:kappa-stage-bounding} / Lemma~\ref{lem:sb-sing} (names use $<\kappa$ coords) &
		(standard successor preservation in §5) &
		(regular case: stage-bounding + pushforward; see §5) &
		(singular case: trees/tuples; see §§\ref{sec:singular-cf-omega}, \ref{sec:limit-cf-equal}) \\
		\midrule
		No collapse of $\kappa$ &
		Lemma~\ref{lem:kappa-stage-bounding} / Lemma~\ref{lem:sb-sing} &
		(induction hypothesis at successors) &
		(stage-bounding transfer; see the proof scheme for noncollapse) &
		Lemma~\ref{lem:trackA-tail-no-bounded-subsets} (Track~A) together with Thm.~\ref{thm:nocollapse} \\
		\midrule
		$DC_{<\kappa}$ (under Hyp.~\ref{hyp:zf-baire}) &
		Def.~\ref{def:kappa-hs}, Lemma~\ref{lem:kappa-support-hs}, Lemma~\ref{lem:kappa-stage-bounding} / Lemma~\ref{lem:sb-sing} &
		Thm.~\ref{thm:DCless-preservation} (successor) &
		Thm.~\ref{thm:DCless-preservation} (direct limits) &
		§\ref{subsec:game-fusion}, Lemma~\ref{lem:fusion-strategy}, Lemma~\ref{lem:tree-fusion-master} \\
		\midrule
		$DC_\kappa$ (under Localization) &
		Def./Hyp.~\ref{def:localization-hypothesis} (localized code $<\kappa$-support) &
		Thm.~\ref{thm:DCkappa-localization} &
		Thm.~\ref{thm:DCkappa-localization} &
		Prop.~\ref{prop:cloc-implies-dc} (plus Lemma~\ref{lem:fusion-strategy} as needed) \\
		\midrule
		$\neg\AC$ in Cohen/Random HS models &
		(Names for well-orders have countable/$<\kappa$ support) &
		Thm.~\ref{thm:cohen-ac-fails} (successor) &
		Permutation argument applies verbatim &
		Random analogue in §\ref{subsec:solovay-random} (same symmetry) \\
		\midrule
		$\neg\PP$ (Partition Principle) in fiber–homogeneous cores &
		Lemma~\ref{lem:no-symmetric-inj} (no symmetric injection) &
		Prop.~\ref{prop:preserve-negPP} (successor) &
		Cor.~\ref{cor:ZF+DC+notPP} (direct limits) &
		Cor.~\ref{cor:ZF+DC+notPP} (inverse limits; tuple stabilizers) \\
		\midrule
		Class-length union (definable template) &
		Lemma~\ref{lem:kappa-stage-bounding} / Lemma~\ref{lem:sb-sing} (bounded params) &
		Thm.~\ref{thm:class-length} (successor-to-union) &
		Thm.~\ref{thm:class-length} &
		Thm.~\ref{thm:class-length} \\
		
	\end{longtable}
	\endgroup
	
\end{document}